\documentclass[preprint, dvips]{article}
\usepackage{amsfonts, amsmath, amssymb, amsthm}
\usepackage[dvips]{graphicx}

\newcommand{\I}{\mathop{\mathrm{I}}}

\newcommand{\fs}{\mathop{\mathrm{\mathfrak{fs}}}}

\newcommand{\FS}{\mathop{\mathrm{\mathcal{FS}}}}

\newcommand{\f}{\mathop{\mathrm{\mathfrak{f}}}}
\newcommand{\p}{\mathop{\mathrm{\mathfrak{p}}}}

\newcommand{\sech}{\mathop{\mathrm{sech}}}
\newcommand{\arcsinh}{\mathop{\mathrm{arcsinh}}}
\newcommand{\arctg}{\mathop{\mathrm{arctg}}}

\newcommand{\csch}{\mathop{\mathrm{csch}}}

\newcommand{\Res}{\mathop{\mathrm{Res}}}

\newtheorem{thm}{Theorem}[section]
\newtheorem{rem}{Remark}[section]
\newtheorem{prop}{Proposition}[section]
\def\D{\Delta} \def\l{\lambda}
\def\s{\sigma} \def\th{\theta}
\def\d{\delta} \def\r{\rho}
\def\ri{\right} \def\le{\left}
\def\no{\nonumber}
\def\Eq{\Leftrightarrow}
\def\R{\mathbb R}
\def\fr{\frac}

\begin{document}

\title{\textbf{On spectral properties and statistical analysis of Fisher-Snedecor diffusion}}

\maketitle

\centerline{\author{\textbf{F. Avram}$^{\textrm{a}}$, \textbf{N. N. Leonenko}$^{\textrm{b},}$\footnote{Corresponding author.} \textbf{and N. \v{S}uvak}$^{\textrm{c}}$}}

{\footnotesize{
$$\begin{tabular}{l}
  $^{\textrm{a}}$ \emph{Department of Mathematics, University of Pau, 64 000 Pau, France, e-mail: florin.avram@univ-pau.fr} \\
  $^{\textrm{b}}$ \emph{School of Mathematics, Cardiff University, Senghennydd Road, Cardiff CF244AG, UK, e-mail: LeonenkoN@Cardiff.ac.uk} \\
  $^{\textrm{c}}$ \emph{Department of Mathematics, University of Osijek, Gajev Trg 6, HR-31 000 Osijek, Croatia, e-mail: nsuvak@mathos.hr}
\end{tabular}$$}}

\bigskip

\noindent\textbf{Abstract.}
We consider the problem of parameter estimation for an ergodic diffusion with Fisher-Snedecor invariant distribution, to be called
Fisher-Snedecor diffusion. We compute the spectral representation of its  transition density, which involves a finite number of discrete eigenfunctions (Fisher-Snedecor polynomials) as well as a continuous part. We propose moments based estimators (related to the Fisher-Snedecor polynomials) and  prove their consistency and asymptotic normality. Furthermore, we propose a statistical test for the distributional assumptions on the marginal distribution of the Fisher-Snedecor diffusion, based on the moment condition derived from the corresponding Stein's equation.

\bigskip

\noindent\textbf{Key words.} Diffusion process, Method of moments, Pearson equation, Fisher-Snedecor polynomials, Stein equation, Transition density.

\bigskip

\noindent\textbf{Mathematics Subject Classification (2010):} 33C47, 60G10, 60J25, 60J60, 62M10, 62M05, 62M02.

\bigskip

\section{Introduction}

{\bf Historical roots.}
The study of diffusion processes with invariant distributions from the Pearson family started in the 1930's, when Kolmogorov \cite{Kolmogorov, Shiryayev} studied the Fokker-Planck or forward Kolmogorov equation
\begin{equation*}
\frac{\partial p}{\partial t} = -\frac{\partial}{\partial x} \left( \left(a_{1}x + a_{0} \right) p \right) + \displaystyle\frac{\partial^{2}}{\partial x^{2}} \left( \left(b_{2}x^{2} + b_{1}x + b_{0} \right) p \right), \quad p = p(x, t), \quad x \in \mathbb{R}, \quad t \geq 0,
\end{equation*}
with a linear drift and a quadratic squared diffusion, and observed that the invariant density $\p(\cdot)$ satisfies the differential
equation
\begin{equation}
\frac{\p^{\prime}(x)}{\p(x)} = \frac{(a_{1} - 2b_{2})x +(a_{0} - b_{1})}{b_{2}x^{2} + b_{1}x + b_{0}} = \frac{c_{1}x + c_{0}}{b_{2}x^{2} + b_{1}x + b_{0}} = \frac{q(x)}{s(x)}, \quad x \in \mathbb{R}, \label{PEq}
\end{equation}
introduced by K. Pearson \cite{Pearson} in 1914 (and of illustrious history - see Diaconis and Zabel \cite{DZ}), in order to  unify some of the most important statistical distributions.

It seems appropriate to call this important class of processes Kolmogorov-Pearson (KP) diffusions, or Gauss-hypergeometric diffusions, due to the appearance of the Gauss $_2 F_1$ function (and its limiting confluent forms) in various explicit formulas.

For a long period of time after that, KP diffusions were neglected, with some notable exceptions like Wong (1964) \cite{Wong}, who reemphasized the importance of this class of models as a most natural extension of the "first order statistical description characterized by $\p(x)$" to a time dependent model, and computed spectral representations of the transition density in some cases (but missed the Fisher-Snedecor case). \\

{\bf  Mathematical finance motivations.}
Recently, the interest in these processes was reawakened in the context of financial modeling. The most famous case is the Merton-Black-Scholes  SDE
\begin{eqnarray}
d X_t &= & r X_t \; d t  + \s  X_t \; d W_t  \Eq \; X_t = X_0
e^{(r- \s^2/2) t + \s W_t}.
\end{eqnarray}
This has had a huge impact in mathematical finance due to its tractability, which lead to a large variety of explicit formulas for the transition and first passage probabilities necessary for the pricing of options. Other tractable diffusions are the Ornstein-Uhlenbeck (Vasicek), and Cox-Ingersoll-Ross (CIR) models.

In the search for more flexible models, financial mathematicians have rediscovered recently the convenience of the Kolmogorov-Pearson diffusions (and others, like diffusions with constant elasticity of variance). Some notable contributions are due to Albanese, Campolieti, Carr and Lipton \cite{ACCL}, Linetsky \cite{Linetsky1, Linetsky2, Linetsky3}, Kuznetsov \cite{AK,Kuz}, Mendoza, Carr and Linetsky \cite{Mendoza} and Shaw \cite{Shaw}. In parallel, the statistical analysis of these processes was developed by K\"uchler and S\o rensen and Forman and S\o rensen \cite{Forman}. Recently, the interest in developing tractable processes resulted in the introduction  by Cuchiero, Keller-Ressel and Teichmann \cite{Teich} of the unifying family of polynomial Markovian processes, whose generator maps polynomials into polynomials of (at most) the same degree. \\

{\bf  The Student parametrization of Kolmogorov-Pearson (KP) diffusions.}
Following the modern formulations of \cite{Shaw} and \cite{Forman}, consider the SDE
\begin{eqnarray}\label{Shaw}&&
dX_{t} =\th (\mu - X_{t} ) \, dt + \s_1 \, dW^{(1)}_{t} + \s_2 \, X_{t} dW^{(2)}_{t}, \quad t \geq 0, \quad \Rightarrow \\&& dX_{t} =\theta \le( \mu -  X_{t} \ri) \, dt + \sqrt{ k \le(\le( X_{t}+ \r \frac{\s_1}{\s_2}\ri)^{2}+ (1-\r^2) \le(\frac{\s_1}{\s_2}\ri)^2\ri)} \,dW_{t},  
\no
\end{eqnarray}
where $
k=\s_2^2 ,$ \; 
$W^{(1)}_{t}$ and $W^{(2)}_{t}$, $t \geq 0$, are standard Brownian motions with correlation $\r,$ and $\{W_{t}, \, t \geq 0\}$ is a standard Brownian motion resulting from combining the two.

{\bf Notes:}
1) The second formulation, to be  called the {\bf skew Student parametrization}, may be used for all KP diffusions (by restricting if necessary to the range of values of $x$ where the square root makes sense). In the  first formulation, one implicitly assumes  $|\r| \leq 1, k  \geq 0$, which characterize the Student subclass of KP diffusions --see below.

2) The linearity of the drift and the quadratic variance  ensure the existence of polynomial eigenfunctions. \\

{\bf The S\o rensen/skew Student parametrization.}
Putting $\mu' = \r \frac{\s_1}{\s_2}, k=\frac{\th}{\nu -1} 
$, 
we arrive at
\begin{equation}
dX_{t} = \theta(\mu- X_{t}) \, dt + \sqrt{\frac{\theta}{ \nu-1} \left[\left({X_{t} - {\mu' }}{}\right)^{2}+\delta^2 \right]} \, dB_{t}, \quad t \geq 0,
\end{equation}
where $\d^2=(1-\r^2) (\frac{\s_1}{\s_2})^2$. This parametrization makes sense for the whole KP family (by allowing $\d^2 \leq 0$), but it will only produce diffusions living on $(-\infty ,\infty)$ when $\delta \in \R, k > 0 .$ \\

{\bf A classification of Kolmogorov-Pearson diffusions} in six basic subfamilies may be achieved by using the criteria based on the degree $\deg(s)$ of the polynomial $s(x)$ from the denominator of the Pearson equation \eqref{PEq} and on the sign of the leading coefficient $k $ and the discriminant $\Delta(s)$ in the quadratic case. The classification is given in the following table

{\footnotesize{
$$\begin{tabular}{|c||c|c|c|c|} \hline {\sf{Pearson diffusion}} & {\sf{Characteristic property}} & {\sf{Invariant density}} & {\sf{Parameter space}} \\ \hline \hline

Ornstein-Uhlenbeck & $\deg(s) = 0$ & {\tiny{$\displaystyle\frac{1}{\sigma \sqrt{2 \pi}} \, e^{-\frac{(x-\mu)^{2}}{2 \sigma^{2}}}$}} & $\mu \in \mathbb{R}$, $\sigma > 0$ \\ \hline

Gamma or CIR & $\deg(s) = 1$ & {\tiny{$\displaystyle\frac{\alpha^{\beta}}{\Gamma(\beta)} \, x^{\beta - 1} e^{-\alpha x} \, {\I}_{\langle 0, \infty \rangle}(x)$}} & $\alpha > 0$, $\beta > 0$ \\ \hline

Beta or Jacobi & $k < 0$, $\deg(s) = 2$, $\Delta(s) > 0$ & {\tiny{$\displaystyle\frac{1}{B\left( \alpha, \beta \right)} \, x^{\alpha - 1} (1-x)^{\beta - 1} \, {\I}_{[0, 1]}(x)$}} & $\alpha > 0$, $\beta > 0$ \\ \hline

Student (see \eqref{StNormConst}) & $k > 0$, $\deg(s) = 2$, $\Delta(s) < 0$ & {\tiny{$c(\mu, \mu', \alpha, \delta) \, \displaystyle\frac{\exp{\left(\frac{\mu - \mu'}{a \delta} \arctg{\left( \frac{x - \mu'}{\delta} \right)}\right)}}{\left[1 + \left( \frac{x - \mu'}{\delta} \right)^{2} \right]^{\frac{1}{2 a} + 1}}$}} & {\large{{$\mu \in \mathbb{R}, \, \, \mu' \in \mathbb{R}, \atop a > 0, \, \, \delta > 0$}}} \\ \hline

Reciprocal gamma & $k > 0$, $\deg(s) = 2$, $\Delta(s) = 0$ & {\tiny{$\displaystyle\frac{\alpha^{\beta}}{\Gamma(\beta)} \, x^{-\beta - 1} e^{-\frac{\alpha}{x}} \, {\I}_{\langle 0, \infty \rangle}(x)$}} & $\alpha > 0$, $\beta > 1$ \\ \hline

Fisher-Snedecor & $k > 0$, $\deg(s) = 2$, $\Delta(s) >  0$ & {\tiny{$\displaystyle\frac{\alpha^{\frac{\alpha}{2}} \beta^{\frac{\beta}{2}}}{B\left( \frac{\alpha}{2}, \frac{\beta}{2} \right)} \, x^{\frac{\alpha}{2} - 1} (\alpha x + \beta)^{-\frac{\alpha + \beta}{2}} \, {\I}_{\langle 0, \infty
\rangle}(x)$}} & $\alpha > 2$, $\beta > 2$ \\ \hline
\end{tabular}$$}}
\begin{center}
{\sf{\footnotesize{{\scshape{Table 1.}} \quad \emph{Classification of ergodic stationary Pearson diffusions.}}}}
\end{center}

{\bf Example:}
In the important Student case, the scale and invariant speed densities are:
\begin{equation*}
\mathfrak{s}(x) = (\d(\tilde{x}^{2}+1))^{\frac{1}{2a}} \; e^{ - \frac{\mu -\mu' }{a \delta} \arctg{(\tilde{x})}}, \quad \mathfrak{m}(x) = \frac{ e^{\frac{\mu -\mu'
}{a \delta} \; \arctg{(\tilde{x})}}} {(\tilde{x}^{2}+1)^{\frac{1}{2a}+1}}, x \in \mathbb{R},
\end{equation*}
where we put $\tilde{x} := (x -\mu')/{\delta}, \; 2 a=\fr{k}{ \th}= \frac{1}{\nu -1}.$

When $ \nu > 1 \Eq a >0,$ the KP diffusions will have Student stationary distribution with $\nu$ degrees of freedom, since the speed density may be normalized, arriving thus to diffusions with the invariant density:
\begin{equation}
f(x) = c(\mu, \mu^{\prime}, a, \delta) \frac{\exp \left\{ \frac{\mu -\mu^{\prime }}{a \delta}\arctg{\left( \frac{x - \mu^{\prime}}{\delta } \right)} \right\}}{\left[ 1 + \left( \frac{x - \mu^{\prime}}{\delta}\right)^{2} \right]^{\frac{1}{2a} + 1}}, \quad x \in \mathbb{R}, \label{STPDF}
\end{equation}
where (see Avram et al. \cite{Avram2})
\begin{equation}
c(\mu, \mu^{\prime}, a, \delta) = \frac{\Gamma \left( 1 + \frac{1}{2a}\right)}{\delta \sqrt{\pi} \Gamma \left( \frac{1}{2} + \frac{1}{2a} \right)} \, \prod\limits_{k=0}^{\infty}\left[ 1 + \left( \frac{\frac{\mu - \mu^{\prime}}{2a \delta}}{1 + \frac{1}{2a} + k}\right)^{2} \right]^{-1}. \label{StNormConst}
\end{equation}


The Ornstein-Uhlenbeck process, CIR process and Jacobi diffusion are well studied and widely applied. However, the first results on statistical analysis of reciprocal gamma, Student and Fisher-Snedecor diffusion, which all have heavy-tailed invariant distributions, are quite new (see \cite{Forman, Leonenko1, Leonenko2}). Their study involves the analysis of the spectrum of the corresponding infinitesimal generators. Namely, the spectrum of the infinitesimal generator of Ornstein-Uhlenbeck process, CIR process and Jacobi diffusion is simple and purely discrete with
classical orthogonal polynomials as corresponding eigenfunctions: Hermite, Laguerre and Jacobi polynomials, respectively. In the case of Pearson diffusions with heavy-tailed invariant distributions, the spectrum of the infinitesimal generator consist of two disjoint parts: the discrete spectrum (consisting of finitely many simple eigenvalues) and the essential spectrum. Furthermore, in all these cases corresponding eigenfunctions are less known finite systems of orthogonal polynomials: Bessel polynomials for reciprocal gamma diffusion (see \cite{Leonenko1}), Routh-Romanovski polynomials for Student diffusion (see \cite{Leonenko2}) and polynomials related to Fisher-Snedecor invariant distribution which have no common name (we will refer to them as to Fisher-Snedecor polynomials).

In this paper we focus on the statistical analysis of the ergodic stationary Fisher-Snedecor diffusion, by which we mean solutions of the non-linear stochastic differential equation
\begin{equation}
dX_{t} = -\theta \left(X_{t} - \mu \right) \, dt + \sqrt{2 \theta X_{t} \left( \frac{ X_{t}}{\beta/2 - 1} +\frac{ \mu}{\alpha/2 } \right)} \, dW_{t}, \quad t \geq 0, \label{fssde1}
\end{equation}
with speed measure/invariant distribution proportional to:
\begin{equation}\frac{(\alpha x)^{\frac{\alpha}{2} - 1}}{(\alpha x + \mu(\beta-2))^{\frac{\alpha}{2} + \frac{\beta}{2}}} \,  {\I}_{\langle 0, \infty \rangle}(x)\end{equation}
By taking $\mu$ as the mean of the Fisher-Snedecor distribution  $\mu=\frac{\beta}{\beta-2},$  (see Section \ref{FSd}) which may be always achieved by scaling, and by assuming $\alpha > 0$ and $\beta > 2$  (ensuring thus the existence of the first moment), we arrive to the classic Fisher-Snedecor distribution $\mathcal{FS}(\alpha, \beta)$ with $(\alpha,\beta)$ degrees of freedom as invariant distribution.

Statistical analysis of stochastic processes often demands the acquaintance of the transition density. Since for the most of diffusion models transition density is unavailable in the closed form, the statistical analysis is often very complicated. A possible solution to this problem is estimation of the transition density, which is used in a recent paper by Beskos \cite{Beskos} in the likelihood based framework. A\"it-Sahalia and Mykland \cite{AitSahalia2, AitSahalia3}, in the approach based on the moment conditions which satisfy certain differentiability conditions, used approximation of the asymptotic variance by the Taylor expansion with terms known in the closed form. However, since for the Fisher-Snedecor diffusion the structure of the spectrum of the infinitesimal generator can be fully determined, the use of the spectral representation of the transition density is implied. Therefore we present in Section 4 a closed expression for the spectral representation of the transition density of Fisher-Snedecor diffusion given in form of the finite sum related to discrete spectrum of the infinitesimal generator (i.e. eigenvalues and eigenfunctions - Fisher-Snedecor polynomials) and the integral which is taken over the absolutely continuous spectrum of the infinitesimal generator. The complicated form of this expression has significant impact on methods used in the statistical part of this paper. However, orthogonality of Fisher-Snedecor polynomials and hypergeometric functions appearing in the continuous part of the spectral representation makes this result applicable in the statistical analysis (see Remark \ref{rem:ortho}). However, in statistical applications treated here in Sections 5 and 6, the discrete part of the spectral representation is used much more frequently than the continuous part. This suggests the importance of Fisher-Snedecor polynomials which have been ignored in standard mathematical books, for instance in Abramowitz and Steghun \cite{Abramowitz}, Chichara \cite{Chichara}, Erdely \cite{Erdely}, Nikiforov and Uvarov \cite{Nikiforov} and Szeg\"o \cite{Szego}. Significant progress in their research is recently made by Masjed-Jamei \cite{Masjed-Jamei1, Masjed-Jamei2} and by Koepf and Masjed-Jamei \cite{Koepf}. Note that since essentially only the
discrete part of the spectrum is used here, readers motivated mainly by statistical applications may skip Sections 4.2. - 4.5.

In statistical part of the paper we observe the problems of parameter estimation and testing hypothesis about invariant distribution of the ergodic Fisher-Snedecor diffusion. Note that statistical inference for ergodic diffusions is widely studied (see the recent book by Kutoyants \cite{Kutoyants2}, the paper by Kutoyants and Yoshida \cite{Kutoyants1} and references therein). The recent review on diffusion estimation can be found in a survey paper by H. S\o rensen \cite{SorensenH}. Significant progress in statistical analysis of Pearson diffusions was recently made by Forman and S\o rensen \cite{Forman}, where the problem of parameter estimation is treated by method based on the martingale estimation equations. It is proved that the resulting estimators are consistent and asymptotically normal with given explicit expressions for the elements of the limiting covariance matrix. However, martingale estimation equation does not always provide estimators in the explicit form and calculation of the covariance matrix can be quite complicated. Furthermore, in our case the maximum likelihood method and the Bayesian estimation theory are not available since the diffusion coefficient depends on unknown parameters. Also, the minimum contrast estimation procedure based on Whittle functionals (see Anh et al. \cite{Anh1} and Leonenko and Sakhno \cite{Sakhno}) is not applicable here, since the fourth order spectral density could not be explicitly calculated. Therefore, parameters of the Fisher-Snedecor diffusion are estimated by the simple method of moments. As shown in Kutoyants and Yoshida \cite{Kutoyants1}, asymptotic efficiency of the method of moments can be quite good. Another reason that makes method of moments suitable for our problem are explicit forms of resulting estimators. Their consistency follows from the ergodic property of the observed diffusion. Furthermore, their asymptotic normality is implied by the functional central limit theorem based on the strong mixing property with an exponentially decaying rate (see Genon-Catalot et al.
\cite{Genon-Catalot}) and the functional delta method (see Serfling \cite{Serfling}). Explicit form of the limiting covariance matrix of the multivariate estimator of parameters of invariant distribution is calculated according to the new method based on the closed form expression for the spectral representation of the transition density and the finite system of orthogonal Fisher-Snedecor polynomials. This method basically enables calculation of moments of the form $E[X_{s + t}^{m} X_{s}^{n}]$, where $m$ and $n$ are at most equal to the number of orthogonal polynomials. This makes an important problem of constructing asymptotic confidence intervals for unknown parameters operational.

Statistical test for Fisher-Snedecor distributional assumptions constructed in this paper is based on the generalized method of moments. This procedure is implied by the Stein equation (see Schoutens \cite{Schoutens}) for the Fisher-Snedecor diffusion and on the results about the Pearson family of distributions of Bontemps and Meddahi \cite{Bontemps2}. Beside these general results, Bontemps and Meddahi proved the robustness of Hermite polynomials for testing normality by using the Stein equation for the normal distribution (see \cite{Bontemps1}). However, we verified that Fisher-Snedecor polynomials are not robust test functions against the parameter uncertainty in our model (see Subsection 6.3.). \\

\textbf{Contents and main results.}
We focus here on statistical analysis of the ergodic diffusion process with Fisher-Snedecor invariant distributions. The paper is structured as follows. In Section 2 we present well known information on Fisher-Snedecor distribution and its moments. In Section 3 we define the Fisher-Snedecor diffusion and present its most important properties. Spectral representation of the transition density of Fisher-Snedecor diffusion is developed in Section 4. In Section 5 we present method of moments estimators for parameters of this process and establish their asymptotic properties: consistency in Theorem \ref{Asymptotic}.(i) and asymptotic normality in Theorem \ref{Asymptotic}.(ii)-(iv). It is important to point out that orthonormality of Fisher-Snedecor polynomials and the closed form expression for spectral representation of transition density made possible calculation of explicit form of the limiting covariance matrix of the bivariate estimator of parameters of invariant distribution (see Theorem \ref{Asymptotic}.(ii)). Section 6 develops methodology for testing statistical hypothesis about Fisher-Snedecor invariant distribution. The established method is based on the Stein equation for Fisher-Snedecor diffusion and the moment condition closely related to the Fisher-Snedecor polynomials. In Appendix A we present properties of the finite orthogonal system of Fisher-Snedecor polynomials, Appendix B contains well known facts about Gauss hypergeometric functions, while in Appendix C we present classification of boundaries of the state space of Fisher-Snedecor diffusion.

\medskip

\section{General information about Fisher-Snedecor distribution \label{FSd}}

Random variable $X$ has Fisher-Snedecor distribution with $\alpha > 0$ and $\beta > 0$ degrees of freedom, i.e. $X \sim \mathop{\mathrm{\mathcal{FS}}}(\alpha, \beta)$, if its probability density function is given by
\begin{equation}
\fs(x) = \frac{ \beta^{\frac{\beta}{2}}}{B\left( \frac{\alpha}{2}, \frac{\beta}{2} \right)} \, \frac{(\alpha x)^{\frac{\alpha}{2} - 1}}{(\alpha x + \beta)^{\frac{\alpha}{2} + \frac{\beta}{2}}} \, \alpha \, {\I}_{\langle 0, \infty \rangle}(x) =\frac{(\frac{\alpha x}{\alpha x+\beta})^{\frac{\alpha}{2}}(\frac{\beta}{\alpha x+ \beta})^{\frac{\beta}{2}}} {x B\left(\frac{\alpha}{2}, \frac{\beta}{2} \right)} \, {\I}_{\langle 0,\infty \rangle}(x),  \label{fs}
\end{equation}
where $B(\cdot, \cdot)$ is the standard Beta function.

\begin{rem}
\emph{While a more general density
$$\frac{ w^{\frac{\beta}{2}}}{B\left( \frac{\alpha}{2}, \frac{\beta}{2} \right)} \, \frac{(\alpha x)^{\frac{\alpha}{2} - 1}}{(\alpha x + w)^{\frac{\alpha}{2} + \frac{\beta}{2}}} \, \alpha \, {\I}_{\langle 0, \infty \rangle}(x) =\frac{(\frac{\alpha x}{\alpha x+ w})^{\frac{\alpha}{2} } (\frac{w}{\alpha x+ w})^{\frac{\beta}{2}}} {x B\left( \frac{\alpha}{2}, \frac{\beta}{2} \right)} \, {\I}_{\langle 0, \infty \rangle}(x), \quad w > 0,$$
might be of interest as well, in this paper we study the Fisher-Snedecor distribution with $w = \beta$.}
\end{rem}

Our distribution belongs to the Pearson family of continuous distributions and is also known as Pearson type VI distribution (see Pearson \cite{Pearson}). In particular, it follows that the tail of the Fisher-Snedecor distribution with density \eqref{fs} decrease like $x^{-(1+\beta/2)}$, and so this distribution is heavy-tailed.

Moment of the $n$-th order of $\mathop{\mathrm{\mathcal{FS}}}(\alpha, \beta)$ distribution is given by the
expression
\begin{equation*}
E\left[ X^{n} \right] = \frac{\beta^{n}}{\alpha^{n-1}} \, \frac{%
\prod_{k=1}^{n-1}(\alpha+2k)}{\prod_{k=1}^{n}(\beta-2k)} = \left(\frac{\beta%
}{\alpha} \right)^{n} \, \frac{\Gamma\left( \frac{\alpha}{2} + n
\right) \, \Gamma\left( \frac{\beta}{2} - n \right)}{\Gamma\left(
\frac{\alpha}{2} \right) \, \Gamma\left( \frac{\beta}{2} \right)},
\quad \beta > 2n, \quad n \in \mathbb{N}.
\end{equation*}

Furthermore, moments could be calculated according to the recurrence relation
\begin{equation*}
\left( \frac{2(k+2)}{\beta+2} - 1 \right) \, E\left[ X^{k+1} \right] = -%
\frac{\beta (\alpha - 2) + 2 \beta (k + 1)}{\alpha (\beta + 2)} \,
E\left[ X^{k} \right], \quad \beta > 2(k+1), \quad k \in \mathbb{N}.
\end{equation*}

In particular, expectation and variance are:
\begin{equation*}
E[X] = \frac{\beta}{\beta - 2}, \quad \beta > 2, \quad \quad %
\mathop{\mathrm{Var}}(X) = \displaystyle\frac{2\beta^{2}(\alpha+\beta-2)}{%
\alpha(\beta - 2)^{2}(\beta-4)}, \quad \beta > 4.
\end{equation*}

\begin{rem}
\emph{Fisher-Snedecor distribution with positive integer degrees of
freedom is frequently used in statistics (e.g. in analysis of
variance). In particular, if $\chi_{n}^{2}$ and $\chi_{m}^{2}$ are
independent chi-square random variables with $n \in \mathbb{N}$ and
$m \in \mathbb{N}$ degrees of freedom, respectively, then the random
variable $\left[ (\chi_{n}^{2}/n) / (\chi_{m}^{2}/m) \right]$ has
Fisher-Snedecor distribution $\mathcal{FS}(n, m)$.}
\end{rem}

\medskip

\section{Fisher-Snedecor diffusion process}

Fisher-Snedecor diffusion is a solution of the non-linear stochastic differential equation
\begin{equation}
dX_{t} = -\theta \left(X_{t} - \frac{\beta}{\beta - 2} \right) \, dt + \sqrt{\frac{4 \theta}{\alpha(\beta - 2)} X_{t} (\alpha X_{t} + \beta) } \, dW_{t}, \quad t \geq 0, \label{fssde}
\end{equation}
(which is equation \eqref{fssde1} for $\mu = \beta/(\beta - 2)$). The infinitesimal parameters, i.e. the drift parameter $\mu(x)$ and the
diffusion parameter $\sigma(x),$ of the Fisher-Snedecor diffusion \eqref{fssde} are respectively given by
\begin{equation}
\mu(x) = -\theta \left(x - \frac{\beta}{\beta-2} \right), \quad \quad \sigma(x) = \sqrt{\displaystyle\frac{4 \theta}{\alpha(\beta - 2)} \, x(\alpha x + \beta)}, \quad \quad \beta > 2, \label{dd}
\end{equation}
where the restriction $\beta > 2$ ensures simultaneously that $\mu = \beta/(\beta - 2)$ is positive and that $\th > 0$, so that $\mu$ is the stationary mean.

For $x > 0$ the corresponding scale density is
\begin{equation}
\mathfrak{s}(x) =\exp\left( -\int_{.}^x \frac{2 \mu(u)}{\sigma^{2}(u)} \, du\right)= x^{-\frac{\alpha}{2}} \, (\alpha x + \beta)^{\frac{\alpha}{2} + \frac{\beta}{2} -1}, \label{scale}
\end{equation}
while the speed density is
\begin{equation}
\mathfrak{m}(x) = \frac{2}{\sigma^{2}(x) \mathfrak{s}(x)} = \frac{\alpha (\beta - 2)}{2 \theta} \, x^{\frac{\alpha}{2}-1} \, (\alpha x + \beta)^{-\frac{\alpha}{2} - \frac{\beta}{2}}. \label{speed}
\end{equation}

For any positive $\alpha$ and $\beta$ the speed density \eqref{speed} is integrable on the diffusion state space, i.e.
\begin{equation}
\int\limits_{0}^{\infty} \mathfrak{m}(x) \, dx = \frac{\alpha(\beta-2)}{2\theta} \, \alpha^{-\frac{\alpha}{2}} \beta^{-\frac{\beta}{2}} B\left( \frac{\alpha}{2} , \frac{\beta}{2} \right) = M < \infty. \label{speed_int}
\end{equation}
However, for the scale density $\mathfrak{s}(x)$ we have $$\int\limits_{x_{0}}^{\infty} \mathfrak{s}(x) \, dx = \infty, \quad \forall x_{0} \in \langle 0, \infty \rangle,$$ which ensures that starting from the arbitrary point $x_{0}$ from the interior of the diffusion state space the boundary $\infty$ almost surely cannot be attained (cf. A\"it-Sahalia \cite{AitSahalia}). A similar statement holds for the boundary $0$, i.e. provided that diffusion starts from the arbitrary point $x_{0} \in \langle 0, \infty \rangle$, the boundary $0$ almost surely cannot be attained if and only if $$\int\limits_{0}^{x_{0}} \mathfrak{s}(x) \, dx = \infty \quad \Leftrightarrow \quad \alpha \geq 2.$$ We restrict ourselves to this case $\alpha \geq 2$ (but note that in the opposite case a stationary diffusion may also be constructed, subject to instantaneous reflection at $0$).

For any positive $\alpha$ and $\beta$ the stochastic differential equation \eqref{fssde} admits a unique strong Markovian solution $\{X_{t}, \, 0 \leq t \leq T\}$ with time-homogenous transition densities, since it satisfies the following sufficient conditions given by A\"it-Sahalia \cite{AitSahalia}:
\begin{itemize}
\item[(C1)] the drift coefficient $\mu(x)$ and the diffusion coefficient $\sigma(x)$ given by expressions \eqref{dd} are continuously differentiable in $x$ on $\langle 0, \infty \rangle$ and $\sigma^{2}(x)$ is strictly positive for all $x \in \langle 0, \infty \rangle$,

\item[(C2)] the speed measure $\eqref{speed}$ has the property given by \eqref{speed_int}.
\end{itemize}
According to A\"it-Sahalia \cite{AitSahalia}, these conditions are considerably less restrictive than the global Lipschitz and the linear growth conditions which are usually imposed on drift and diffusion coefficients to obtain existence and uniqueness of a strong solution (see Mikosch \cite{Mikosch}).

For any  $\alpha > 2$ and $\beta > 2$ the Fisher-Snedecor diffusion is ergodic (see for example Genon-Catalot et al. \cite{Genon-Catalot}, or S\o rensen \cite{SorensenM}). If furthermore $X_{0} \sim \FS(\alpha, \beta)$, then the Fisher-Snedecor diffusion is strictly stationary. For $\beta > 2$ the conditional expectation satisfy
\begin{equation*}
E \left[ X_{s+t} | X_{s} = x \right] = xe^{-\theta t} + \frac{\beta}{\beta-2} (1 - e^{-\theta t}),
\end{equation*}
and if $\beta > 4$, i.e. if the invariant distribution has finite variance, the autocorrelation function is given by
\begin{equation}
\rho(t) = \mathop{\mathrm{Corr}}(X_{s + t}, X_{s}) = e^{-\theta t}, \quad t \geq 0, \quad s \geq 0 \label{acf}
\end{equation}
(see Bibby et al. \cite{Bibby1}, Theorem 2.3.(iii)).

\begin{rem} \label{rem:mixing}
\emph{According to the general result by Genon-Catalot et al. (see \cite{Genon-Catalot}, Corollary 2.1), the Fisher-Snedecor diffusion
is an $\alpha$-mixing process with an exponentially decaying rate, i.e.
\begin{equation*}
\alpha_{X}(t) = \sup\limits_{s \geq 0}{\alpha(\mathcal{F}_{s}, \mathcal{F}^{s+t})} \leq \frac{1}{4} \, e^{-\delta t}, \quad \delta > 0,
\end{equation*}
where
\begin{equation*}
\alpha(\mathcal{F}_{s}, \mathcal{F}^{s+t}) = \sup\limits_{A \in \mathcal{F}_{s}, B \in \mathcal{F}^{s+t}}{\left| P(A \cap B) - P(A)P(B) \right|}, \end{equation*}
\begin{equation*}
A \in \mathcal{F}_{s} = \sigma\{X_{u}, \, u \leq s\}, \quad B \in \mathcal{F}^{s+t} = \sigma\{X_{u}, \, u \geq s+t\}.
\end{equation*}
Indeed, Fisher-Snedecor diffusion satisfies the sufficient properties by Genon-Catalot et al.:
\begin{itemize}
\item[(i)] the drift coefficient $\mu(x)$ and the squared diffusion coefficient $\sigma^{2}(x)$, given by expressions \eqref{dd}, satisfy condition (C1). Furthermore, there exists a strictly positive constant $$K = \displaystyle\frac{\theta (\beta + 2\sqrt{\alpha})^{2}}{\alpha (\beta - 2)}$$ such that
\begin{equation*}
\begin{tabular}{ccc}
$\left| \mu(x) \right| \leq K \, (1 + |x|)$ & \textrm{and} & $\sigma^{2}(x) \leq K \, (1 + x^{2})$. \\
&  &
\end{tabular}
\end{equation*}
\item[(ii)] The speed density $\mathfrak{m}(x)$ given by \eqref{speed} is integrable on the diffusion state space (see expression \eqref{speed_int}). For $\alpha > 2$ the scale density $\mathfrak{s}(x)$ given by \eqref{scale} is non-integrable in the neighborhood of boundary points $0$ and $\infty$.
\item[(iii)] The random variable $X_{0}$ has (Fisher-Snedecor) density function which is proportional to the speed density $\mathfrak{m}(x)$, i.e.
\begin{equation*}
\mathop{\mathrm{\mathfrak{fs}}}(x) = \frac{\mathfrak{m}(x)}{M} \, {\mathop{\mathrm{I}}}_{\langle 0, \infty \rangle}(x).
\end{equation*}
\item[(iv)] The product of the diffusion coefficient $\sigma(x)$ and the speed density $\mathfrak{m}(x)$ converge to $0$ as $x \to 0$ and $x \to \infty$, i.e.  $$\lim\limits_{x \to 0} \sigma(x)\mathfrak{m}(x) = \lim\limits_{x \to \infty} \sigma(x)\mathfrak{m}(x) = 0.$$
\item[(v)] If we define
\begin{equation*}
\gamma(x) = \sigma^{\prime}(x) - \frac{2\mu(x)}{\sigma(x)} = \frac{\beta\left[ \alpha(x-1) + 1 \right]}{x(\alpha x + \beta)} \, \sqrt{\frac{\theta x (\alpha x + \beta)}{\alpha (\beta - 2)}}, \quad x \in \langle 0, \infty \rangle,
\end{equation*}
then it follows
\begin{equation*}
\lim\limits_{x \to 0} \frac{1}{\gamma(x)} = \sqrt{\frac{\alpha(\beta - 2)}{\beta^{2} \theta}} \lim\limits_{x \to 0} \frac{\sqrt{x(\alpha x + \beta)}}{\alpha(x-1) + 1} = 0 < \infty,
\end{equation*}
\begin{equation*}
\lim\limits_{x \to \infty} \frac{1}{\gamma(x)} = \sqrt{\displaystyle\frac{\alpha(\beta - 2)}{\beta^{2} \theta}} \lim\limits_{x \to \infty} \displaystyle\frac{\sqrt{x(\alpha x + \beta)}}{\alpha(x-1) + 1} = \frac{1}{\beta} \, \sqrt{\frac{\beta - 2}{\theta}} < \infty.
\end{equation*}
\end{itemize} \vspace*{3mm}
Detailed exposition of the theory of mixing processes and additional information on corresponding central limit theorems are given by Doukhan \cite{Doukhan}.}
\end{rem}

\begin{rem}
\emph{The autonomous stochastic differential equation \eqref{fssde} is a special case of the generalized non-linear mean reversion A\"it-Sahalia model
\begin{equation*}
dX_{t} = (c_{-1}X_{t}^{-1}+c_{0}+c_{1}X_{t}+c_{2}X_{t}^{2}) \, dt + \sqrt{\sigma^2_{0}+\sigma^2_{1}X_{t}+\sigma^2_{2}X_{t}^{\delta}} \, dW_{t}
\end{equation*}
for following parameter values:
\begin{equation*}
\begin{tabular}{lcl}
$c_{-1} = 0$, & \phantom{aaa} & $\sigma^2_{0} = 0$, \\
$c_{0} = \frac{\theta \beta}{\beta-2} > 0$, & \phantom{aaa} &
$\sigma^2_{1}
= \frac{4 \theta \beta}{\alpha(\beta-2)} > 0$, \\
$c_{1} = -\theta < 0$, & \phantom{aaa} & $\sigma^2_{2} = \frac{4 \theta}{\beta-2} > 0$, \\
$c_{2} = 0$, & \phantom{aaa} & $\delta = 2$. \\
&  &
\end{tabular}
\end{equation*}
The generalized non-linear mean reversion A\"it-Sahalia model is frequently used for interest rates modeling. For more details about this model see \cite{AitSahalia} or \cite{Iacus}.}
\end{rem}

\begin{rem}
\emph{The Lamperti transform $\{Y_{t}, \, t \geq 0\}$ of the Fisher-Snedecor diffusion is a solution of the stochastic differential equation
\begin{equation*}
dY_{t} = \frac{\alpha - \beta - 1}{2} \, \sqrt{\frac{\theta}{\beta - 2}} \tanh{\left( Y_{t} \sqrt{\frac{\beta \theta}{\beta - 2}} \right)} \, dt + dW_{t}, \quad t \geq 0
\end{equation*}
with unit diffusion parameter. According to Iacus \cite{Iacus}, Lamperti transforms of diffusion processes are recommended for use in simulation studies.}
\end{rem}

\medskip

\section{Spectral representation of the transition density of Fisher-Snedecor diffusion}

\subsection{Fokker-Planck equation}

In this section we present the spectral representation of the transition density
\begin{equation}
p = p(x, t) = p(x; x_{0}, t) = \frac{d}{dx} \, P(X_{t} \leq x \mid X_{0} = x_{0}), \quad x > 0, \quad t \geq 0, \label{trden}
\end{equation}
of the Fisher-Snedecor diffusion in terms of solutions of the corresponding Sturm-Liouville equation. Note that the corresponding Sturm-Liouville operator is closely related to the infinitesimal generator of the Fisher-Snedecor diffusion.

According to Karlin and Taylor \cite{Karlin2}, the transition density \eqref{trden} is the principal solution of the Fokker-Planck or the forward Kolmogorov equation
\begin{equation*}
\frac{\partial p}{\partial t} = -\frac{\partial}{\partial x} \left(-\theta \left(x - \frac{\beta}{\beta - 2} \right) \, p \right) + \frac{1}{2} \displaystyle\frac{\partial^{2}}{\partial x^{2}} \left( \displaystyle\frac{4 \theta}{\alpha(\beta - 2)} x(\alpha x + \beta) \, p \right), \quad x > 0, \quad t \geq 0, \label{fokkerplanck}
\end{equation*}
and its Laplace transform is known explicitly (see \eqref{green}).

Inverting the Laplace transform of the transition density yields its spectral decomposition - see for example the classical books by Titchmarsh \cite{Titchmarsh}, Karlin and Taylor \cite{Karlin2} and It\^o and McKean \cite{Ito}, as well as the paper by McKean \cite{McKean} (who call this "eigendifferential expansion").

The discrete part of the spectrum can be treated in a unified manner for all Pearson diffusions - see Wong \cite{Wong} and Forman and Sorensen \cite{Forman}. The exact formulas for the continuous part of the spectrum require however computations which need to be performed on a case by case basis. Some of such computations for special cases of Pearson diffusions were performed recently, motivated by applications in finance, by Linetsky \cite{Linetsky1, Linetsky2, Linetsky3}, Davydov and Linetsky \cite{Davydov} and Carr et al. \cite{Mendoza}. We treated here the Fisher-Snedecor case,  and we chose, for self containedness, to include a detailed account of the inversion of the Laplace transform of the corresponding transition density.

\subsection{Sturm-Liouville equation}

The infinitesimal generator of the Fisher-Snedecor diffusion is given by the expression
\begin{equation}
(\mathcal{G}f)(x) = \frac{2 \theta}{\alpha(\beta - 2)} \, x^{1 - \frac{\alpha}{2}}(\alpha x + \beta)^{\frac{\alpha}{2} + \frac{\beta}{2}} \, \frac{d}{dx} \, \left(x^{\frac{\alpha}{2}}(\alpha x + \beta)^{1 - \frac{\alpha}{2} - \frac{\beta}{2}} f^{\prime}(x) \right) = \label{infgen}
\end{equation}
\begin{equation*}
\hspace*{7.5mm} = \frac{2\theta}{\alpha(\beta - 2)} \, x(\alpha x + \beta) f^{\prime \prime}(x) - \theta \left(x - \frac{\beta}{\beta - 2} \right)f^{\prime}(x), \quad x > 0,
\end{equation*}
acting on the domain
$$D(\mathcal{G}) = \left\{ f \in L^{2}(\langle 0, \infty \rangle, \fs(x)) \cap C^{2}(\langle 0, \infty \rangle): \mathcal{G}f \in L^{2}(\langle 0, \infty \rangle, \fs(x)), \right.$$ $$\left. \lim\limits_{x \to 0} \frac{f^{\prime}(x)}{\mathfrak{s}(x)} = \lim\limits_{x \to \infty} \frac{f^{\prime}(x)}{\mathfrak{s}(x)} = 0 \right\},$$ where $\mathfrak{s}(x)$ is the scale density \eqref{scale}. The negative of the infinitesimal generator $\mathcal{G}$ is called the Sturm-Liouville operator. The associated Sturm-Liouville differential equation $(-\mathcal{G}f)(x) = \lambda f(x)$ for the Fisher-Snedecor diffusion takes the form
\begin{equation}
\frac{2\theta}{\alpha(\beta - 2)} \, x(\alpha x + \beta) f^{\prime \prime}(x) - \theta \left(x - \frac{\beta}{\beta - 2} \right) f^{\prime}(x) + \lambda f(x) = 0, \label{sleq}
\end{equation}
where $\lambda \geq 0$ is the spectral parameter (see Karlin and Taylor \cite{Karlin2} and Linetsky \cite{Linetsky1}).

\begin{rem}
\emph{We consider also the equation
\begin{equation}
\mathcal{G} f(x) = s f(x), \quad s > 0,  \label{sleqs}
\end{equation}
(with $s = -\lambda$, where $\lambda$ is the spectral parameter from the equation \eqref{sleq}). According to It\^o and McKean (see \cite{Ito}, page 128) or Borodin and Salminen (see \cite{BS}, page 18), for any $s > 0$ it admits two positive, linearly independent and increasing/decreasing solutions with Wronskian proportional to the scale density, called fundamental solutions. We will denote these solutions respectively by $f_{1}(x) = f_{1}(x, s)$ and $f_{4}(x) = f_{4}(x, s)$ (see \eqref{f1l} and \eqref{f34l}). They intervene in the expression of Green's function, as well as in several first-passage problems.}

\begin{equation*}
\begin{tabular}{c}
{\includegraphics[scale=0.7]{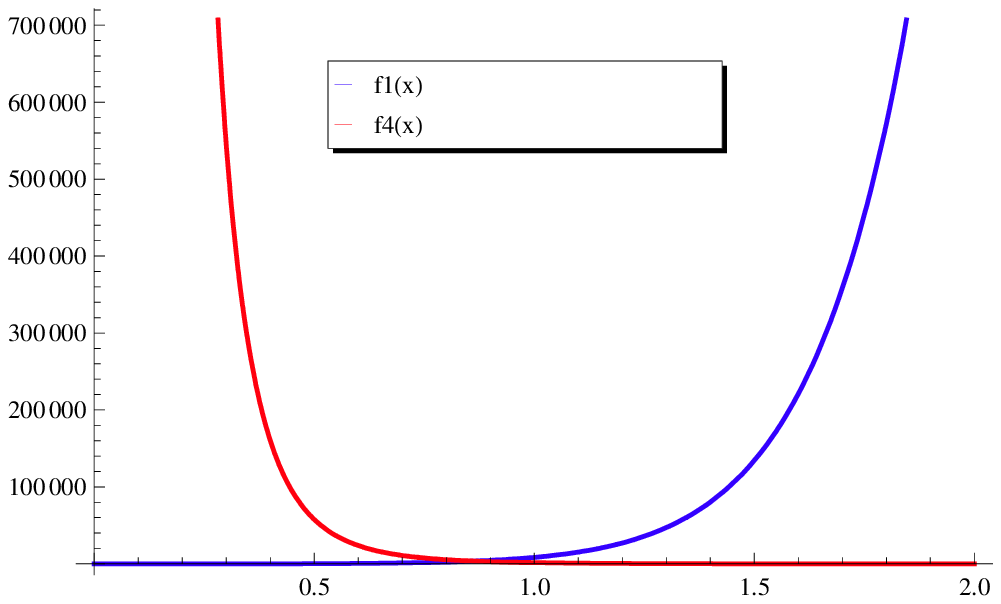}}
\end{tabular}
\end{equation*}
\begin{center}
{\footnotesize{\scshape{Figure 1.}} \, \, {\textit{Graphs of fundamental solutions of equation \eqref{sleqs} for $\alpha = 5$, $\beta = 20$ and $\theta = 0.05$.}}}
\end{center}
\end{rem}

The Sturm-Liouville equation \eqref{sleq} may be transformed by using the substitution
\begin{equation*}
f(x) = v\left( y(x) \right) = v \left( \frac{\alpha}{\beta} \, x \right),
\end{equation*}
to an equation of the hypergeometric type (see Appendix B):
\begin{equation}
y(y+1)v^{\prime \prime}(y) + \left(\frac{\alpha}{2} + \left( 1 - \frac{\beta}{2} \right) \, y \right)v^{\prime}(y) + \l^*v(y) = 0, \label{hgeq}
\end{equation}
where $\l^* = \frac{\lambda(\beta-2)}{2\theta}$.

Denote the roots of the quadratic equation $z^2 + \frac{\beta}{2} z + \l^*=0$ by
\begin{equation} \label{Delta}
z_\pm = z_{\pm, \l} = -\frac{\beta}{4} \pm \D_\l, \quad \D = \D_\l = \sqrt{\frac{\beta^{2}}{16} - \frac{\lambda(\beta-2)}{2\theta}}.
\end{equation}

According to Titchmarsh (see \cite{Titchmarsh}, Example 4.18., page 100), two linearly independent solutions of equation \eqref{hgeq}
for $|y| < 1$ are
\begin{equation*}
v_{1}(y) = \phantom{,}_{2}F_{1} \left(z_+, z_-; \frac{\alpha}{2}; -y \right),
\end{equation*}
\begin{equation*}
v_{2}(y) = y^{1 - \frac{\alpha}{2}} \, _{2}F_{1} \left(u_+, u_-; 2 - \frac{\alpha}{2}; -y \right), \quad \alpha > 2, \quad \alpha \notin \{2(m+1), \, m \in \mathbb{N}\},
\end{equation*}
where $u_\pm = u_{\pm,\l} = 1 - \frac{\alpha}{2} +z_\pm $. Furthermore, two linearly independent solutions for $|y| > 1$ are
\begin{equation*}
v_{3, 4}(y) = y^{-z_\mp } \, \phantom{,}_{2}F_{1} \left(z_ \mp , u_\mp ;  1 \mp 2 \D_\l; -\frac{1}{y} \right),
\end{equation*}
where the upper sign refers to the solution $v_{3}(y)$ and the lower sign refers to the solution $v_{4}(y)$.

\begin{rem} \label{f2EvenAlpha}
\emph{Cf. Abramowitz and Stegun (see \cite{Abramowitz}, Section 15, Expression 15.5.21.), the existence of the solution $v_{2}(y)$ follows from the condition $\alpha \notin \{2(m+1), \, m \in \mathbb{N}\}$ imposed on the value of parameter $\alpha$. However, if $\alpha > 2$ is an even integer, then the corresponding solution is given by much more complicated expression in terms of the digamma function $\psi(\cdot)$ which is hardly evaluated in explicit calculations.}
\end{rem}

\begin{rem}
\emph{According to Luke (see \cite{Luke}, Section 3.9., Expression (1); see also Slater \cite{Slater}, Section 1.8.1., Expression 1.8.1.11.), the solution $v_{1}(y)$ valid in the region $|y| < 1$ can be analytically continued to the whole complex plane cut along the interval $\langle -\infty, -1]$ (see Appendix B). This continuation is provided by the following expression which represents the solution $v_{1}(y)$ as the linear combination of the solutions $v_{3}(x)$ and $v_{4}(x)$, i.e.
\begin{equation}
v_{1}(x) = B_\l v_{3}(x) + A_\l v_{4}(x) \label{f1=Bf3+Af4}
\end{equation}
where
\begin{equation}
B = B_\l = \frac{\Gamma\left( \frac{\alpha}{2} \right) \, \Gamma\left( 2 \D_\l \right)}{%
\Gamma\left( z_{+, \l} \right) \, \Gamma\left( 1 - u_{-, \l}
\right)}, \quad
A = A_\l =\frac{\Gamma\left( \frac{\alpha}{2} \right) \, \Gamma\left( -2 \D_\l \right)}{%
\Gamma\left( z_{-, \l} \right) \, \Gamma\left( 1 - u_{+, \l}
\right)}. \label{B}
\end{equation}
Similarly, solutions $v_{3}(y)$ and $v_{4}(y)$, valid in the region $|y| > 1$, could be analytically continued to the whole complex plane cut along the negative real axis. For details on the corresponding continuation formulas we refer to Luke (see \cite{Luke}, Section 3.9., Expressions (3) and (4)). This discussion shows that we are dealing with functions which are all analytic on the complex plane cutted along the subset of the negative real axis.}
\end{rem}

Changing back variables, the solutions of Sturm-Liouville equation \eqref{sleq} are
\begin{equation}
f_{1}(x) = f_{1}(x, -\lambda) \phantom{a} = \phantom{,}_{2}F_{1} \left(z_+, z_-; \frac{\alpha}{2}; -\frac{\alpha}{\beta} \, x \right),  \label{f1l}
\end{equation}
\begin{equation}
f_{2}(x) = f_{2}(x, -\lambda) = \left( \frac{\alpha}{\beta} \, x \right)^{1 - \frac{\alpha}{2}} \, _{2}F_{1} \left(u_+, u_-; 2-\frac{\alpha}{2}; -\frac{\alpha}{\beta} \, x \right),
\end{equation}
\centerline{$\alpha > 2, \quad \alpha \notin \{2(m+1), \, m \in \mathbb{N}\},$}
\begin{equation}
f_{3, 4}(x) = f_{3, 4}(x, -\lambda) = \left( \frac{\alpha}{\beta} \,
x \right)^{-z_{\mp}} \, \phantom{,}_{2}F_{1} \left(z_ \mp , u_\mp ;
1 \mp 2 \D_\l; -\frac{\beta}{\alpha x} \right), \label{f34l}
\end{equation}
where the upper sign refers to the solution $f_{3}(x, -\lambda)$ and the lower sign refers to the solution $f_{4}(x, -\lambda)$.

Due to the analytic continuation of the hypergeometric functions, solutions $f_{1}(x, -\lambda)$, $f_{3}(x, -\lambda)$ and $f_{4}(x, -\lambda)$ are all analytic on the state space of the Fisher-Snedecor diffusion. Therefore, fundamental solutions of equation \eqref{sleqs} on $\langle 0, \infty \rangle$ could be identified with the increasing solution $f_{1}(x, s)$ and the decreasing solution $f_{4}(x, s)$, see Figure 1.

\begin{rem}
\emph{According to Borodin and Salminen (\cite{BS}, page 19) and Buchholz \cite{Buchholz}, Wronskian of linearly independent solutions $f_{1}(x, s)$ and $f_{4}(x, s)$ of Sturm-Liouville equation \eqref{sleq} is proportional to the scale density \eqref{scale}, i.e.
\begin{equation}
W(f_{4}, f_{1}) = W_{s}(f_{4}, f_{1})(x) = W_{s}(f_{4}, f_{1})(x_0) \, \mathfrak{s}(x),
\end{equation}
where the value $x_{0} > 0$ could be chosen so that $f_{1}(x_{0}, s)$, $f_{1}^{\prime}(x_{0}, s)$, $f_{4}(x_{0}, s)$ and $f_{4}^{\prime}(x_{0}, s)$ or the limit $\lim_{x \to x_{0}} W \left( f_{4}(x, s), f_{1}(x, s) \right)$ are easy to calculate (see Titchmarsh \cite{Titchmarsh}). According to expression \eqref{Wf1f4gen} from Appendix C, Wronskian of solutions $f_{4}(x, s)$ and $f_{1}(x, s)$ is
\begin{equation}
W_{s}(f_{4}, f_{1}) = 2 \mathfrak{s}(x) B_s \alpha^{1 - \frac{\alpha}{2}} \beta^{-\frac{\beta}{2}} \, \sqrt{\frac{\beta^{2}}{16} + \frac{s(\beta - 2)}{2 \theta}} = 2 \mathfrak{s}(x) \alpha^{1 - \frac{\alpha}{2}} \beta^{-\frac{\beta}{2}} B_s \D_s, \label{Wf1f4s}
\end{equation}
where $B_s$ is given by \eqref{B} and $\D_s$ by \eqref{Delta} for $s = -\l$.}
\end{rem}

\medskip

\subsection{Spectral representation of transition density}

As well-known (see e.g. Borodin and Salminen  \cite{BS}, page 19), the Laplace transform of transition density $p(x; x_{0}, t),$ also known as the resolvent kernel or Green's function
\begin{equation*}
G_{s}(x_{0}, x) = \int\limits_{0}^{\infty} e^{-st} p(x; x_{0}, t) \, dt,
\end{equation*}
admits the explicit representation
\begin{equation}
G_{s}(x_{0}, x) = \frac{\mathfrak{m}(x)}{w_{s}(\varphi_{s}, \psi_{s})} \, \psi_{s}(x_{0} \wedge x) \varphi_{s}(x_{0} \vee x), \label{green}
\end{equation}
where $\mathfrak{m}(\cdot)$ is the speed density, $\psi_{s}(\cdot)$ and $\varphi_{s}(\cdot)$ are the fundamental solutions of the equation  \eqref{sleqs}, and $w_{s}(\varphi_{s}, \psi_{s})$ is their Wronskian (cf. Linetsky \cite{Linetsky1}) with respect to the scale density $\mathfrak{s}(x)$, i.e.
\begin{equation*}
w_{s}(\varphi_{s}, \psi_{s}) = \frac{1}{\mathfrak{s}(x)} \left( \psi_{s}^{\prime}(x) \varphi_{s}(x) - \psi_{s}(x) \varphi^{\prime}_{s}(x) \right) = \frac{W_{s}(\varphi_{s}, \psi_{s})}{\mathfrak{s}(x)}.
\end{equation*}
Therefore, if it is possible to determine the explicit form of the fundamental solutions, then the transition density can be obtained by the
Laplace inversion formula
\begin{equation}
p(x; x_{0}, t) = \frac{1}{2 \pi i} \int\limits_{c - i \infty}^{c + i \infty} e^{st} \, G_{s}(x_{0}, x) \, ds, \quad t > 0, \label{invLap}
\end{equation}
where $s > 0$ is the spectral parameter from equation \eqref{sleqs} and where integration is performed along a line $\mathrm{Re}(s) = c$, $c > 0$, leaving all the singularities of the Green's function to its left. This approach implies the importance of understanding the nature of the spectrum of the Sturm-Liouville operator $(-\mathcal{G})$. Therefore, clarification of spectral properties regarding operator $(-\mathcal{G})$ for Fisher-Snedecor diffusion is given in the following remark.

\begin{rem} \label{SpecCat}
\emph{According to the Theorem \ref{boundary_class}. from Appendix C, boundaries of the state space of Fisher-Snedecor diffusion are classified as follows: $0$ is non-oscillatory boundary (regular for $\alpha \leq 2$ and entrance for $\alpha > 2$), while $\infty$ is oscillatory/non-oscillatory singular boundary (natural for all $\alpha > 0$) with unique cutoff
\begin{equation}
\Lambda = \frac{\theta \beta^{2}}{8(\beta-2)}, \quad \beta > 2.  \label{cutoff}
\end{equation}
Furthermore, $\infty$ is non-oscillatory for $\lambda = \Lambda$ (see Theorem \ref{boundary_class}.(ii)). According to this classification of boundaries of the diffusion state space, the Sturm-Liouville operator $(-\mathcal{G})$ has a finite set of simple eigenvalues
\begin{equation}
\lambda_{n} = \frac{\theta}{\beta - 2} \, n (\beta - 2n), \quad n = 0, \ldots, \left\lfloor \frac{\beta}{4} \right\rfloor, \quad \beta > 2,
\label{diseigen}
\end{equation}
in $[0, \Lambda]$. Therefore, the discrete spectrum of the operator $(-\mathcal{G})$ is the finite set $\sigma_{d}(-\mathcal{G}) = \left\{\lambda_{n}, \, n = 0, \ldots, \left\lfloor \beta/4 \right\rfloor \right\}$. Furthermore, the operator $(-\mathcal{G})$ has the essential spectrum $\sigma_{ess}(-\mathcal{G}) = [\Lambda, \infty \rangle$ which contains the purely absolutely continuous spectrum of multiplicity one in $\langle \Lambda, \infty \rangle$, i.e. $\sigma_{ac}(-\mathcal{G}) = \langle \Lambda, \infty \rangle$.
The elements of the absolutely continuous spectrum could be parameterized as
\begin{equation}
\lambda = \Lambda + \frac{2\theta k^{2}}{\beta-2} = \frac{2\theta}{\beta-2} \left( \frac{\beta^{2}}{16} + k^{2} \right), \quad \beta > 2, \quad k > 0.
\label{conteigen}
\end{equation}
According to Linetsky's \cite{Linetsky1} spectral classification of one-dimensional diffusions based on the nature of the spectrum of Sturm-Liouville operator $(-\mathcal{G})$, Fisher-Snedecor diffusion belongs to the spectral category II.}
\end{rem}

According to It\^o and McKean \cite{Ito} (see also Linetsky \cite{Linetsky1}), the general form of the spectral representation of transition density \eqref{trden} for diffusions belonging to spectral category II is of the following form:
\begin{equation}
p(x; x_{0}, t) = \mathfrak{m}(x) \, \displaystyle\sum_{n=0}^{N} e^{-\lambda_{n} t} \varphi_{n}(x_{0}) \varphi_{n}(x) + \mathfrak{m}(x) \, \int\limits_{\Lambda}^{\infty} e^{-\lambda t} \varphi(x_{0}, -\lambda) \varphi(x, -\lambda) \, d\lambda, \label{gsr}
\end{equation}
where $\lambda_{n}$, $n \in \{0, \ldots, N\}$, are elements of the discrete spectrum $\sigma_{d}(-\mathcal{G})$ and $\varphi_{n}(\cdot)$ are corresponding eigenfunctions normalized with respect to the speed density $\mathfrak{m}(x)$, while $\lambda > \Lambda$ are elements of the absolutely continuous spectrum $\sigma_{ac}(-\mathcal{G})$ and $\varphi(\cdot, -\lambda)$ are solutions of the equation \eqref{sleq} also normalized with respect to the speed density $\mathfrak{m}(x)$.

Now we proceed to compute the spectral representation of the transition density for our process of interest following the general procedure from It\^o and McKean \cite{Ito} and Linetsky \cite{Linetsky1} and using the results of the analysis of the spectrum of the corresponding Sturm-Liouville operator (see Appendix C and Remark \ref{SpecCat}).

\newpage

\begin{thm} \label{thm:SpecRep} Spectral representation of the transition density of ergodic stationary diffusion with marginal Fisher-Snedecor distribution with parameters $\alpha > 2$, $\alpha \notin \{2(m+1), \, m \in \mathbb{N}\}$, and $\beta > 2$ is of the form
\begin{equation}
p(x; x_{0}, t) = p_{d}(x; x_{0}, t) + p_{c}(x; x_{0}, t) \label{d+c}.
\end{equation}
The discrete part of the spectral representation
\begin{equation}
p_{d}(x; x_{0}, t) = \fs(x) \, \sum\limits_{n=0}^{\left\lfloor \frac{\beta}{4} \right\rfloor}e^{-\lambda_{n} t} \, F_{n}(x_{0}) \, F_{n}(x) \label{disspec}
\end{equation}
is given in terms of the eigenvalues $\lambda_{n}$ given by \eqref{diseigen} and the normalized Fisher-Snedecor polynomials $F_{n}(\cdot)$ given by \eqref{fphyp}. The continuous part of the spectral representation
\begin{equation}
p_{c}(x; x_{0}, t) = \fs(x) \, \frac{1}{\pi} \int\limits_{\frac{\theta \beta^{2}}{8(\beta-2)}}^{\infty} e^{-\lambda t} \, k(\lambda) \times \label{contspec}
\end{equation}
\begin{equation*}
\times \left| \frac{ B^{\frac{1}{2}}\left( \frac{\alpha}{2}, \frac{\beta}{2} \right) \, \Gamma\left( -\frac{\beta}{4} + ik(\lambda) \right) \Gamma\left( \frac{\alpha}{2} + \frac{\beta}{4} + ik(\lambda) \right)}{\Gamma\left(\frac{\alpha}{2} \right) \, \Gamma\left(1 + 2ik(\lambda) \right)} \right|^{2} \, f_{1}(x_{0}, -\lambda) f_{1}(x, -\lambda) \, d\lambda
\end{equation*}
is given in terms of the elements $\lambda$ of the absolutely continuous spectrum of the operator $(-\mathcal{G})$ given by \eqref{conteigen}, solution $f_{1}(\cdot, -\lambda)$ of the Sturm-Liouville equation \eqref{sleq} given by \eqref{f1l} and parameter $k(\lambda) = -i \D_\l$, where $\D_\l$ is given in \eqref{Delta}.
\end{thm}

\begin{proof}
Since Wronskian of fundamental solutions $f_{1}(x, s)$ and $f_{4}(x, s)$ of equation \eqref{sleqs} is given by expression \eqref{Wf1f4s}, it is possible to obtain the explicit form of the Green's function for Fisher-Snedecor diffusion. In particular,
\begin{equation}
G_{s}(x_{0}, x) = \frac{\mathfrak{m}(x)}{2 \alpha^{1-\frac{\alpha}{2}}
\beta^{-\frac{\beta}{2}} B_{s} \D_{s}} \, f_{1}(x \wedge x_{0}, s) f_{4}(x \vee x_{0}, s), \label{greenfs}
\end{equation}
where $\mathfrak{m}(\cdot)$ is the speed density given by \eqref{speed}, $f_{1}(\cdot, s)$ and $f_{4}(\cdot, s)$ are non-normalized (with respect to the speed density) linearly independent solutions of the differential equation \eqref{sleqs} and
\begin{equation*}
W_{s}(f_{4}, f_{1}) = 2 B_{s} \alpha^{1 - \frac{\alpha}{2}} \beta^{-\frac{\beta}{2}} \, \sqrt{\frac{\beta^{2}}{16} + \frac{s(\beta - 2)}{2 \theta}} = 2 \alpha^{1 - \frac{\alpha}{2}} \beta^{-\frac{\beta}{2}} B_s \D_s
\end{equation*}
is their Wronskian with respect to the scale density \eqref{scale}, where $B_{s}$ is given by \eqref{B} and $\D_{s}$ is given by \eqref{Delta}. Observe now the Green's function \eqref{greenfs} as the function of the complex variable $s$. We conclude that the term
\begin{equation*}
\Gamma\left( z_{+, s} \right) = \Gamma\left( -\frac{\beta}{4} + \sqrt{\frac{\beta^{2}}{16} + \frac{s(\beta-2)}{2\theta}} \right)
\end{equation*}
from $B_{s}$ in expression \eqref{greenfs} has simple poles in
\begin{equation*}
s = -\lambda_{n} = -\frac{\theta}{\beta - 2} \, n (\beta - 2n), \quad n = 0, \ldots, \left\lfloor \frac{\beta}{4} \right\rfloor,
\end{equation*}
since for that value of $s$ we have $\Gamma\left( z_{+, \lambda_{n}} \right) = \Gamma \left( -n \right)$, $n = 0, \ldots, \left\lfloor \beta/4 \right\rfloor$. Note that these simple poles coincide with the negative simple eigenvalues of the Sturm-Liouville operator $(-\mathcal{G})$ and that, since $\alpha > 2$, $\alpha \notin \{2(m+1), \, m \in \mathbb{N}\}$, and $\beta > 2$, these are the only poles of the Green's function \eqref{greenfs}.

According to Abramowitz and Stegun \cite{Abramowitz}, $f_{4}(x, s) = C_{s} f_{1}(x, s) + D_{s} f_{2}(x, s)$, where
\begin{equation*}
C_{s} = \frac{\Gamma\left(1 - \frac{\alpha}{2} \right) \, \Gamma\left(1 + 2 \D_s \right)}{\Gamma\left(u_{+, s} \right) \, \Gamma\left(1 - z_{-, s} \right)} \quad \textrm{and} \quad D_{s} = \frac{\Gamma\left(\frac{\alpha}{2} - 1 \right) \, \Gamma\left(1 + 2 \D_s \right)}{\Gamma\left(z_{+, s} \right) \, \Gamma\left(1 - u_{-, s} \right)}.
\end{equation*}
Therefore, regarding the fact that Green's function has simple poles at $s = -\lambda_{n}$, expression \eqref{greenfs} can be written in terms of solutions $f_{1}(\cdot, s)$ and $f_{2}(\cdot, s)$, i.e.
\begin{equation}
G_{s = -\lambda_{n}}(x_{0}, x) = \mathfrak{m}(x) \, \left\{ \frac{\Gamma\left( -n \right) \Gamma\left( \frac{\alpha}{2} + \frac{\beta}{2} - n \right)}{\alpha^{1-\frac{\alpha}{2}} \beta^{-\frac{\beta}{2}} \Gamma\left( \frac{\alpha}{2} \right) \Gamma\left( 1 + \frac{\beta}{2} - n \right)} \, (-1)^{n} \, \prod_{j=0}^{n-1} \left( \frac{\alpha}{2} + j \right) \, f_{1}(x_{0}, s) f_{1}(x, s) + \right. \label{greenf1f2}
\end{equation}
\begin{equation*}
\hspace*{1.3cm} \left. + \frac{\Gamma\left( \frac{\alpha}{2} - 1 \right)}{\alpha^{1-\frac{\alpha}{2}} \beta^{-\frac{\beta}{2}} \Gamma\left( \frac{\alpha}{2} \right)} \, f_{1}(x \wedge x_{0}, s) f_{2}(x \vee x_{0}, s) \right\}.
\end{equation*}
Due to supposed values of parameter $\alpha$, the second term in expression \eqref{greenf1f2} has no poles and therefore does not contribute to the residues of the Green's function. However, due to simple poles at $s=-\lambda_{n}$ and according to Titchmarsh (see \cite{Titchmarsh}, Chapter IV, Example 4.19.), the residues of Green's function at these poles are given by
\begin{equation*}
{\Res}_{s=-\lambda_{n}}G_{s}(x_{0}, x) = \mathfrak{m}(x) \, \frac{(-1)^{2n}}{n!} \, \frac{\left( \frac{\beta}{2} - 2n \right) \Gamma\left( \frac{\alpha}{2} + \frac{\beta}{2} -n \right)}{\alpha^{1-\frac{\alpha}{2}} \beta^{-\frac{\beta}{2}} \Gamma\left( \frac{\alpha}{2} \right) \Gamma\left( 1 + \frac{\beta}{2} - n \right)} \, \prod\limits_{j=0}^{n-1}\left( \frac{\alpha}{2} + j \right) \times
\end{equation*}
\begin{equation}
\hspace*{3.5cm} \times \, \phantom{,}_{2}F_{1} \left(-n, n - \frac{\beta}{2}; \frac{\alpha}{2}; -\frac{\alpha}{\beta} \, x_{0} \right) \, \phantom{,}_{2}F_{1} \left(-n, n - \frac{\beta}{2}; \frac{\alpha}{2}; -\frac{\alpha}{\beta} \, x \right). \label{greenres}
\end{equation}
In expression \eqref{greenres} we can recognize non-normalized (with respect to the speed density) Fisher-Snedecor polynomials given by formula \eqref{fphyp}.

Furthermore, Green's function \eqref{greenfs} has a branch point at $s = -\Lambda = -\frac{\theta \beta^{2}}{8(\beta - 2)}$, since for $s < -\Lambda$ the gamma function $\Gamma\left( z_{+, s} \right)$ has the argument with non-zero imaginary part. Note that this branch point coincides with the negative cutoff between the discrete and the continuous part of the spectrum of the Sturm-Liouville operator $(-\mathcal{G})$. The branch cut of discontinuity is placed from $-\Lambda$ to $-\infty$ on the negative part of the real axis and is parameterized as
\begin{equation*}
s = -\lambda = -\frac{2 \theta}{\beta - 2} \left( \frac{\beta^{2}}{16} + k^{2} \right), \quad k > 0.
\end{equation*}

After this analysis we can start with the evaluation of the inverse Laplace transform of Green's function \eqref{greenfs} using the inversion formula \eqref{invLap}. This procedure results in the spectral representation of transition density of Fisher-Snedecor diffusion. First we observe the Bromwich contour $\mathcal{C}$ as shown in the following figure.

\begin{equation*}
\begin{tabular}{c}
{\includegraphics[scale=0.28]{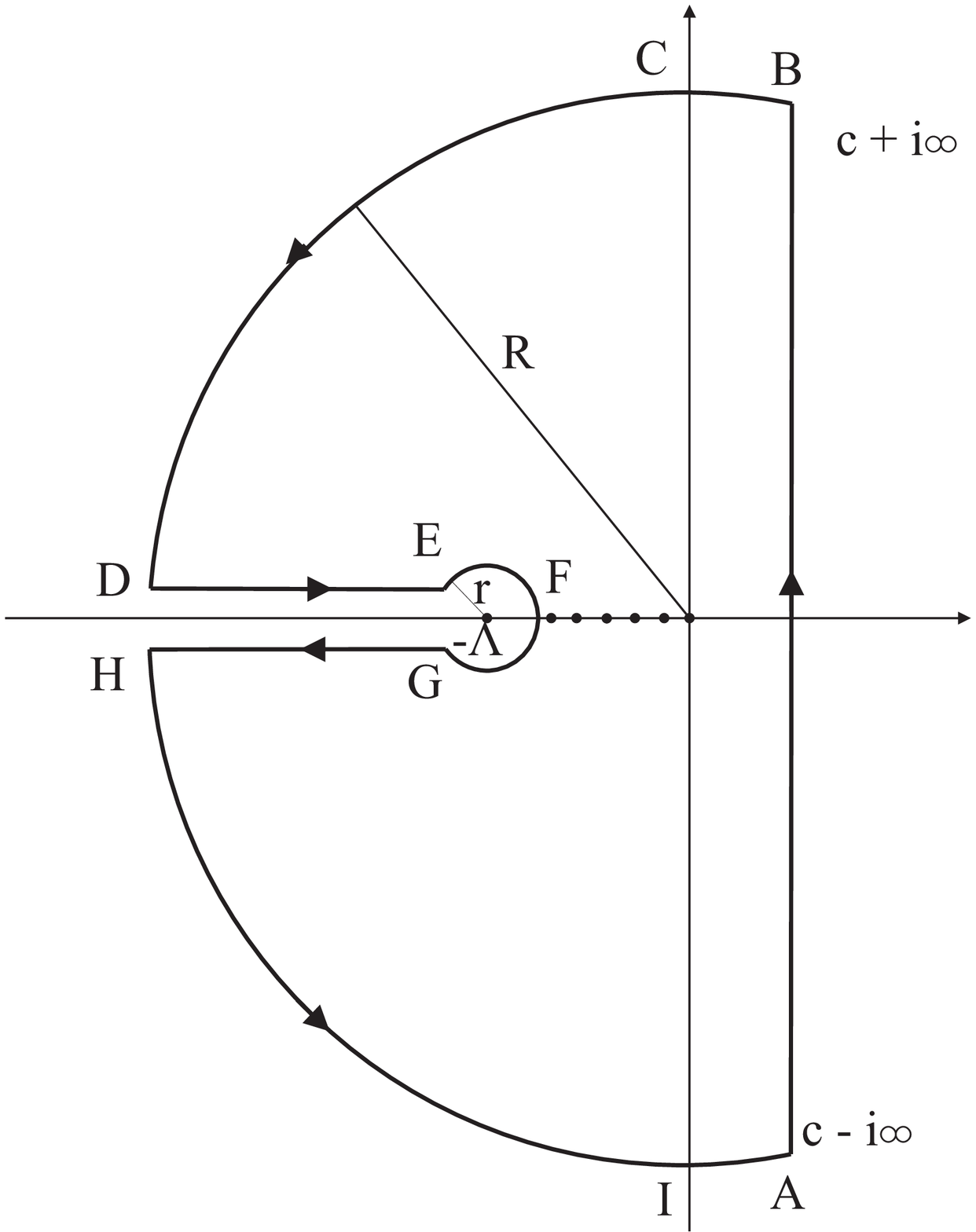}}
\end{tabular}%
\end{equation*}
\begin{center}
{\footnotesize {\scshape{Figure 2.}} \quad \textit{The Bromwich
contour $\mathcal{C}$, simple poles, branch point and branch cut of
Green's function.}}
\end{center}

Since simple poles $(-\lambda_{n})$ of the Green's function are placed inside the contour $\mathcal{C}$, according to the Cauchy Residue Theorem it follows that
\begin{equation}
\frac{1}{2 \pi i} \int\limits_{\mathcal{C}} e^{st} \, G_{s}(x_{0}, x) \, ds = \sum_{n=0}^{\left\lfloor \frac{\beta}{4} \right\rfloor} {\Res}_{s=-\lambda_{n}}e^{st} G_{s}(x_{0}, x). \label{BromwichC}
\end{equation}

On the other hand, the integral around the contour $\mathcal{C}$ is equal to the sum of the integral along the line $AB$, the integrals along the arcs $BCD$ and $HIA$, the integrals along the lines $DE$ and $GH$ on each side of the branch cut, and the integral along the arc $EFG$ around the branch point $s = -\Lambda$.

Asymptotic behavior of the Green's function \eqref{greenfs} as the radius $r$ of the arc $EFG$ around the branch point tends to zero implies that the integral along this arc vanishes. To verify this, first substitute $s = r e^{i \gamma} - \Lambda$, $\gamma \in \langle -\pi, \pi \rangle$. Then the integral along the arc $EFG$ has the following form:
\begin{equation}
\int\limits_{EFG} e^{st} G_{s}(x_{0}, x) \, ds = -i e^{-\Lambda t} \int\limits_{-\pi}^{\pi} e^{r e^{i \gamma}t} G_{r e^{i \gamma} - \Lambda}(x_{0}, x) r e^{i \gamma} \, d \gamma, \label{int:EFG}
\end{equation}
where
\begin{equation}
G_{r e^{i \gamma} - \Lambda}(x_{0}, x) = \frac{\mathfrak{m}(x)}{2 \alpha^{1-\frac{\alpha}{2}} \beta^{-\frac{\beta}{2}} B_{r e^{i \gamma} - \Lambda} \D_{r e^{i \gamma} - \Lambda}} \, f_{1}(x_{0} \wedge x, r e^{i \gamma} - \Lambda) f_{4}(x_{0} \vee x, r e^{i \gamma} - \Lambda). \label{BCDgreen}
\end{equation}
The expression under the integral vanishes as radius $r \to 0$ if $$\lim\limits_{r \to 0} G_{r e^{i \gamma} - \Lambda}(x_{0}, x) < \infty,$$ which can be verified by observing the power series expansion of $G_{r e^{i \gamma} - \Lambda}(x_{0}, x)$ around the point zero:
\begin{equation}
G_{r e^{i \gamma} - \Lambda}(x_{0}, x) = c_{\alpha, \beta, (x_{0} \wedge x), (x_{0} \vee x)}^{(1)} + c_{\alpha, \beta, (x_{0} \wedge x), (x_{0} \vee x)}^{(2)} \, e^{i \frac{\gamma}{2}} \sqrt{r} + c_{\alpha, \beta, (x_{0} \wedge x), (x_{0} \vee x)}^{(3)} \, e^{i \gamma} r + O(r^{2}), \label{BCDseries0}
\end{equation}
where quantities $c_{\alpha, \beta, (x_{0} \wedge x), (x_{0} \vee x)}^{(i)}$, $i \in \{1, 2, 3\}$, are independent of $r$. Expression \eqref{BCDseries0} implies that $\lim\limits_{r \to 0} G_{r e^{i \gamma} - \Lambda}(x_{0}, x) < \infty$. Therefore, the integral \eqref{int:EFG} vanishes as $r \to 0$.

Furthermore, asymptotic behavior of the Green's function \eqref{greenfs} as the radius $R$ of the arcs $BCD$ and $HIA$ tends to infinity implies that the integrals along these arcs also vanish. Since these integrals are treated similarly, we explain the procedure regarding the arc $BCD$. Substitution $s = R e^{i \gamma} - \Lambda$, where $\gamma \in \langle \pi/2, \pi \rangle$, transforms the integral along the arc $BCD$ in the following form:
\begin{equation}
\int\limits_{BCD} e^{st} G_{s}(x_{0}, x) \, ds = i e^{-\Lambda t} \int\limits_{\frac{\pi}{2}}^{\pi} R e^{R e^{i \gamma} t} G_{R e^{i \gamma} - \Lambda}(x_{0}, x) e^{i \gamma} \, d \gamma, \label{FSD_IntBCD}
\end{equation}
where
\begin{equation*}
G_{R e^{i \gamma} - \Lambda}(x_{0}, x) = \mathfrak{m}(x) \, \frac{\Gamma\left( -\frac{\beta}{4} + \sqrt{\frac{R (\beta-2)}{2\theta}} \, e^{i \frac{\gamma}{2}} \right) \Gamma\left( \frac{\alpha}{2} + \frac{\beta}{4} + \sqrt{\frac{R (\beta-2)}{2\theta}} \, e^{i \frac{\gamma}{2}} \right)}{\alpha^{1-\frac{\alpha}{2}} \beta^{-\frac{\beta}{2}} \Gamma\left( \frac{\alpha}{2} \right) \Gamma\left( 1 + 2\sqrt{\frac{R (\beta-2)}{2\theta}} \, e^{i \frac{\gamma}{2}} \right)} \times
\end{equation*}
\begin{equation}
\times f_{1}(x_{0} \wedge x, R e^{i \gamma} - \Lambda) f_{4}(x_{0} \vee x, R e^{i \gamma} - \Lambda). \label{FSD_BCDgreen}
\end{equation}
The expression under the integral on the right side of the expression \eqref{FSD_IntBCD} vanishes as the radius $R \to \infty$ if $$\lim\limits_{R \to \infty} R e^{R e^{i \gamma} t} G_{R e^{i \gamma} - \Lambda}(x_{0}, x) = 0.$$ According to the formulas for asymptotic behavior of the Gauss hypergeometric function $_{2}F_{1}(a, b; c; z)$ as two or more of its parameters tend to infinity (cf. Erdely \cite{Erdely}, page 77, expressions (16) and (17)) it follows that
\begin{equation*}
f_{1}(x_{0} \wedge x, R e^{i \gamma} - \Lambda) f_{4}(x_{0} \vee x, R e^{i \gamma} - \Lambda) =
\end{equation*}
\begin{equation*}
= \frac{\Gamma \left( \frac{\alpha}{2} \right)}{2^{\frac{\beta}{2} + 1} \, \sqrt{\pi}} (1 - k^{-1}_{\alpha, \beta, (x_{0} \wedge x)})^{\frac{1 - \alpha}{2}} (1 + k^{-1}_{\alpha, \beta, (x_{0} \wedge x)})^{\frac{\alpha + \beta - 2}{2}} (1 - k^{-1}_{\alpha, \beta, (x_{0} \vee x)})^{\frac{1 - \alpha}{2}} (1 + k^{-1}_{\alpha, \beta, (x_{0} \vee x)})^{\frac{\alpha + \beta - 2}{2}} \times
\end{equation*}
\begin{equation*}
\times \frac{\Gamma\left(1 + \frac{\beta}{4} + \sqrt{\frac{R (\beta-2)}{2\theta}} \, e^{i \frac{\gamma}{2}} \right) \Gamma\left(1  + 2 \sqrt{\frac{R (\beta-2)}{2\theta}} \, e^{i \frac{\gamma}{2}} \right)}{\Gamma\left( \frac{\alpha}{2} + \frac{\beta}{4} + \sqrt{\frac{R (\beta-2)}{2\theta}} \, e^{i \frac{\gamma}{2}} \right) \Gamma\left(1 - \frac{\alpha}{2} - \frac{\beta}{4} + \sqrt{\frac{R (\beta-2)}{2\theta}} \, e^{i \frac{\gamma}{2}} \right) \Gamma\left(\frac{\alpha}{2} - \frac{\beta}{4} + \sqrt{\frac{R (\beta-2)}{2\theta}} \, e^{i \frac{\gamma}{2}} \right)} \times
\end{equation*}
\begin{equation}
\times \sqrt{\frac{2 \theta}{R (\beta - 2)}} \, e^{-i \frac{\gamma}{2}} k^{\frac{\beta}{4} - \sqrt{\frac{R(\beta - 2)}{2 \theta}} \, e^{i \frac{\gamma}{2}}}_{\alpha, \beta, (x_{0} \vee x)} \left[ k^{\frac{\beta}{4} + \sqrt{\frac{R(\beta - 2)}{2 \theta}} \, e^{i \frac{\gamma}{2}}}_{\alpha, \beta, (x_{0} \wedge x)} + e^{i \pi \frac{\alpha - 1}{2}} k^{\frac{\beta}{4} - \sqrt{\frac{R(\beta-2)}{2 \theta}} \, e^{i \frac{\gamma}{2}}}_{\alpha, \beta, (x_{0} \wedge x)} \right] \left[1 + O\left( \frac{1}{R e^{i \gamma}} \right) \right], \label{FSD_BCDf1f4}
\end{equation}
for large values of $R$, where
\begin{equation*}
k_{\alpha, \beta, (x_{0} \wedge x)} = 1 + \frac{2 \alpha}{\beta}(x_{0} \wedge x) + \sqrt{\frac{4 \alpha}{\beta} (x_{0} \wedge x) \left(1 + (x_{0} \wedge x) \right)} > 1,
\end{equation*}
\begin{equation*}
k_{\alpha, \beta, (x_{0} \vee x)} = 1 + \frac{2 \alpha}{\beta}(x_{0} \vee x) + \sqrt{\frac{4 \alpha}{\beta} (x_{0} \vee x) \left(1 + (x_{0} \vee x) \right)} > 1
\end{equation*}
and $k_{\alpha, \beta, (x_{0} \wedge x)} < k_{\alpha, \beta, (x_{0} \vee x)}$. From expression \eqref{FSD_BCDf1f4} it follows that
\begin{equation*}
R e^{R e^{i \gamma} t} G_{R e^{i \gamma} - \Lambda}(x_{0}, x) =
\end{equation*}
\begin{equation*}
= \frac{\alpha^{\frac{\alpha}{2}-1} \beta^{\frac{\beta}{2}}}{2^{\frac{\beta}{2} + 1} \, \sqrt{\pi}} (1 - k^{-1}_{\alpha, \beta, (x_{0} \wedge x)})^{\frac{1 - \alpha}{2}} (1 + k^{-1}_{\alpha, \beta, (x_{0} \wedge x)})^{\frac{\alpha + \beta - 2}{2}} (1 - k^{-1}_{\alpha, \beta, (x_{0} \vee x)})^{\frac{1 - \alpha}{2}} (1 + k^{-1}_{\alpha, \beta, (x_{0} \vee x)})^{\frac{\alpha + \beta - 2}{2}} \times
\end{equation*}
\begin{equation*}
\times R e^{R e^{i \gamma} t} \, \left\{ \frac{\Gamma\left(1 + \frac{\beta}{4} + \sqrt{\frac{R (\beta-2)}{2\theta}} \, e^{i \frac{\gamma}{2}} \right)}{\Gamma\left(1 - \frac{\alpha}{2} - \frac{\beta}{4} + \sqrt{\frac{R (\beta-2)}{2\theta}} \, e^{i \frac{\gamma}{2}} \right)} \, \frac{\Gamma\left(-\frac{\beta}{4} \sqrt{\frac{R (\beta-2)}{2\theta}} \, e^{i \frac{\gamma}{2}} \right)}{\Gamma\left( \frac{\alpha}{2} - \frac{\beta}{4} + \sqrt{\frac{R (\beta-2)}{2\theta}} \, e^{i \frac{\gamma}{2}} \right)} \times \right.
\end{equation*}
\begin{equation*}
\times \sqrt{\frac{2 \theta}{R (\beta - 2)}} \, e^{-i \frac{\gamma}{2}} \left[ (k_{\alpha, \beta, (x_{0} \wedge x)} k_{\alpha, \beta, (x_{0} \vee x)})^{\frac{\beta}{4}} \left( \frac{k_{\alpha, \beta, (x_{0} \wedge x)}}{k_{\alpha, \beta, (x_{0} \vee x)}} \right)^{\sqrt{\frac{R(\beta - 2)}{2 \theta}} \cos{\left(\frac{\gamma}{2} \right)}} \, \left( \frac{k_{\alpha, \beta, (x_{0} \wedge x)}}{k_{\alpha, \beta, (x_{0} \vee x)}} \right)^{i \, \sqrt{\frac{R(\beta - 2)}{2 \theta}} \sin{\left(\frac{\gamma}{2} \right)}} \, + \right.
\end{equation*}
\begin{equation*}
\left. \left. + \, e^{i \pi \frac{\alpha}{2}} \, (k_{\alpha, \beta, (x_{0} \wedge x)} k_{\alpha, \beta, (x_{0} \vee x)})^{\frac{\beta}{4} - \sqrt{\frac{R(\beta - 2)}{2 \theta}} \cos{\left(\frac{\gamma}{2} \right)}} \, (k_{\alpha, \beta, (x_{0} \wedge x)} k_{\alpha, \beta, (x_{0} \vee x)})^{-i \, \sqrt{\frac{R(\beta - 2)}{2 \theta}} \sin{\left(\frac{\gamma}{2} \right)}} \right] \right\} \left[ 1 + O\left( \frac{1}{R} \right) \right]
\end{equation*}
for large values of $R$. Since $k_{\alpha, \beta, (x_{0} \wedge x)} < k_{\alpha, \beta, (x_{0} \vee x)}$ and since for $\pi/2 < \gamma < \pi$ we know that $\cos{(\frac{\gamma}{2})} > 0$ and $\cos{(\gamma)} < 0$ it follows that $R e^{R e^{i \gamma} t} G_{R e^{i \gamma} - \Lambda}(x_{0}, x)$ tends to zero as $R \to \infty$. Therefore, the integral \eqref{FSD_IntBCD} vanishes as $R \to \infty$. The same procedure is used for verifying that the integral along the arc $HIA$ also vanishes as $R \to \infty$. Problem observed in a recent paper by Carr, Linetsky and Mendoza \cite{Mendoza} is treated with the similar technique.

Finally, we conclude that the inversion formula \eqref{invLap} for the Green's function \eqref{greenfs} reduces to the sum of the integral along
the line $\mathrm{Re}(x) = c$, $c > 0$, and the integral of the jump across the branch cut, i.e.
\begin{equation}
\frac{1}{2 \pi i} \int\limits_{\mathcal{C}} e^{st} \, G_{s}(x_{0}, x) = \frac{1}{2 \pi i} \int\limits_{c - i \infty}^{c + i \infty} e^{st} \, G_{s}(x_{0}, x) \, ds + \frac{1}{2 \pi i} \int\limits_{-\infty}^{-\Lambda} e^{st} \, \left( G_{s}(f_{1}, f_{4}) - \overline{G}_{s}(f_{1}, f_{4}) \right) \, ds, \label{Bromwichcut}
\end{equation}
where $\overline{G}_{s}(f_{1}, f_{4})$ is the complex conjugate of the Green's function \eqref{greenfs}. From expressions \eqref{BromwichC} and \eqref{Bromwichcut} it follows:
\begin{equation}
\frac{1}{2 \pi i} \int\limits_{c - i \infty}^{c + i \infty} e^{st} \, G_{s}(x_{0}, x) \, ds = \sum_{n=0}^{\left\lfloor \frac{\beta}{4} \right\rfloor} {\Res}_{s=-\lambda_{n}} e^{st} G_{s}(x_{0}, x) + \frac{1}{2 \pi i} \int\limits_{-\infty}^{-\Lambda} e^{st} \, \left( \overline{G}_{s}(f_{1}, f_{4}) - G_{s}(f_{1}, f_{4}) \right) \, ds. \label{Bromwich}
\end{equation}
Now it remains to determine the explicit form of the expression $\overline{G}_{s}(f_{1}, f_{4}) - G_{s}(f_{1}, f_{4})$ for the jump across the branch cut. Since $\overline{f}_{1}(\cdot, s) = f_{1}(\cdot, s)$, $\overline{f}_{4}(\cdot, s) = f_{3}(\cdot, s)$ and $f_{1}(\cdot, s) = B f_{3}(\cdot, s) + A f_{4}(\cdot, s)$ (see Remark \ref{f1=Bf3+Af4}.), for $s = -\lambda$ the jump across the branch cut is given by expression
\begin{equation}
\overline{G}_{-\lambda}(x_{0}, x) - G_{-\lambda}(x_{0}, x) =
\label{greenjump}
\end{equation}
\begin{equation*}
= \mathfrak{m}(x) \, 2i k(\lambda) \,\left| \frac{\alpha^{\frac{\alpha - 2}{4}} \beta^{\frac{\beta}{4}} \, \Gamma\left( -\frac{\beta}{4} + ik(\lambda) \right) \Gamma\left( \frac{\alpha}{2} + \frac{\beta}{4} + ik(\lambda) \right)}{\Gamma\left( \frac{\alpha}{2} \right) \Gamma\left(1 + 2ik(\lambda) \right)} \right|^{2} \, f_{1}(x_{0}) f_{1}(x).
\end{equation*}

Substitution of the expressions \eqref{greenres} and \eqref{greenjump} in the relation \eqref{Bromwich} results in the spectral representation of the transition density of Fisher-Snedecor diffusion. However, Fisher-Snedecor polynomials and solutions of the Sturm-Liouville equation \eqref{sleq} for $\lambda > \Lambda$ are not normalized with respect to the Fisher-Snedecor density \eqref{fs}. Now the uniqueness of solutions of the Sturm-Liouville equation up to a constant factor ensures that instead of the solution $f_{1}(x, -\lambda)$ we can use the solution $\widetilde{f}_{1}(x, -\lambda) = \sqrt{\frac{2 \theta}{\beta - 2}} \, f_{1}(x, -\lambda)$. This procedure results in the expression for the spectral representation of transition density in terms of functions normalized with respect to the Fisher-Snedecor density, i.e. the result is exactly the expression $p(x; x_{0}, t) = p_{d}(x; x_{0}, t) + p_{c}(x; x_{0}, t)$, where $p_{d}(x; x_{0}, t)$ is given by \eqref{disspec} and $p_{c}(x; x_{0}, t)$ by \eqref{contspec}.
\end{proof}

\begin{rem}
\emph{In the previous theorem we supposed that $\alpha > 2$, $\alpha \notin \{2(m+1), \, m \in \mathbb{N}\}$. If this is not the case, i.e. if $\alpha > 2$ is an even integer, then in \eqref{greenf1f2} expression for the solution $f_{2}(x, s)$ discussed in the Remark \ref{f2EvenAlpha}. should be used. However, due to its complicated form it is hard to expect that a closed form expression for spectral representation of transition density is obtainable in such an elegant form as in the Theorem \ref{thm:SpecRep}.}
\end{rem}

\begin{rem}
\emph{The general approach to the spectral theory of one-dimensional diffusions is treated in a number of papers by Linetsky (see \cite{Linetsky1, Linetsky2, Linetsky3, Mendoza}). Beside the concise and informative introduction to basics of the spectral theory, he gives the procedure for calculation of the spectral representation of the transition density for one-dimensional diffusions (see \cite{Linetsky3, Mendoza}). We already used his approach for transition densities of reciprocal gamma and Student diffusion (see \cite{Leonenko1, Leonenko2, Avram3}). Since these processes, together with the Fisher-Snedecor diffusion, make a class of heavy-tailed Pearson diffusions, it is natural that spectral representations of their transition densities have similar forms. Namely, for all three heavy-tailed Pearson diffusions spectral representations of transition densities are given in terms of the finite sum of residues and an integral of the jump across the branch cut of the corresponding Green's function. However, the differences between mentioned spectral representations are generated by different types of the corresponding Sturm-Liouville equations. In particular, Sturm-Liouville equation related to reciprocal gamma diffusion can be transformed to the Whittaker equation, while the Sturm-Liouville equations for Student and Fisher-Snedecor diffusions can be transformed to the hypergeometric equations. According to the qualitative nature of the spectrum of the corresponding Sturm-Liouville operator $(-\mathcal{G})$, the continuous part of the spectral representation of transition density of reciprocal gamma diffusion is written in terms of Whittaker functions or hypergeometric function $\phantom{.}_{2}F_{0}(a, b; ; \cdot)$ (see \cite{Linetsky2, Linetsky3, Leonenko1, Wong}), while the continuous part of the spectral representation for Student and Fisher-Snedecor diffusions is written in terms of the hypergeometric function $\phantom{.}_{2}F_{1}(a, b; c; \cdot)$, see Appendix B.}
\end{rem}

\begin{rem} \label{rem:ortho}
\emph{Fisher-Snedecor polynomials $F_{n}(x)$, $n = 0, \ldots, \left\lfloor \beta/4 \right\rfloor$, and functions $f_{1}(x, -\lambda)$ for $\lambda > \Lambda$ belong to orthogonal subspaces of the Hilbert space $\mathcal{H} = L^{2}(\langle 0, \infty \rangle, \fs(\cdot))$, i.e.
$$F_{n}(x) \in \mathcal{H}_{pp}, \quad n \in \{ 0, 1, \ldots, \left\lfloor \beta/4 \right\rfloor\}, \quad \quad f_{1}(x, -\lambda) \in \mathcal{H}_{ac}, \quad \forall \lambda > \Lambda.$$ Here $\mathcal{H}_{pp}$ denotes the subspace of the Hilbert space $L^{2}(\langle 0, \infty \rangle, \fs(\cdot))$ containing functions having only the pure point spectral measure, while $\mathcal{H}_{ac}$ denotes the subspace of the Hilbert space $L^{2}(\langle 0, \infty \rangle, \fs(\cdot))$ containing functions having only the spectral measure which is absolutely continuous with respect to the Lebesgue measure (cf. Reed and Simon \cite{Reed} or Linetsky \cite{Linetsky1}). From here it follows that Fisher-Snedecor polynomials and functions $f_{1}(x, -\lambda)$ for $\lambda > \Lambda$ are orthogonal with respect to the invariant density of Fisher-Snedecor diffusion, i.e.
\begin{equation}
\int\limits_{0}^{\infty} F_{n}(x) f_{1}(x, -\lambda) \fs(x) \, dx = 0. \label{FSD_Fnf1Ortho}
\end{equation}
This orthogonality relation plays a crucial role in the statistical analysis of this diffusion process.}
\end{rem}

\begin{rem}
\emph{Two-dimensional density of the stationary diffusion process $X = \{X_{t}, \, t \geq 0\}$ that solves stochastic differential equation \eqref{fssde} and has marginal density function \eqref{fs} is given by the expression
\begin{equation}
p(x, y, t) = \frac{\partial^{2}}{\partial x \partial y} P(X_{s + t} \leq x, X_{s} \leq y) = \fs(y) \, p(x; y, t) = \fs(y) \, \left( p_{d}(x; y, t) + p_{c}(x; y, t) \right) \label{twodens}
\end{equation}
where $p_{d}(x; y, t)$ is given by \eqref{disspec} while $p_{c}(x; y, t)$ is given by \eqref{contspec}.}
\end{rem}

\medskip

\section{Parameter estimation for Fisher-Snedecor diffusion}

\subsection{Estimation of autocorrelation parameter $\protect\theta$}

The autocorrelation parameter $\theta > 0$ will be estimated under the assumption that parameters $\alpha$ and $\beta$ of the Fisher-Snedecor
diffusion are known. The ergodicity of the Fisher-Snedecor diffusion is ensured by the assumption $\alpha > 2$ while the existence of the first and the second moment of the invariant distribution is ensured by the assumption $\beta > 4$.

Let us consider the sample $X_{1}, \ldots, X_{n}$ from the ergodic stationary Fisher-Snedecor diffusion and the corresponding sample of paired observations $(X_{1}, X_{t+1}), (X_{2}, X_{t+2}), \ldots, (X_{n-t}, X_{n})$, where $t < n$. Empirical counterpart of the autocorrelation function is given by the Pearson's sample correlation coefficient
\begin{equation}
\widehat{\rho}_{n}(t) = \displaystyle\frac{\frac{1}{n-t} \sum_{i = 1}^{n - t} X_{i} X_{t + i} - \frac{1}{n-t} \sum_{i = 1}^{n - t} X_{i} \cdot \frac{1}{n-t} \sum_{i = 1}^{n - t} X_{t + i}}{\sqrt{\frac{1}{n-t} \sum_{i = 1}^{n - t} X_{i}^{2} - \left( \frac{1}{n-t} \sum_{i = 1}^{n - t} X_{i} \right)^{2}} \, \sqrt{\frac{1}{n-t} \sum_{i = 1}^{n - t} X_{t + i}^{2} - \left( \frac{1}{n-t} \sum_{i = 1}^{n - t} X_{t + i} \right)^{2}}},  \label{pearcorr}
\end{equation}
where the term in the numerator represents the empirical covariance of random variables $X_{s}$ and $X_{t + s}$, while the term in the denominator represents the product of the empirical standard deviations of random variables $X_{s}$ and $X_{t+s}$, respectively.

Expression \eqref{pearcorr} implies that there exists a continuous function $g \colon \mathbb{R}^{5} \setminus \{(x, y, z, u, v): x^{2} = z \text{ or } y^{2} = u \} \to \mathbb{R}$ defined by the expression
\begin{equation*}
g(x, y, z, u, v) = \frac{v - xy}{\sqrt{(z - x^{2})(u - y^{2})}},
\end{equation*}
such that
\begin{equation*}
\begin{tabular}{rcl}
$\widehat{\rho}_{n}(t)$ & = & $g \left(\frac{1}{n-t} \sum\limits_{i = 1}^{n - t} X_{i}, \frac{1}{n-t} \sum\limits_{i = 1}^{n - t} X_{t + i}, \frac{1}{n-t} \sum\limits_{i = 1}^{n - t} X_{i}^{2}, \frac{1}{n-t} \sum\limits_{i = 1}^{n - t} X_{t + i}^{2}, \frac{1}{n-t} \sum\limits_{i = 1}^{n - t} X_{i} X_{t + i} \right)$.
\end{tabular}
\end{equation*}
Since $X_{1}, \ldots, X_{n}$ is the sample from the ergodic stationary Markov process, according to ergodic theorem for stationary sequences (see Karlin and Taylor \cite{Karlin1}, Theorems 5.4. and 5.6.), it follows that
\begin{equation}
\left( \frac{1}{n-t} \sum\limits_{i = 1}^{n - t} X_{i}, \frac{1}{n-t} \sum\limits_{i = 1}^{n - t} X_{t + i}, \frac{1}{n-t} \sum\limits_{i = 1}^{n
- t} X_{i}^{2}, \frac{1}{n-t} \sum\limits_{i = 1}^{n - t} X_{t + i}^{2}, \frac{1}{n-t} \sum\limits_{i = 1}^{n - t} X_{i} X_{t + i} \right) \label{moments_est}
\end{equation}
is the $P$-consistent estimator of $E\left[ (X_{0}, X_{t}, X_{0}^{2}, X_{t}^{2}, X_{0}X_{t}) \right]$.

From expression \eqref{acf} for the autocorrelation function $\rho(t)$ of the Fisher-Snedecor diffusion we see that it takes only positive values. However, according to Hassani (cf. \cite{Hassani}, Theorem 2.1.) the sum of the sample autocorrelation function at lag $t \geq 1$ is always $-1/2$ for any stationary time series with arbitrary length $n \geq 2$, i.e. $$S_{acf} = \sum\limits_{t=1}^{n-1} \widehat{\rho}_{n}(t) = -\frac{1}{2}.$$ Therefore, instead of the Pearson's sample correlation coefficient \eqref{pearcorr} we observe here its absolute value $\vert \widehat{\rho}_{n}(t) \vert$. Since $\vert \widehat{\rho}_{n}(t) \vert$ is the continuous transformation of the $P$-consistent estimator \eqref{moments_est}, according to continuity mapping theorem it is the $P$-consistent estimator of the autocorrelation function \eqref{acf} for any fixed $t > 0$, i.e.
\begin{equation}
\vert \widehat{\rho}_{n}(t) \vert \overset{P}{\to} \rho(t), \quad n \to \infty, \quad \text{ for any fixed } t > 0.  \label{rhohat}
\end{equation}
Therefore, here we observe the equation $\vert \widehat{\rho}_{n}(t) \vert = \rho(t)$, where $\rho(t) = e^{-\theta t}$. Solving this estimation equation in terms of the unknown parameter $\theta$ results in the estimator $\widehat{\theta}$ of the autocorrelation parameter $\theta$, i.e.
\begin{equation}
\widehat{\theta} = -\displaystyle\frac{1}{t} \, \log{\vert \widehat{\rho}_{n}(t) \vert}.
\label{thetaest}
\end{equation}
The continuity of the logarithmic function and the relation \eqref{rhohat} imply that
\begin{equation*}
-\displaystyle\frac{1}{t} \, \log{\vert \widehat{\rho}_{n}(t) \vert} \overset{P}{\to} - \displaystyle\frac{1}{t} \, \log{\rho(t)}, \quad n \to \infty, \quad \text{ for any fixed } t > 0,
\end{equation*}
which means that $\widehat{\theta}$ is the $P$-consistent estimator of the autocorrelation parameter $\theta$, i.e.
\begin{equation}
\widehat{\theta} \overset{P}{\to} \theta, \quad n \to \infty.  \label{thetacons}
\end{equation}

\medskip

\subsection{Estimation of parameters $\protect\alpha$ and $\protect\beta$}

In this section we present method of moments estimators of unknown parameters $\alpha$ and $\beta$ of invariant distribution of ergodic Fisher-Snedecor diffusion and analyze their asymptotic properties under the assumption that $\theta$ is the known value of the autocorrelation parameter. Estimation of the unknown parameter $(\alpha, \beta, \theta)$ by the method based on the martingale estimation functions (cf. Forman and S\o rensen \cite{Forman}) is briefly discussed in Section \ref{MartingaleEstimation}. However, in contrast to the method which will be developed here, their method for calculation of limiting covariance matrix in asymptotic normality framework seems to be more complicated.

A new method for calculation of $E[X_{s+t}^{n}X_{s}^{m}]$ based on the orthogonality of the normalized Fisher-Snedecor polynomials with respect to the density \eqref{fs} will be introduced. The existence of Fisher-Snedecor polynomials $F_{1}(x)$ and $F_{2}(x)$ needed for mentioned calculations is ensured by the assumption that $\beta > 8$.

Equating the first and the second theoretical moments with the corresponding empirical counterparts
\begin{equation}
\overline{m}_{1} = \displaystyle\frac{1}{n} \sum_{t=1}^{n}X_{t} \quad \quad \textrm{and} \quad \quad \overline{m}_{2} = \displaystyle\frac{1}{n}
\sum_{t=1}^{n}X^{2}_{t} \label{emp12}
\end{equation}
and solving the resulting system of equations in terms of the unknown parameters $\alpha$ and $\beta$ results in the method of moments estimators for these parameters. They are given by following expressions:
\begin{equation}
\widehat{\alpha} = \frac{2 \, \overline{m}_{1}^{\, 2}}{\overline{m}_{2}(2 - \overline{m}_{1}) - \overline{m}_{1}^{\, 2}},  \label{alphaest}
\end{equation}
\begin{equation}
\widehat{\beta} = \frac{2\overline{m}_{1}}{\overline{m}_{1}-1}.  \label{betaest}
\end{equation}

The statement of the following proposition follows from the general result for the spectral representation of transition density for one-dimensional diffusions. This special case for Fisher-Snedecor diffusion is treated here in details to demonstrate the use of a closed form expression for spectral representation of transition density (see Theorem \ref{thm:SpecRep}.). Analogue results for reciprocal gamma and Student diffusions are presented in \cite{Leonenko1, Leonenko2}. \newline

\begin{prop} \label{ExpFmFn}
Fix $i$ and $j$ such that $i, j \in \{ 1, \ldots, \left\lfloor \beta/4 \right\rfloor \}$ and suppose that $\beta > 2(i + j)$, $\alpha > 2$, $\alpha \notin \{2(m+1), \, m \in \mathbb{N}\}$ and $\theta > 0$. If $(X_{s+t}, X_{s})$ is a two-dimensional random vector with the density function \eqref{twodens} and $F_{n}(x)$, $n \in \{ 1, \ldots, \left\lfloor \beta/4 \right\rfloor \}$, are orthonormal Fisher-Snedecor polynomials given by Rodrigues formula \eqref{rodfsnorm}, then
\begin{equation*}
E[F_{i}(X_{s + t}) \, F_{j}(X_{s})] = e^{-\lambda_{j} t} \, \delta_{ij},
\end{equation*}
where $\lambda_{j}$ is the eigenvalue given by \eqref{diseigen} and $\delta_{ij}$ the standard Kronecker symbol.
\end{prop}

\begin{proof}
The two-dimensional density function $p(x, y, t)$ given by \eqref{twodens} could be written in the following form:
\begin{equation} p(x, y, t) = \sum\limits_{n=0}^{\left\lfloor \frac{\beta}{4} \right\rfloor} e^{-t \lambda_{n}} \Phi_{n}(x) \Phi_{n}(y) + \displaystyle\frac{1}{\pi} \, \int\limits_{\frac{\theta \beta^{2}}{8(\beta-2)}}^{\infty} e^{-\lambda t} \, \nu^{2}(\lambda) \Psi(x, \lambda) \Psi(y, \lambda) \, d\lambda, \label{TwoDens}
\end{equation} where $$\nu^{2}(\lambda) = k(\lambda) \, \left| \frac{B^{\frac{1}{2}} \left( \frac{\alpha}{2}, \frac{\beta}{2} \right) \, \Gamma\left( -\frac{\beta}{4} + ik(\lambda) \right) \Gamma\left( \frac{\alpha}{2} + \frac{\beta}{4} + ik(\lambda) \right)}{\Gamma \left( \frac{\alpha}{2} \right) \, \Gamma\left( 1 + 2ik(\lambda) \right)} \right|^{2},$$ and
where
$$\begin{tabular}{rcl}
    $\Phi_{n}(\cdot)$      & = & $\fs(\cdot) \, F_{n}(\cdot)$, \\
    $\Psi(\cdot, \lambda)$ & = & $\fs(\cdot) \, f_{1}(\cdot) \, \, \, = \, \, \, \fs(\cdot) \, f_{1}(\cdot, -\lambda)$.
\end{tabular}$$
Here $F_{n}(\cdot)$, $n \in \{ 1, \ldots, \left\lfloor \beta/4 \right\rfloor \}$, are Fisher-Snedecor polynomials and $\sqrt{\nu^{2}(\lambda)} \, f_{1}(\cdot, -\lambda)$ are solutions of the Sturm-Liouville equation \eqref{sleq} for $\lambda > \Lambda$, both normalized with respect to the invariant density of the Fisher-Snedecor diffusion.

Orthogonality of Fisher-Snedecor polynomials $F_{n}(\cdot)$, $n \in \{ 1, \ldots, \left\lfloor \beta/4 \right\rfloor \}$, and functions $f_{1}(\cdot, -\lambda)$, $\lambda > \Lambda$, (see the Remark \ref{rem:ortho}) together with the representation \eqref{TwoDens} of two-dimensional density \eqref{twodens} imply the following:
$$E[F_{i}(X_{s + t}) \, F_{j}(X_{s})] = \displaystyle\int\limits_{0}^{\infty}\displaystyle\int\limits_{0}^{\infty} F_{i}(x) F_{j}(y) p(x, y, t) \, dx \, dy \, = $$ $$ = \, \sum\limits_{n=0}^{\left\lfloor \frac{\beta}{4} \right\rfloor} e^{-\lambda_{n} t} \left(\displaystyle\int\limits_{0}^{\infty} F_{i}(x) F_{n}(x) \fs(x) \, dx \right) \left(\displaystyle\int\limits_{0}^{\infty} F_{j}(y) F_{n}(y) \fs(y) \, dy \right) \, +$$ $$ + \, \displaystyle\frac{1}{\pi} \, \displaystyle\int\limits_{\frac{\theta \beta^{2}}{8(\beta-2)}}^{\infty} \nu^{2}(\lambda) e^{-\lambda t} \, \left( \displaystyle\int\limits_{0}^{\infty} F_{i}(x) f_{1}(x, \lambda) \fs(x) \, dx \right)\left( \displaystyle\int\limits_{0}^{\infty} F_{j}(y) f_{1}(y, \lambda) \fs(y) \, dy \right) \, d\lambda \, =$$ $$= \, e^{- \lambda_{j} t} \delta_{ij} = \exp{\left\{-\theta j \frac{\beta - 2j}{\beta - 2} t \right\}} \, \delta_{ij},$$
where, according to the Remark \ref{rem:ortho}, the term related to the continuous part of the spectral representation \eqref{d+c} disappears. In order to apply the Fubini's theorem in the above computation we used the property that the solution $f_{1}(z, -\lambda) = \,_{2}F_{1}(a, b; c; -z)$ of the hypergeometric equation can be written in the form $$f_{1}(z, -\lambda) = C_{1} (-z)^{-a} +C_{2}(-z)^{-b} + O((-z)^{-a-1}) + O((-z)^{-b-1}),$$ where $C_{1}$ and $C_{2}$ are constants and $(a-b)$ is not an integer.
\end{proof}

Expressions \eqref{alphaest} and \eqref{betaest} imply that the estimators $\widehat{\alpha}$ and $\widehat{\beta}$ are continuous transformations of the empirical counterparts of the first and the second moment given by \eqref{emp12}. In particular, $\widehat{\alpha} = g_{1}(\overline{m}_{1}, \overline{m}_{2})$ and $\widehat{\beta} = g_{2}(\overline{m}_{1}, \overline{m}_{2})$. Therefore, asymptotic properties of the estimators $\widehat{\alpha}$ and $\widehat{\beta}$ are directly implied by asymptotic properties of the empirical counterpsrts $\overline{m}_{1}$ and $\overline{m}_{2}$, which will be presented now. Since Fisher-Snedecor diffusion has the $\alpha$-mixing property with the exponentially decaying rate (see Remark \ref{rem:mixing}), in further analysis we can use the following central limit theorem for $\alpha$-mixing sequences (see Hall and Heyde \cite{Hall} or Billingsley \cite{Billingsley2}; for multidimensional version see Genon-Catalot et al. \cite{Genon-Catalot}).

\begin{thm} \label{MFCLT}
\emph{Let $0 < \delta \leq \infty$. Suppose that $X = (X_{n}, \, n \in \mathbb{N})$ is a stationary and $\alpha$-mixing sequence of random variables such that $E[|X_{0}|] < \infty$ and let $f_{1}, \ldots, f_{m} \colon R \to R$, $m \in \mathbb{N}$, be Borel functions such that
\begin{itemize}
  \item $E[|f_{k}(X_{0})|^{2+\delta}] < \infty$ if $0 < \delta < \infty$,
  \item $|f_{k}(X_{0})| \leq c < \infty$ if $\delta = \infty$,
  \item $\sum\limits_{n=1}^{\infty}(\alpha_{f_{k}(X)}(n))^{\frac{\delta}{2 + \delta}} < \infty$.
\end{itemize}
Let $\mathbb{S}_{n} = \displaystyle\sum\limits_{k=1}^{n} (f_{1}(X_{k}), \ldots, f_{m}(X_{k}))$. Then
\begin{equation*}
\sigma^{2}_{ij} = \mathop{\mathrm{Cov}}(f_{i}(X_{0}), f_{j}(X_{0})) + \sum\limits_{k=1}^{\infty} \mathop{\mathrm{Cov}}(f_{i}(X_{0}),
f_{j}(X_{k})) + \sum\limits_{k=1}^{\infty} \mathop{\mathrm{Cov}}(f_{i}(X_{k}), f_{j}(X_{0}))
\end{equation*}
exists $\forall i, j \in \{1, \ldots, m\}$. If the matrix $\Sigma = [\sigma^{2}_{ij}]_{(m \times m)}$ is positive semidefinite, then
\begin{equation*}
\frac{\mathbb{S}_{n} - E[\mathbb{S}_{n}]}{\sqrt{n}} \overset{d}{\to} \mathcal{N}(0, \Sigma), \quad n \to \infty.
\end{equation*}}
\end{thm}

\begin{rem}
\emph{As an alternative, the central limit theorem can be formulated in terms of the $\beta$-mixing property of stationary diffusion. Beside the approach of Genon-Catalot et al. \cite{Genon-Catalot}, $\beta$-mixing property of Fisher-Snedecor diffusion can be proved by the methodology proposed by Abourashchi and Veretennikov. In recent papers (see \cite{Abourashchi1, Abourashchi2}) these authors proved that reciprocal gamma and Student diffusion models posses the exponential bounds for the $\beta$-mixing coefficient both in time and in space. The latter result is particulary valuable since it implies that functional central limit theorem can be also formulated for non-stationary reciprocal gamma and Student diffusion models.}
\end{rem}

\begin{rem}
\emph{An alternative approach to the central limit theorem for stationary diffusions is based on the theory of associated random systems. This approach is developed by Bulinski and Shashkin (see \cite{Bulinski}, Corollary 2.30. and Theorem 2.29.) and takes its origin in the paper by Barbato (see \cite{Barbato}, Lemma 8) where he proved that the stationary diffusions satisfy association conditions.}
\end{rem}

The ergodic theorem for stationary sequences (see Karlin and Taylor \cite{Karlin1} or the Lemma 20.3. from Billingsley \cite{Billingsley}) imply that
$\left( \overline{m}_{1}, \overline{m}_{2} \right)$ is the $P$-consistent estimator of $E[(X_{t}, X_{t}^{2})]$, while Theorem \ref{MFCLT} implies that $\left( \overline{m}_{1}, \overline{m}_{2} \right)$ is asymptotically normal, i.e.
\begin{equation}
\sqrt{n} \left( \overline{m}_{1} - E[\overline{m}_{1}], \overline{m}_{2} - E[\overline{m}_{2}] \right) \overset{d}{\to} \mathcal{N}(\mathbf{0}, \mathbf{\Sigma}),  \label{ANm1m2}
\end{equation}
where
\begin{equation*}
E[\overline{m}_{1}] = \frac{\beta}{\beta - 2}, \quad E[\overline{m}_{2}] = \frac{\beta^{2}(\alpha + 2)}{\alpha(\beta - 2)(\beta-4)}.
\end{equation*}

To obtain the explicit form of the $(2 \times 2)$ covariance matrix $\mathbf{\Sigma}$ we must first calculate the expectations $E[X_{s + t} X_{s}]$, $E[X_{s + t} X_{s}^{2}]$ and $E[X_{s + t}^{2} X_{s}^{2}]$, where $X_{s}$ and $X_{s+t}$ are random variables from the Fisher-Snedecor diffusion. These moments can be calculated by using the known expectations of the products of Fisher-Snedecor transformations of random variables $X_{s}$ and $X_{s+t}$. In particular, we observe random variables $F_{n}(X_{s})$ and $F_{n}(X_{s+t})$, where $n \in \{1, 2\}$ and $F_{n}(\cdot)$ is normalized Fisher-Snedecor polynomial \eqref{rodfsnorm}. According to the Proposition \ref{ExpFmFn}, expectations of their products are given by
\begin{equation*}
\begin{tabular}{l}
$E \left[ F_{1}(X_{s + t}) F_{1}(X_{s}) \right] = e^{-\theta t}$, \\
$E \left[ F_{1}(X_{s + t}) F_{2}(X_{s}) \right] = 0$, \\
$E \left[ F_{2}(X_{s + t}) F_{2}(X_{s}) \right] = e^{-2 \theta \, \frac{\beta - 4}{\beta - 2} \, t}$.
\end{tabular}
\end{equation*}

On the other hand, the same expectations can be calculated regarding the explicit definitions of the Fisher-Snedecor polynomials $F_{1}(x)$ and $F_{2}(x)$ (see Appendix A, section A.2.) and written in terms of the unknown expectations $E[X_{s + t} X_{s}]$, $E[X_{s + t} X_{s}^{2}]$ and $E[X_{s + t}^{2} X_{s}^{2}]$. In particular, \\

$E \left[ F_{1}(X_{s + t}) F_{1}(X_{s}) \right] = \displaystyle\frac{\alpha(\beta - 4)}{2 \beta^{2} (\alpha + \beta - 2)} \left[(\beta - 2)^{2}
E[X_{s+t} X_{s}] - \beta^{2} \right]$, \\

$E \left[ F_{1}(X_{s + t}) F_{2}(X_{s}) \right] = \displaystyle\frac{1}{4 \alpha \beta^{3} (\alpha + \beta - 2)} \, \sqrt{\displaystyle\frac{%
(\beta-2)(\beta-4)(\beta-8)}{(\alpha+2)(\alpha+\beta-4)}} \left[\alpha^{3}(\beta-2)(\beta-4)(\beta-6) \times \right.$ \\ \\
$\hspace*{4.2cm} \left. \times E[X_{s+t}X_{s}^{2}] - \displaystyle\frac{4 \alpha \beta^{3} (\alpha+2) (\alpha + \beta - 2)}{\beta - 2} \, e^{-\theta t} - \displaystyle\frac{\alpha^{2} \beta^{3} (\alpha+2) (\beta - 6)}{\beta - 2} \right]$, \\

$E \left[ F_{2}(X_{s + t}) F_{2}(X_{s}) \right] = \displaystyle\frac{\alpha (\beta - 2) (\beta - 8)}{8 \beta^{4} (\alpha + 2) (\alpha + \beta -
2) (\alpha + \beta - 4)} \left[ \alpha^{2}(\beta-4)^{2}(\beta-6)^{2} E[X_{s+t}^{2} X_{s}^{2}] - \right.$ \\ \\
$\hspace*{4.2cm} \left. - \displaystyle\frac{8 \beta^{4} (\alpha + 2)^{2} (\beta - 4)
(\alpha + \beta - 2)}{\alpha (\beta - 2)^{2}} \, e^{-\theta t} - \displaystyle\frac{\beta^{4} (\alpha + 2)^{2} (\beta - 6)^{2}}{(\beta -
2)^{2}} \right].$ \\

This approach results in the system of equations whose unknowns are expectations $E[X_{s + t} X_{s}]$, $E[X_{s + t} X_{s}^{2}]$ and $E[X_{s +
t}^{2} X_{s}^{2}]$. Explicitly, these expectations are given by the following expressions: \\

$E \left[X_{s + t} \, X_{s}\right] = \displaystyle\frac{1}{(\beta - 2)^{2}} \, \left[\displaystyle\frac{2 \beta^{2} (\alpha + \beta - 2)}{\alpha(\beta - 4)} \, e^{-\theta t} + \beta^{2} \right]$, \\

$E \left[X_{s + t} \, X_{s}^{2}\right] = \displaystyle\frac{\beta^{3} (\alpha + 2)}{\alpha^{2} (\beta - 2)^{2} (\beta - 4) (\beta - 6)} \, \left[4(\alpha + \beta - 2) e^{-\theta t} + \alpha (\beta - 6) \right] \, \, \, = \, \, \, E \left[X_{s + t}^{2} \, X_{s}\right]$, \\

$E \left[X_{s + t}^{2} \, X_{s}^{2}\right] = \displaystyle\frac{8 \beta^{4} (\alpha + 2) (\alpha + \beta - 2) (\alpha + \beta - 4)}{\alpha^{3}(\beta - 2) (\beta - 4)^{2} (\beta - 6)^{2} (\beta - 8)} \, e^{-2\theta \, \frac{\beta - 4}{\beta - 2} \, t} +$ \\ \\
$\hspace*{3.2cm} + \displaystyle\frac{8 \beta^{4} (\alpha + 2)^{2} (\alpha + \beta - 2)}{\alpha^{3} (\beta - 2)^{2} (\beta - 4) (\beta - 6)^{2}} \, e^{-\theta t} + \displaystyle\frac{\beta^{4} (\alpha + 2)^{2}}{\alpha^{2} (\beta - 2)^{2} (\beta - 4)^{2}}$. \\

Elements $\sigma^{2}_{ij}$, $i, j \in \{1, 2\}$, of the symmetric covariance matrix $\mathbf{\Sigma}$ now can be determined according to the formula for covariances given in the Theorem \ref{MFCLT}. In particular, \\

$\sigma_{11}^{2} = \displaystyle\frac{2 \beta^{2} (\alpha + \beta - 2)}{\alpha (\beta - 2)^{2} (\beta - 4)} \, \displaystyle\frac{e^{\theta} + 1}{%
e^{\theta} - 1}$, \\

$\sigma_{12}^{2} = \displaystyle\frac{4 \beta^{2} (\alpha + 2) (\alpha + \beta - 2)}{\alpha^{2} (\beta - 2)^{2} (\beta - 4) (\beta - 6)} \, %
\displaystyle\frac{e^{\theta} + 1}{e^{\theta} - 1}  = \sigma_{21}^{2}$, \\

$\sigma_{22}^{2} = \displaystyle\frac{\beta^{4} (\alpha + 2) \left[ (\alpha + 4) (\alpha + 6) (\beta - 2) (\beta - 4) - \alpha (\alpha + 2)
(\beta - 6) (\beta - 8) \right]}{\alpha^3 (\beta - 2)^2 (\beta - 4)^2 (\beta - 6) (\beta - 8)} +$ \\ \\
$\hspace*{1.5cm} + \displaystyle\frac{16 \beta^{4} (\alpha + 2) (\alpha + \beta - 2)(\alpha + \beta - 4)}{\alpha^{3} (\beta - 2) (\beta- 4)^{2} (\beta - 6)^{2} (\beta - 8)} \, \frac{1}{e^{2 \theta \frac{\beta - 4}{\beta - 2}} - 1} + \displaystyle\frac{16 \beta^{4} (\alpha + 2)^{2} (\alpha + \beta - 2)}{\alpha^{3} (\beta - 2)^{2} (\beta- 4) (\beta - 6)^{2}} \, \frac{1}{e^{\theta} - 1}.$ \\

Since $\alpha > 2$, $\alpha \notin \{2(m+1), \, m \in \mathbb{N}\}$, $\beta > 8$ and $\theta > 0$, the covariance matrix $\mathbf{\Sigma}$ is positive definite. This concludes the analysis of asymptotic properties of the estimator $(\overline{m}_{1}, \overline{m}_{2})$. Now we can start with the analysis of  asymptotic properties of the estimator $\left( \widehat{\alpha}, \widehat{\beta} \right)$, given in the following theorem. \\

\begin{thm} \label{Asymptotic}
Let $\{X_{t}, \, t \geq 0\}$ be the Fisher-Snedecor diffusion with the unknown parameter $(\alpha, \beta)$, where $\alpha > 2$, $\alpha \notin \{2(m+1), \, m \in \mathbb{N}\}$, and $\beta > 8$.

\begin{itemize}
\item[(i)] $(\widehat{\alpha}, \widehat{\beta}) \overset{P}{\to} (\alpha, \beta)$, i.e. $(\widehat{\alpha}, \widehat{\beta})$ is the $P$-consistent estimator of the unknown parameter $(\alpha, \beta)$.

\item[(ii)] $\sqrt{n} (\widehat{\alpha} - \alpha, \widehat{\beta} - \beta) \overset{d}{\to} \mathcal{N}(\mathbf{0}, \mathbf{\Sigma(\alpha, \beta, \theta)})$, where $\mathbf{\Sigma(\alpha, \beta, \theta)}$ is the positive definite covariance matrix.

\item[(iii)] $\sqrt{n} \left[ \mathbf{\Sigma(\widehat{\alpha}, \widehat{\beta}, \theta)} \right]^{-1/2} (\widehat{\alpha} - \alpha, \widehat{\beta} - \beta) \overset{d}{\to} \mathcal{N}(\mathbf{0}, \mathbf{I}) $, where $\theta$ is the known value of the autocorrelation parameter.

\item[(iv)] $\sqrt{n} \left[ \mathbf{\Sigma(\widehat{\alpha}, \widehat{\beta}, \widehat{\theta})} \right]^{-1/2} (\widehat{\alpha} - \alpha, \widehat{\beta} - \beta) \overset{d}{\to} \mathcal{N}(\mathbf{0}, \mathbf{I}) $, where $\widehat{\theta}$ is the $P$-consistent estimator \eqref{thetaest} of the autocorrelation parameter $\theta$.
\end{itemize}
\end{thm}

\begin{proof}
\begin{itemize}
\item[]
\item[(i)] Estimators $\overline{m}_{1}$ and $\overline{m}_{2}$ given by \eqref{emp12} are the $P$-consistent estimators of the first and the second moment, respectively. Since estimators $\widehat{\alpha}$ and $\widehat{\beta}$ are their continuous transformations, i.e. $$(\widehat{\alpha}, \widehat{\beta}) = \left( g_{1}(\overline{m}_{1}, \overline{m}_{2}), g_{2}(\overline{m}_{1}, \overline{m}_{2}) \right),$$ according to the Theorem 1.7.1. from Serfling \cite{Serfling}, it follows that
    $$\widehat{\alpha} = g_{1}(\overline{m}_{1}, \overline{m}_{2}) \overset{P}{\to} g_{1} \left( E[\overline{m}_{1}, \overline{m}_{2}] \right)= \alpha, \quad n \to \infty,$$
    $$\widehat{\beta} = g_{2}(\overline{m}_{1}, \overline{m}_{2}) \overset{P}{\to} g_{2} \left( E[\overline{m}_{1}, \overline{m}_{2}] \right) = \beta, \quad n \to \infty.$$
    From here it follows that $$(\widehat{\alpha}, \widehat{\beta}) \overset{P}{\to} (\alpha, \beta), \quad n \to \infty,$$ i.e. $(\widehat{\alpha}, \widehat{\beta})$ is the $P$-consistent estimator of the unknown parameter $(\alpha, \beta)$.

\item[(ii)] Since $\widehat{\alpha}$ and $\widehat{\beta}$ are the continuous transformations of $\overline{m}_{1}$ and $\overline{m}_{2}$, according to the multivariate delta method (see Serfling \cite{Serfling}, Theorem 3.3.A.), the estimator $(\widehat{\alpha}, \widehat{\beta})$ is asymptotically normal, i.e.
    \begin{equation}
    \sqrt{n} (\widehat{\alpha} - \alpha, \widehat{\beta} - \beta) \overset{d}{\to} \mathcal{N} \left( \mathbf{0}, D \mathbf{\Sigma} D^{\tau} \right), \label{ANab}
    \end{equation}
    where $\alpha = g_{1}\left( E[\overline{m}_{1}, \overline{m}_{2}] \right)$, $\beta = g_{2}\left( E[\overline{m}_{1}, \overline{m}_{2}] \right)$ and $(2 \times 2)$ matrix $D = \left[ \displaystyle\frac{\partial g_{i}}{\partial x_{j}} \right]_{\mathbf{x} = E[\overline{m}_{1}, \overline{m}_{2}]}$, $i \in \{1, 2\}$, has the following form:
    $$D = \left[ \begin{array}{cc}
    \displaystyle\frac{\alpha (\alpha + 2) (\beta - 2) (3\beta - 8)}{2 \beta (\beta - 4)} & -\displaystyle\frac{\alpha^{2} (\beta - 2) (\beta - 4)}{2 \beta^{2}} \\ \\
    -\displaystyle\frac{(\beta - 2)^{2}}{2} & 0 \\
    \end{array} \right]. $$

    Explicit forms of the elements $\sigma^{2}_{ij}(\alpha, \beta, \theta)$ of the covariance matrix $D \mathbf{\Sigma} D^{\tau} = \mathbf{\Sigma(\alpha, \beta, \theta)}$ are given by the following expressions: \\

    $\sigma^{2}_{11}(\alpha, \beta, \theta) = \displaystyle\frac{\alpha (\alpha + 2) (\beta - 2) (\alpha + \beta - 2)}{2 \beta (\beta - 4)^{3} (\beta - 6)^{2}} \times$ \\ \\
    $\hspace*{1.5cm} \times \left[ \displaystyle\frac{(\alpha + 2)(-3072 + \beta (6528 + \beta (-4736 + \beta (1548 + \beta (-232 + 13 \beta)))))}{\beta - 2} \, \coth{\left( \displaystyle\frac{\theta}{2}\right)} + \right.$ \\ \\
    $\hspace*{1.5cm} \left. + \displaystyle\frac{4 \beta (\beta - 4)^{3} (\alpha + \beta - 4)}{\beta - 8} \, \coth{\left( \displaystyle\frac{\theta (\beta - 4)}{\beta - 2}\right)} \right],$ \\

    $\sigma^{2}_{12}(\alpha, \beta, \theta) = -\displaystyle\frac{(\alpha + 2) (\beta - 2) (\alpha + \beta - 2) \left[\beta (\beta - 4) (3 \beta - 16) - 32 \right]}{2 (\beta - 4)^{2} (\beta - 6)} \, \coth{\left( \displaystyle\frac{\theta}{2} \right)} = \sigma^{2}_{21}(\alpha, \beta, \theta),$ \\

    $\sigma^{2}_{22}(\alpha, \beta, \theta) = \displaystyle\frac{\beta^{2} (\beta - 2)^{2} (\alpha + \beta - 2)}{2 \alpha (\beta - 4)} \, \coth{\left(  \displaystyle\frac{\theta}{2}\right)}$, \quad where  $\coth{\left( \displaystyle\frac{\theta}{2}\right)} = \displaystyle\frac{e^{\theta} + 1}{e^{\theta} - 1}$. \\

\item[(iii)] Part (ii) of this theorem and the positive definiteness of the covariance matrix $\mathbf{\Sigma}(\alpha, \beta, \theta)$ imply that \begin{equation}
    \left[ \mathbf{\Sigma}(\alpha, \beta, \theta) \right]^{-1/2} \, \sqrt{n} \, (\widehat{\alpha} - \alpha, \widehat{\beta} - \beta) \overset{d}{\to} \mathcal{N}(\mathbf{0}, \mathbf{I}), \label{ANab2}
    \end{equation} where $\mathbf{I}$ is the $(2 \times 2)$ identity matrix. Since $(\widehat{\alpha}, \widehat{\beta})$ is the $P$-consistent estimator of the parameter $(\alpha, \beta)$, it follows that
    \begin{equation}
    \left[ \mathbf{\Sigma}(\widehat{\alpha}, \widehat{\beta}, \theta) \right]^{-1/2} \, \left[ \mathbf{\Sigma}(\alpha, \beta, \theta) \right]^{1/2} \overset{P}{\to} \mathbf{I}. \label{I}
    \end{equation}
    If we multiply $$\left[ \mathbf{\Sigma}(\alpha, \beta, \theta) \right]^{-1/2} \, \sqrt{n} \, (\widehat{\alpha} - \alpha, \widehat{\beta} - \beta)$$ from the left by the matrix $\left[ \mathbf{\Sigma}(\widehat{\alpha}, \widehat{\beta}, \theta) \right]^{-1/2} \, \left[ \mathbf{\Sigma}(\alpha, \beta, \theta) \right]^{1/2}$, then \eqref{ANab2}, \eqref{I} and the Slutsky Lemma (see Serfling \cite{Serfling}, theorem 1.5.4.1.) imply that
    \begin{equation*}
    \left[ \mathbf{\Sigma}(\widehat{\alpha}, \widehat{\beta}, \theta) \right]^{-1/2} \, \left[ \mathbf{\Sigma}(\alpha, \beta, \theta) \right]^{1/2} \, \left[ \mathbf{\Sigma}(\alpha, \beta, \theta) \right]^{-1/2} \, \sqrt{n} \, (\widehat{\alpha} - \alpha, \widehat{\beta} - \beta) \overset{d}{\to} \mathcal{N}(\mathbf{0}, \mathbf{I}).
    \end{equation*}
    Since $\left[ \mathbf{\Sigma}(\alpha, \beta, \theta) \right]^{1/2} \, \left[ \mathbf{\Sigma}(\alpha, \beta, \theta) \right]^{-1/2} = \mathbf{I}$, it follows that
    \begin{equation*}
    \left[ \mathbf{\Sigma}(\widehat{\alpha}, \widehat{\beta}, \theta) \right]^{-1/2} \, \sqrt{n} \, (\widehat{\alpha} - \alpha, \widehat{\beta} - \beta) \overset{d}{\to} \mathcal{N}(\mathbf{0}, \mathbf{I}).
    \end{equation*}

\item[(iv)] Since $\widehat{\theta}$ given by \eqref{thetaest} is the $P$-consistent estimator of the unknown autocorrelation parameter $\theta$, the proof follows by the same argument as the proof of (iii).
\end{itemize}
\end{proof}

\begin{rem} \label{rem:Sqrt(n)Con}
\emph{Since $\left[ \mathbf{\Sigma}(\alpha, \beta, \theta) \right]^{-1/2} \sqrt{n} \left( \widehat{\alpha} - \alpha, \widehat{\beta} - \beta \right) \overset{d}{\to} \mathcal{N}(\mathbf{0}, \mathbf{I})$, from expressions for elements of the asymptotic covariance matrix it follows that $\left (\widehat{\alpha}, \widehat{\beta} \right)$ is the $\sqrt{n}$-consistent estimator of the unknown parameter $\left( \alpha, \beta  \right)$.}
\end{rem}

\begin{rem}
\emph{Estimators \eqref{alphaest} and \eqref{betaest} of parameters $\alpha$ and $\beta$ are the method of moments estimators based on the empirical counterparts $\overline{m}_{1}$ and $\overline{m}_{2}$ given by \eqref{emp12}. It is interesting that if instead of $\overline{m}_{2}$ the empirical counterpart of the second central moment
\begin{equation*}
\overline{M}_{2} = \displaystyle\frac{1}{n} \sum_{t=1}^{n} (X_{t} - \overline{m}_{1})^{2},
\end{equation*}
is used, the method of moments results in the same estimator for parameter $\alpha$, i.e. the result is the estimator \eqref{alphaest}. In particular, the use of $\overline{M}_{2}$ results in the estimator of parameter $\alpha$ which is given by the expression
\begin{equation}
\widetilde{\alpha} = \frac{2 \overline{m}_{1}^{\, 2}}{\overline{M}_{2} (2 - \overline{m}_{1}) - \overline{m}_{1}^{\, 2} (\overline{m}_{1} - 1)}. \label{alphaestC}
\end{equation}
Since $\overline{M}_{2} = \overline{m}_{2} - \overline{m}_{1}^{2}$, substituting this in the expression \eqref{alphaestC} results in the expression \eqref{alphaest}.}
\end{rem}

\medskip

\subsection{Martingale estimating functions} \label{MartingaleEstimation}

In this section we briefly discuss estimation of the unknown parameter $\psi = (\alpha, \beta, \theta)$ of the Fisher-Snedecor diffusion by the optimally weighted martingale estimating function approach due to Forman and S\o rensen \cite{Forman}. The result presented here is just a special case of the more general result for estimation of unknown parameters of Pearson diffusions, specialized for the Fisher-Snedecor diffusion.

In particular, the martingale estimating function for the Fisher-Snedecor diffusion is given by $$G_{n}(\psi) = \sum\limits_{i=1}^{n} \sum\limits_{j = 1}^{N} w_{j}(X_{i-1}, \psi) \left( F_{j}(X_{i}) -e^{-\lambda_{j}} F_{j}(X_{i-1}) \right), \quad N \in \{1, \ldots \left\lfloor \beta/4 \right\rfloor\},$$ where
\begin{itemize}
  \item $w_{j}(\cdot, \psi)$ are the optimal weight functions chosen to minimize the asymptotic variance of the estimator (cf. expressions (3.4) and (3.5) in \cite{Forman}),
  \item eigenvalues $\lambda_{j}$ given by expression \eqref{diseigen} depend on the unknown parameters $\beta$ and $\theta$,
  \item Fisher-Snedecor polynomials $F_{j}(\cdot)$ defined by Rodrigues formula \eqref{rodfsnorm} depend on the unknown parameters $\alpha$ and $\beta$.
\end{itemize}
According to the Theorem 3.1. from \cite{Forman}, the system of martingale estimating equations $G_{n}(\psi) = 0$ provides consistent and asymptotically normal estimator $\widetilde{\psi} = (\widetilde{\alpha}, \widetilde{\beta}, \widetilde{\theta})$ of the parameter $\psi = (\alpha, \beta, \theta)$ of the ergodic Fisher-Snedecor diffusion with the restriction $\beta > 2N$ on the value of the parameter $\beta$, ensuring the existence of the first $N$ Fisher-Snedecor polynomials, if the mapping $\psi \mapsto \varphi = (\theta, \mu, b_{2}, b_{1}, b_{0})$ is differentiable and the matrix $\partial \varphi/\partial \psi^{\tau}$ has the full rank, i.e. the rank three. Here $$\varphi = (\theta, \mu, b_{2}, b_{1}, b_{0}) = \left( \theta, \frac{\beta}{\beta - 2}, \frac{2}{\beta-2}, \frac{2 \beta}{\alpha(\beta - 2)}, 0 \right)$$ is the canonical parameter of the ergodic Fisher-Snedecor diffusion (cf. equation (2.1.)) and the corresponding matrix of partial derivatives
$$\frac{\partial \varphi}{\partial \psi^{\tau}} =
    \left[ \begin{array}{ccccc}
        0 & 0                          & 0                          & - \frac{2 \beta}{\alpha^{2}(\beta - 2)} & 0 \\
        0 & -\frac{2}{(\beta - 2)^{2}} & -\frac{2}{(\beta - 2)^{2}} & -\frac{4}{\alpha(\beta - 2)^{2}}        & 0 \\
        1 & 0                          & 0                          & 0                                       & 0
    \end{array} \right]
$$
has the rank three. Therefore, Theorem 3.1. from Forman and S\o rensen \cite{Forman} ensures the existence of consistent and asymptotically normal estimator of the parameter $\psi = (\alpha, \beta, \theta)$ of the Fisher-Snedecor diffusion. The procedure for calculating the explicit form of the asymptotic covariance matrix in the asymptotic normality framework is also given in the Theorem 3.1., but this procedure seems to be more complicated than the procedure based on the orthogonality of Fisher-Snedecor polynomials and related to method of moments estimators of the parameter $(\alpha, \beta)$ of the invariant distribution of the Fisher-Snedecor diffusion (cf. Subsection 5.2.).

In conclusion, martingale estimation function approach provides optimal estimators of the unknown parameters of the Fisher-Snedecor diffusion, in the sense that the variances of estimators are minimized, but the calculation of the asymptotic covariance matrix could be complicated. On the other hand, method of moments provides estimators in the explicit form and, since these estimators are continuous transformations of the empirical counterparts of the first and the second moment, there is an elegant procedure for calculation of the asymptotic covariance matrix in the asymptotic normality framework which makes this method more applicable and easier to use.

\medskip

\section{Testing Fisher-Snedecor distributional assumptions}

\subsection{Stein equation for Fisher-Snedecor diffusion}

The results of this section are based on the general theory developed by Barbour \cite{Barbour} and Schoutens \cite{Schoutens} and derived under the following assumptions on values of parameters $\alpha$ and $\beta$: $$\alpha > 2, \quad \alpha \notin \{2(m+1), \, m \in \mathbb{N}\}, \quad \beta > 4.$$

The general form of the Stein equation is
\begin{equation}
h(x) - E[h(Z)] = \mathcal{A}f_{h}(x),  \label{steineq}
\end{equation}
where $h$ and $f_{h}$ are smooth functions, $Z$ is the random variable with the known distribution and $\mathcal{A}$ is the Stein operator closely related to the distribution of the random variable $Z$. For every smooth function $h$ the wanted solution of the equation \eqref{steineq} is the function $f_{h}$ with the property that for each random variable $Y$ the equality
\begin{equation*}
E[h(Y)] - E[h(Z)] = E[\mathcal{A}f_{h}(Y)]
\end{equation*}
holds. According to Schoutens \cite{Schoutens}, distribution of the random variable $Y$ is close to the known distribution of the random variable $Z$ if and only if $E[\mathcal{A}f_{h}(Y)]$ is close to zero for a large class of functions $f_{h}$ and a proper choice of the Stein operator $\mathcal{A}$. According to Barbour \cite{Barbour}, the best choice for the Stein operator $\mathcal{A}$ related to the invariant distribution of the diffusion process is the corresponding infinitesimal generator.

Substitution of the infinitesimal generator \eqref{infgen} of the Fisher-Snedecor diffusion and the random variable $X_{t} \sim \mathop{\mathrm{\mathcal{FS}}}(\alpha, \beta)$ in the equation \eqref{steineq} results in the Stein-Markov equation
\begin{equation}
h(x) - E[h(X_{t})] = \frac{2\theta}{\alpha (\beta - 2)} \, x (\alpha x + \beta) g^{\prime \prime}(x) - \theta \left( x - \frac{\beta}{\beta - 2} \right) g^{\prime}(x),  \label{stmeq}
\end{equation}
where $g(\cdot)$ is the smooth function. Substitution of the function $f_{h}(x) = g^{\prime}(x)$ into the equation \eqref{stmeq} results in the equation
\begin{equation}
h(x) - E[h(X_{t})] = \frac{2\theta}{\alpha (\beta - 2)} \, x (\alpha x + \beta) f_{h}^{\prime}(x) - \theta \left( x - \frac{\beta}{\beta - 2} \right)
f_{h}(x)  \label{FSsteineq}
\end{equation}
which will be called the Stein equation for Fisher-Snedecor diffusion. Considering Barbour's method, the corresponding Stein operator is given by the expression
\begin{equation}
\mathcal{A}f_{h}(x) = \frac{2\theta}{\alpha (\beta - 2)} \, x (\alpha x + \beta) f_{h}^{\prime}(x) - \theta \left( x - \frac{\beta}{\beta - 2} \right)
f_{h}(x), \quad x \in \langle 0, \infty \rangle.
\end{equation}
According to the general result from Schoutens \cite{Schoutens} for Pearson diffusions, solution $f_{h}$ of the Stein equation \eqref{FSsteineq} is of the form
\begin{equation*}
\begin{tabular}{rcl}
$f_{h}(x)$ & = & $\displaystyle\frac{2}{\sigma^{2}(x) {\mathop{\mathrm{\mathfrak{fs}}}}(x)} \, \int\limits_{-\infty}^{x} \left(
h(y) - E[h(X_{t})] \right) {\mathop{\mathrm{\mathfrak{fs}}}}(y) \, dy \, =$
\\
& = & $\displaystyle\frac{-2}{\sigma^{2}(x) {\mathop{\mathrm{\mathfrak{fs}}}}(x)} \, \int\limits_{x}^{\infty} \left( h(y) - E[h(X_{t})]
\right) {\mathop{\mathrm{\mathfrak{fs}}}}(y) \, dy$, \\
&  &
\end{tabular}
\end{equation*}
where $\sigma^{2}(x) = \frac{4 \theta}{\alpha(\beta - 2)} \, x(\alpha x + \beta)$ is the squared diffusion coefficient of the Fisher-Snedecor diffusion.

\begin{rem}
\emph{Solution $f_{h}(x)$ of the equation \eqref{FSsteineq} for the identity function $h(x) = x$ has the surprisingly simple form:
\begin{equation*}
f_{h}(x) = \frac{1}{\alpha(2-\beta)}.
\end{equation*}}
\end{rem}

If $Y$ is the random variable with the unknown distribution and $X_{t} \sim \FS(\alpha, \beta)$, the Stein equation \eqref{FSsteineq}
takes the form
\begin{equation}
\frac{2\theta}{\alpha (\beta - 2)} \, Y (\alpha Y + \beta) f_{h}^{\prime}(Y) - \theta \left( Y - \frac{\beta}{\beta - 2} \right) f_{h}(Y) = h(Y) -
E[h(X_{t})].  \label{steineqY}
\end{equation}
Applying the expectation on the both sides of \eqref{steineqY} results in the expression
\begin{equation*}
E \left[ \frac{2\theta}{\alpha (\beta - 2)} \, Y (\alpha Y + \beta) f_{h}^{\prime}(Y) - \theta \left( Y - \frac{\beta}{\beta - 2} \right) f_{h}(Y) \right] = E[h(Y)] - E[h(X_{t})].
\end{equation*}
In particular, according to the general result from Schoutens \cite{Schoutens}, random variable $Y$ on $\langle e_{1}, e_{2} \rangle$ with the density function $f_{Y}(\cdot)$ and the finite second moment has Fisher-Snedecor distribution if and only if
\begin{equation}
E \left[ \frac{2\theta}{\alpha (\beta - 2)} \, Y (\alpha Y + \beta) f_{h}^{\prime}(Y) - \theta \left( Y - \frac{\beta}{\beta - 2} \right)
f_{h}(Y) \right] = 0  \label{momentcond}
\end{equation}
for all real continuous and piecewise continuously differentiable functions $f_{h}$ on the interval $\langle e_{1}, e_{2} \rangle$ for which the function $g(x) \equiv |\frac{1}{2} \sigma^{2}(x) f_{h}^{\prime}(x)| + |\mu(x) f_{h}(x)|$ is bounded, where $\mu(x)$ and $\sigma(x)$ are infinitesimal parameters of the Fisher-Snedecor diffusion given by \eqref{dd}. If we define function $f_{h}(\cdot)$ by the expression
\begin{equation*}
f_{h}(x) = \displaystyle\frac{\beta - 2}{\theta n (2n - \beta)} \, F^{\prime}_{n}(x), \quad n = 1, \ldots, \left\lfloor \frac{\beta}{4} \right\rfloor,
\end{equation*}
where $F^{\prime}_{n}(x)$ is the first derivative of the $n$-th Fisher-Snedecor polynomial $F_{n}(x)$, the moment condition \eqref{momentcond} takes the form
\begin{equation}
\displaystyle\frac{\beta - 2}{\theta n (2n - \beta)} \, E \left[ \frac{2\theta}{\alpha (\beta - 2)} \, Y (\alpha Y + \beta) F_{n}^{\prime \prime}(Y) - \theta \left( Y - \frac{\beta}{\beta - 2} \right) F_{n}^{\prime}(Y) \right] = 0, \quad n = 1, \ldots, \left\lfloor \frac{\beta}{4} \right\rfloor.  \label{momentcondfh}
\end{equation}
According to the Sturm-Liouville differential equation \eqref{saeq}, the final form of the moment condition \eqref{momentcondfh} is
\begin{equation}
E[F_{n}(Y)] = 0, \quad n = 1, \ldots, \left\lfloor \frac{\beta}{4} \right\rfloor.
\end{equation}
Results of the previous analysis and the following remark will be used for construction of the statistical test for Fisher-Snedecor distributional
assumptions.

\begin{rem} \label{rem:MomentCond}
\emph{If random variable $Y$ has Fisher-Snedecor distribution with parameters $\alpha > 2$ and $\beta > 4$ then
\begin{equation*}
E[F_{n}(Y)] = 0, \quad n = 1, \ldots, \left\lfloor \frac{\beta}{4} \right\rfloor.
\end{equation*}}
\end{rem}

\medskip

\subsection{Statistical test for known parameters $\protect\alpha$ and $\protect\beta$}

Null hypothesis of the form
\begin{equation*}
H_{0}: \quad {\text{The probability density function of random variable $X_{t}$ is}} \, \, \fs(x)
\end{equation*}
is considered here under assumption that parameters $\alpha$ and $\beta$ are known. Furthermore, statistical test for testing this distributional assumption about invariant distribution of Fisher-Snedecor diffusion is constructed. The basic idea for construction of such a test follows from the GMM approach for testing distributional assumptions developed by Bontemps and Maddahi \cite{Bontemps2}. However, statistical hypothesis concerning invariant Fisher-Snedecor distribution is not treated in that paper.

We consider the sample $X_{1}, \ldots, X_{n}$ from the ergodic diffusion process $\{X_{t}, \, t \geq 0\}$ with invariant Fisher-Snedecor distribution.
Autocorrelation function \eqref{acf} of Fisher-Snedecor diffusion implies that this sample has a decreasing exponential correlation structure. Let us define the function $h \colon \mathbb{R} \to \mathbb{R}^{m}$ which for $x \in \mathbb{R}$ has a value of the form
\begin{equation}
h(x) = \left[
\begin{array}{c}
2 x (\alpha x + \beta) f_{1}^{\prime}(x) + \alpha \left[\beta - (\beta - 2)x \right] f_{1}(x) \\
\vdots \\
2 x (\alpha x + \beta) f_{m}^{\prime}(x) + \alpha \left[\beta - (\beta - 2)x \right] f_{m}(x) \\
\end{array}
\right],  \label{h}
\end{equation}
where $f_{1}(\cdot), \ldots, f_{m}(\cdot)$ are real continuous and piecewise continuously differentiable functions. For example, we can define the function $h(\cdot)$ in terms of functions $f_{1}(\cdot) = F_{1}(\cdot), \ldots, f_{m}(\cdot) = F_{m}(\cdot)$, i.e. in terms of the Fisher-Snedecor polynomials or even their linear combinations, e.g.
\begin{equation}
h(x) = [F_{1}(x), \ldots, F_{m}(x)]^{\tau}, \quad m \in \left\{1, \ldots, \left\lfloor \beta/4 \right\rfloor \right\}. \label{hfs}
\end{equation}
In this case the Remark \ref{rem:MomentCond}. yields the moment condition
\begin{equation}
E[h(X_{k})] = 0, \quad  X_{k} \sim \mathop{\mathrm{\mathcal{FS}}}(\alpha, \beta), \label{exphfs}
\end{equation}
which is the basis for construction of the test for Fisher-Snedecor distributional assumptions. The simplest case of the moment condition \eqref{exphfs} is based on the function
\begin{equation*}
h(x) = F_{j}(x), \quad j \in \left\{1, \ldots, \left\lfloor \beta/4 \right\rfloor \right\},
\end{equation*}
i.e. the moment condition is of the form $E[F_{j}(X_{t})] = 0$, $X_{t} \sim \FS(\alpha, \beta)$. According to the Theorem \ref{MFCLT}. it
follows that
\begin{equation}
\frac{1}{\sqrt{n}} \sum_{t=1}^{n} F_{j}(X_{t}) \overset{d}{\to} \mathcal{N}(0, \sigma_{j}^{2}),  \label{ANhfj}
\end{equation}
where
\begin{equation*}
E[F_{j}(X_{t})] = 0, \quad \sigma_{j}^{2} = \mathop{\mathrm{Var}} \left(F_{j}(X_{0}) \right) + 2 \sum\limits_{k=1}^{\infty} \mathop{\mathrm{Cov}} \left(F_{j}(X_{0}), F_{j}(X_{k}) \right) = \coth{\left(\frac{\theta}{2(\beta - 2)} \, j (\beta - 2j)\right)}.
\end{equation*}

Asymptotic normality \eqref{ANhfj} implies that for testing Fisher-Snedecor distributional assumptions we could use the following test statistic:
\begin{equation}
\left[\coth{\left(\frac{\theta}{2(\beta - 2)} \, j (\beta - 2j)\right)}\right]^{-1} \left( \frac{1}{\sqrt{n}} \sum_{t=1}^{n} F_{j}(X_{t}) \right)^{2} \sim \chi^{2}_{1}.  \label{TShfj}
\end{equation}

Furthermore, we can define the test statistic based on the function \eqref{hfs}, i.e. in terms of $m$, $\in \left\{1, \ldots, \left\lfloor \beta/4 \right\rfloor \right\}$, Fisher-Snedecor polynomials. According to the Theorem \ref{MFCLT}. it follows that
\begin{equation}
\frac{1}{\sqrt{n}} \sum_{t=1}^{n} h(X_{t}) \overset{d}{\to} \mathcal{N} \left( \mathbf{0}, \mathbf{\Sigma(\beta, \theta)} \right),  \label{ANhfs}
\end{equation}
where the asymptotic expectation is $\mathbf{0}$ and the $(m \times m)$ covariance matrix $\mathbf{\Sigma(\beta, \theta)} = [\sigma^{2}_{ij}(\beta,
\theta)]_{m \times m}$ is the diagonal matrix whose elements are of the form
\begin{equation*}
\sigma^{2}_{ij}(\beta, \theta) = \mathop{\mathrm{Cov}}(F_{i}(X_{0}), \, F_{j}(X_{0})) + \sum_{k=1}^{\infty} \mathop{\mathrm{Cov}}(F_{i}(X_{0}), \, F_{j}(X_{k})) + \sum_{k=1}^{\infty} \mathop{\mathrm{Cov}}(F_{i}(X_{k}), \, F_{j}(X_{0})) \, =
\end{equation*}
\begin{equation*}
= \, \coth{\left(\frac{\theta}{2(\beta - 2)} \, j (\beta - 2j)\right)} \, \, \delta_{ij}, \quad i, j \in \left\{1, \ldots, \left\lfloor \beta/4 \right\rfloor \right\}.
\end{equation*}
Then the test statistic for testing Fisher-Snedecor distributional assumptions has the following asymptotic distribution:
\begin{equation}
\left( \frac{1}{\sqrt{n}} \sum_{t=1}^{n} h(X_{t}) \right)^{\tau} \mathbf{\Sigma(\beta, \theta)}^{-1} \left( \frac{1}{\sqrt{n}} \sum_{t=1}^{n}
h(X_{t}) \right) \sim \chi^{2}_{m}. \label{TShfs}
\end{equation}

\begin{rem} \label{rem:AsympNorm}
\emph{Since $\widehat{\beta}$ given by \eqref{betaest} and $\widehat{\theta}$ given by \eqref{thetaest} are the $P$-consistent estimators of parameters $\beta$ and $\theta$ and $\sigma^{2}_{ij}(\beta, \theta)$, $i, j \in \left\{1, \ldots, \left\lfloor \beta/4 \right\rfloor \right\}$, is continuous function of $\beta$ and $\theta$ it follows that
\begin{equation*}
[\mathbf{\Sigma(\widehat{\beta}, \widehat{\theta})}]^{-1/2} \, \left[ \mathbf{\Sigma(\beta, \theta)} \right]^{1/2} \overset{P}{\to} \mathbf{I}, \quad n \to \infty.
\end{equation*}
Now positive definiteness of the covariance matrix $\mathbf{\Sigma(\beta, \theta)}$, asymptotic normality \eqref{ANhfs} and Slutsky Lemma imply
\begin{equation}
[\Sigma(\widehat{\beta}, \widehat{\theta})]^{-1/2} \, \frac{1}{\sqrt{n}} \sum_{t=1}^{n} h(X_{t}) \overset{d}{\to} \mathcal{N} \left( \mathbf{0}, \mathbf{I} \right). \label{ANhfs2}
\end{equation}}
\end{rem}

\begin{rem}
\emph{Proposed goodness of fit test is not consistent against certain fixed alternatives. For instance, one can consider an alternative related to the
probability density function $\f(x)$ with support $\langle e_{1}, e_{2} \rangle$ which is different from $\fs(x)$, but such that
\begin{equation*}
\int\limits_{e_{1}}^{e_{2}} F_{n}(x) \f(x) \, dx = 0, \quad n \in \left\{1, \ldots, \left\lfloor \beta/4 \right\rfloor \right\}.
\end{equation*}
Probability density function $\f(x)$ could also be used for construction of such a test. According to the nature of our test, it is difficult to expect that it would be consistent against such an alternative. However, the goal of this paper is the statistical analysis of Fisher-Snedecor diffusion which is one of six Pearson diffusions (see the classification scheme in Introduction). This fact implies that our test is surely consistent in this class of stochastic processes. Moreover, since we are dealing with the process with known dependence structure it is hard to obtain the distribution of the test statistics that correspond to some classical tests. However, these remarks on consistency of our test and its advantages/disadvatages in comparison with some classical tests demand further investigation.}
\end{rem}

\medskip

\subsection{Testing hypothesis with parameter uncertainty}

The test statistic given in terms of unknown parameters and test statistic given in terms of consistent estimators of these parameters can have different asymptotic distributions. Therefore we observe an alternative of the statistical test constructed in the Subsection 6.2. in the framework of the parameter uncertainty. The approach to this problem is based on the moment condition \eqref{exphfs}, where instead of the unknown parameter $\omega = (\alpha, \beta)$ its $\sqrt{n}$-consistent estimator $\widehat{\omega} = (\widehat{\alpha}, \widehat{\beta})$ will be used (see Remark \ref{rem:Sqrt(n)Con}.).

The basic idea is to modify the test statistic \eqref{TShfs} by replacing $\frac{1}{\sqrt{n}} \sum_{t=1}^{n} h\left( X_{t} \right)$ with
\begin{equation}
\frac{1}{\sqrt{n}} \sum_{t=1}^{n} h\left( X_{t}, (\widehat{\alpha}, \widehat{\beta}) \right) = \frac{1}{\sqrt{n}} \sum_{t=1}^{n} h\left( X_{t}, \widehat{\omega} \right). \label{est}
\end{equation}
It means that in the parameter uncertainty framework we would like to observe the test statistic
\begin{equation}
\left( \frac{1}{\sqrt{n}} \sum_{t=1}^{n} h(X_{t}, \widehat{\omega}) \right)^{\tau} \mathbf{\Sigma(\widehat{\beta}, \widehat{\theta})}^{-1} \left( \frac{1}{\sqrt{n}} \sum_{t=1}^{n} h(X_{t}, \widehat{\omega}) \right).  \label{TShfsEst}
\end{equation}
Asymptotic distributions of the test statistics \eqref{TShfs} and \eqref{TShfsEst} coincide only if $$\frac{1}{\sqrt{n}} \sum_{t=1}^{n} h\left( X_{t} \right) \quad \textrm{and} \quad \frac{1}{\sqrt{n}} \sum_{t=1}^{n} h\left( X_{t}, \widehat{\omega} \right)$$ have the same asymptotic distributions. According to the Remark \ref{rem:AsympNorm}. it means that we need to check if
\begin{equation}
\left[ \mathbf{\Sigma(\widehat{\beta}, \widehat{\theta})} \right]^{-1/2} \frac{1}{\sqrt{n}} \sum_{t=1}^{n} h\left( X_{t}, \widehat{\omega} \right) \overset{d}{\to} \mathcal{N}
\left( \mathbf{0}, \mathbf{I} \right)  \label{ANhfsEst}
\end{equation}
holds. The classical method for checking the asymptotic normality of \eqref{est} is the Taylor approximation around the unknown parameter $\omega = (\alpha, \beta)$. In particular,
\begin{equation}
\frac{1}{\sqrt{n}} \sum_{t=1}^{n} h\left( X_{t}, \widehat{\omega} \right) = \frac{1}{\sqrt{n}} \sum_{t=1}^{n} h(X_{t}, \omega) + \left[ \frac{1}{n} \sum_{t=1}^{n} \frac{\partial h(X_{t}, \omega)}{\partial \omega^{\tau}} \right] \, \sqrt{n} (\widehat{\omega} - \omega) + o_{m}(1).  \label{TA}
\end{equation}
where $\frac{1}{n} \sum_{t=1}^{n} \frac{\partial h(X_{t}, \omega)}{\partial \omega^{\tau}}$ stands for the sample counterpart of the first moment of the random variable $\frac{\partial h(X_{t}, \omega)}{\partial \omega^{\tau}}$.

If we define
\begin{equation*}
M_{h} = \lim_{n \to \infty} \displaystyle\frac{1}{n} \sum_{t=1}^{n} \displaystyle\frac{\partial h(X_{t}, \omega)}{\partial \omega^{\tau}} = E \left[\displaystyle\frac{\partial h(X_{t}, \omega)}{\partial \omega^{\tau}} \right],
\end{equation*}
the Taylor approximation \eqref{TA} could be written as follows:
\begin{equation}
\frac{1}{\sqrt{n}} \sum_{t=1}^{n} h\left( X_{t}, \widehat{\omega} \right) = [I_{m \times m} \, \, M_{h}] \, \left[ \begin{array}{c}
\frac{1}{\sqrt{n}} \sum_{t=1}^{n} h\left( X_{t}, \omega \right) \\ \sqrt{n}(\widehat{\omega} - \omega) \\ \end{array} \right] + o_{m}(1), \label{TA2}
\end{equation}
where $m$ is the number of the variables of function $h$ and $I_{m \times m}$ is the identity matrix.

Relation \eqref{TA2} implies that the asymptotic distribution of \eqref{est} depends on the asymptotic distributions of $\frac{1}{\sqrt{n}} \sum_{t=1}^{n} h\left( X_{t}, \omega \right) = \frac{1}{\sqrt{n}} \sum_{t=1}^{n} h\left( X_{t} \right)$ and $\sqrt{n}(\widehat{\omega} - \omega) = \sqrt{n}(\widehat{\alpha} - \alpha, \widehat{\beta} - \beta)$ which are given by \eqref{ANhfs2} and the Theorem \ref{Asymptotic}.(iv), respectively. According to Bontemps and Meddahi \cite{Bontemps2}, asymptotic distributions of $\frac{1}{\sqrt{n}} \sum_{t=1}^{n} h\left( X_{t} \right)$ and $\frac{1}{\sqrt{n}} \sum_{t=1}^{n} h\left( X_{t}, \widehat{\omega} \right)$ coincide only if $M_{h} = 0$. In that case we say that the test statistic \eqref{TShfsEst} is robust against the parameter uncertainty.

To check if this is true, we use the approach based on the generalized information matrix equality (see Bontemps and Meddahi \cite{Bontemps2}). Since the moment condition \eqref{exphfs} holds for the test function $h(\cdot)$ defined in terms of Fisher-Snedecor polynomials, the generalized information matrix equality is
\begin{equation}
E \left[ \frac{\partial h(X_{t}, \omega)}{\partial \omega^{\tau}} \right] + E[h^{\tau}(X_{t}, \omega) \, s(X_{t}, \omega)] = M_{h} + E[h^{\tau}(X_{t},
\omega) \, s(X_{t}, \omega)] = 0,
\end{equation}
where $s(X_{t}, \omega)$ is the score function of the random variable $X_{t} \sim \FS(\alpha, \beta)$. It follows that the condition $M_{h} = 0$, which ensures the robustness of the test statistic \eqref{TShfsEst} against the parameter uncertainty, holds if and only if $E[h^{\tau}(X_{t}, \omega) \, s(X_{t}, \omega)] = 0$, i.e. if $h^{\tau}(X_{t}, \omega)$ and $s(X_{t}, \omega)$ are orthogonal. When $X_{t} \sim {\FS}(\alpha, \beta)$, the score function is
\begin{equation}
s(x, \omega) = \left[
\begin{array}{c}
\frac{1}{2(\alpha x + \beta)} \left[\beta(1 - x) + (\alpha x + \beta)\log{\frac{\alpha x}{\alpha x + \beta}} \right] - (\alpha x + \beta) \left[ \frac{\Gamma^{\prime}\left( \frac{\alpha}{2} \right)}{\Gamma\left( \frac{\alpha}{2} \right)} - \frac{\Gamma^{\prime}\left( \frac{\alpha + \beta}{2} \right)}{\Gamma\left( \frac{\alpha + \beta}{2} \right)} \right]
\\ \\
\frac{1}{2(\alpha x + \beta)} \left[\alpha(x - 1) + (\alpha x + \beta)\log{\frac{\beta}{\alpha x + \beta}} \right] + (\alpha x + \beta) \left[ \frac{\Gamma^{\prime}\left( \frac{\alpha + \beta}{2} \right)}{\Gamma\left( \frac{\alpha + \beta}{2}\right)} - \frac{\Gamma^{\prime}\left( \frac{\beta}{2} \right)}{\Gamma\left( \frac{\beta}{2} \right)} \right] \\
\end{array}
\right].  \label{score}
\end{equation}

From the definition of the test function \eqref{hfs} in terms of only two Fisher-Snedecor polynomials and the expression \eqref{score} for the score function we conclude that $h^{\tau}(X_{t}, \omega)$ and $s(X_{t}, \omega)$ are not orthogonal, i.e. Fisher-Snedecor polynomials are not orthogonal to the components of the score function $s(x, \omega)$. It implies that Fisher-Snedecor polynomials are not robust against parameter uncertainty for testing hypothesis about Fisher-Snedecor invariant distribution.

\bigskip \bigskip

\textbf{Acknowledgements} \newline

The research is partly supported by the EPSRC grant EP/D057361 and Marie Curie grant of European Communities PIRSES-GA-2008-230804, the Welsh Institute of Mathematics and Computer Science and the grant of Croatian Foundation for Science, Higher Education and Technological Development.

\bigskip \bigskip

\begin{center}
\textbf{Appendix A}

\textbf{Fisher-Snedecor Polynomials}
\end{center}

In this section we refer to some properties of the finite system of polynomials related to the Fisher-Snedecor distribution which are called Fisher-Snedecor polynomials. The general theory of orthogonal systems of polynomials can be found in the book by Chihara \cite{Chichara}.

\medskip

\textbf{A.1. Pearson equation}

The Pearson equation (see Pearson \cite{Pearson}) for Fisher-Snedecor distribution is the differential equation
\begin{equation}
\frac{{\fs}^{\prime}(x)}{{\fs}(x)} = \frac{\beta(\alpha - 2) - \alpha(\beta + 2)x}{2x(\alpha x + \beta)} = \frac{q(x)}{s(x)}. \tag{A.1} \label{pearson}
\end{equation}

Polynomial solutions of the hypergeometric equation
\begin{equation*}
s(x) \, y^{\prime \prime}(x) + [q(x) + s^{\prime}(x)] \, y^{\prime}(x) + \lambda \, y(x) = 0,
\end{equation*}
where $s(x)$ and $q(x)$ are polynomials from the Pearson equation \eqref{pearson}, are known as Fisher-Snedecor polynomials.

\medskip

\textbf{A.2. Rodrigues formula}

There exists the unique integer $N \in \mathbb{N}_{0}$ such that
\begin{equation*}
\int\limits_{0}^{\infty} s^{N}(x) \fs(x) \, dx < \infty, \quad n \in \{0, \ldots, N\} \text{\quad \quad and \quad} \int\limits_{-\infty}^{\infty} s^{n}(x) \fs(x) \, dx = \infty, \quad n \in \{N+k, \, k \in \mathbb{N}\}.
\end{equation*}
In particular, for the Fisher-Snedecor distribution $N = \left\lfloor \frac{\beta}{4} \right\rfloor$, $\beta > 2$, and for each $n \in \left\{0, \ldots, \left\lfloor \beta/4 \right\rfloor \right\}$ the
function
\begin{equation}
\widetilde{F}_{n}(x) = x^{1-\frac{\alpha}{2}} \, (\alpha x + \beta)^{\frac{\alpha}{2} + \frac{\beta}{2}} \, \frac{d^{n}}{dx^{n}} \, \left\{ 2^{n} \, x^{\frac{\alpha}{2} + n - 1} \, (\alpha x + \beta)^{n - \frac{\alpha}{2} - \frac{\beta}{2}} \right\}, \tag{A.2} \label{rodfs}
\end{equation}
is the polynomial of at most $n$-th degree. This formula is known as the Rodrigues formula while the polynomials defined by \eqref{rodfs} are non-normalized (with respect to the Fisher-Snedecor density \eqref{fs}) Fisher-Snedecor polynomials. The first six polynomials of this type are: \\

$\widetilde{F}_{0}(x) = 1$, \vspace*{1mm}

$\widetilde{F}_{1}(x) = -\alpha(\beta - 2)x + \alpha \beta$, \vspace*{1mm}

$\widetilde{F}_{2}(x) = \alpha^{2}(\beta - 4)(\beta - 6)x^{2} - 2\alpha \beta(\alpha + 2)(\beta - 4)x + \alpha \beta^{2}(\alpha + 2)$, \vspace*{1mm}

$\widetilde{F}_{3}(x) = -\alpha^{3}(\beta - 6)(\beta - 8)(\beta - 10)x^{3} + 3\alpha^{2} \beta(\alpha + 4)(\beta - 6)(\beta - 8)x^{2} - 3\alpha
    \beta^{2}(\alpha + 2)(\alpha + 4)(\beta - 6)x \, +$ \\
    $\hspace*{1.5cm} + \, \alpha \beta^{3}(\alpha + 2)(\alpha + 4)$, \vspace*{1mm}

$\widetilde{F}_{4}(x) = \alpha^{4}(\beta - 8)(\beta - 10)(\beta - 12)(\beta - 14)x^{4} - 4\alpha^{3} \beta(\alpha + 6)(\beta - 8)(\beta - 10)(\beta -
    12)x^{3} \, +$ \\
    $\hspace*{1.5cm}+ \, 6\alpha^{2} \beta^{2}(\alpha + 4)(\alpha + 6)(\beta - 8)(\beta - 10)x^{2} - 4\alpha \beta^{3}(\alpha + 2)(\alpha + 4)(\alpha + 6)(\beta - 8)x \, +$ \\
    $\hspace*{1.5cm}+ \, \alpha \beta^{4}(\alpha + 2)(\alpha + 4)(\alpha + 6)$, \vspace*{1mm}

$\widetilde{F}_{5}(x) = -\alpha^{5}(\beta - 10)(\beta - 12)(\beta - 14)(\beta - 16)(\beta - 18)x^{5} + 5\alpha^{4} \beta(\alpha + 8)(\beta - 10)(\beta
    - 12)(\beta - 14)(\beta - 16)x^{4} \, -$ \\
    $\hspace*{1.5cm}- \, 10\alpha^{2} \beta^{2}(\alpha + 6)(\alpha + 8)(\beta - 10)(\beta - 12)(\beta - 14)x^{3} + 10\alpha^{2} \beta^{3}(\alpha + 4)(\alpha + 6)(\alpha + 8)(\beta - 10)(\beta - 12)x^{2} \, -$ \\
    $\hspace*{1.5cm}- \, 5\alpha \beta^{4}(\alpha + 2)(\alpha + 4)(\alpha + 6)(\alpha + 8)(\beta - 10)x + \alpha \beta^{5}(\alpha + 2)(\alpha + 4)(\alpha + 6)(\alpha + 8)$. \\

Fisher-Snedecor polynomials could also be defined in terms of the Gauss hypergeometric function, i.e.
\begin{equation}
F_{n}(x) = \alpha \beta^{n} 2^{n-1} \phantom{,}_{2}F_{1} \left(-n, n - \frac{\beta}{2}; \frac{\alpha}{2}; -\frac{\alpha}{\beta} \, x \right), \quad n \in \left\{0, \ldots, \left\lfloor \beta/4 \right\rfloor \right\}. \tag{A.3} \label{fphyp}
\end{equation}

\medskip

\textbf{A.3. Orthogonality property and normalizing constant}

Fisher-Snedecor polynomials defined by the Rodrigues formula \eqref{rodfs} are orthogonal in the general sense for $n \in \left\{0, \ldots, \left\lfloor \beta/4 \right\rfloor \right\}$, $\beta > 4$. Precisely, there is no set of parameters $\alpha$ and $\beta$ which ensure the orthogonality in the restricted sense (see Beale \cite{Beale}), i.e.
\begin{equation*}
\int\limits_{0}^{\infty} \widetilde{F}_{m}(x) \widetilde{F}_{n}(x) \fs(x) \, dx = 0, \quad m, n \in \left\{0, \ldots, \left\lfloor \beta/4 \right\rfloor \right\}, \quad m \neq n.
\end{equation*}
To normalize them, we must multiply each of them by the normalizing constant
\begin{equation*}
\alpha_{n} = \frac{(-1)^{n}}{\sqrt{\left( -1 \right)^{n} \, n! \, d_{n} \, I_{n}}},
\end{equation*}
where
\begin{equation*}
d_{n} = \prod_{k=0}^{n-1} \left( q^{\prime}(x) + \frac{n+k+1}{2} \, s^{\prime \prime}(x)\right), \quad I_{n} = \int\limits_{-\infty}^{\infty} s^{n}(x) \fs(x) \, dx.
\end{equation*}
In this particular case $d_{n}$ and $I_{n}$ are given by the expressions
\begin{equation*}
d_{n} = \prod_{k=0}^{n-1} \left[ \alpha (2n+2k-\beta) \right] =  (2 \alpha)^{n} \, \frac{\Gamma\left( 2n-\frac{\beta}{2} \right)}{\Gamma\left( n-\frac{\beta }{2} \right)}, \quad n \in \left\{0, \ldots, \left\lfloor \beta/4 \right\rfloor \right\},
\end{equation*}
\begin{equation*}
I_{n} = \left( \frac{2\beta^{2}}{\alpha} \right)^{n} \, \frac{B\left(n+\frac{\alpha}{2}, -2n+\frac{\beta}{2} \right)}{B(\frac{\alpha}{2}, \frac{\beta }{2})}, \quad n \in \left\{0, \ldots, \left\lfloor \beta/4 \right\rfloor \right\},
\end{equation*}
while the normalizing constant $\alpha_{n}$ for Fisher-Snedecor polynomials is given by expression
\begin{equation*}
\alpha_{n} = (-1)^{n} \sqrt{\frac{\Gamma\left( n-\frac{\beta}{2} \right) B(\frac{\alpha }{2}, \frac{\beta}{2})}{n! (-1)^{n} (2 \beta)^{2n} \Gamma\left( 2n-\frac{\beta}{2} \right) B\left( n+\frac{\alpha}{2}, -2n+\frac{\beta}{2} \right)}} =
\end{equation*}
\begin{equation}
\hspace*{1cm} = (-1)^{n} \sqrt{\frac{B(\frac{\alpha}{2}, \frac{\beta}{2})}{n! (2 \beta)^{2n} B\left( \frac{\alpha}{2} + n, \frac{\beta }{2} - 2n \right)} \, \left[ \prod_{k=1}^{n} \left( \frac{\beta}{2} + k - 2n \right) \right]^{-1}}.  \tag{A.4} \label{normfs}
\end{equation}

According to \eqref{fphyp} and \eqref{normfs}, the Rodrigues formula for the normalized (with respect to the Fisher-Snedecor density \eqref{fs}) Fisher-Snedecor polynomials is given by the expression
\begin{equation}
F_{n}(x) = \alpha_{n} \, \widetilde{F}_{n}(x) = \alpha_{n} \, x^{1-\frac{\alpha}{2}} \, (\alpha x+\beta)^{\frac{\alpha}{2} + \frac{\beta}{2}} \, \frac{d^{n}}{dx^{n}} \, \left\{2^{n} \, x^{\frac{\alpha}{2}+n-1} \, (\alpha x+\beta)^{n-\frac{\alpha}{2}-\frac{\beta}{2}} \right\}. \tag{A.5} \label{rodfsnorm}
\end{equation}
In conclusion, polynomials $F_{n}(x)$, $n \in \left\{0, \ldots, \left\lfloor \beta/4 \right\rfloor \right\}$, $\beta > 4$, are normalized and orthogonal, i.e.
\begin{equation}
\int\limits_{0}^{\infty} F_{m}(x) F_{n}(x) \fs(x) \, dx = \delta_{mn}, \tag{A.6} \label{orthofs}
\end{equation}
where $\delta_{mn}$ is the standard Kronecker symbol. Relation \eqref{orthofs} implies interesting properties of the random variables $F_{n}(Y)$, where $Y \sim \FS(\alpha, \beta)$. In particular, random variables $F_{n}(Y)$ are orthonormal, i.e.
\begin{equation*}
E[F_{n}(Y) F_{m}(Y)] = \delta_{nm}, \quad n, m \in \left\{0, \ldots, \left\lfloor \beta/4 \right\rfloor \right\}.
\end{equation*}
Since $F_{0}(x) = 1$, substituting $n \neq 0$ and $m = 0$ to the previous expression we see that random variables $F_{n}(Y)$ are centered, i.e.
\begin{equation*}
E[F_{n}(Y)] = 0,\quad n \in \left\{1, \ldots, \left\lfloor \beta/4 \right\rfloor \right\}.
\end{equation*}

\medskip

\textbf{A.4. Recurrence relation}

Fisher-Snedecor polynomials satisfy the following Favard-Jacobi recurrent relation:
\begin{equation*}
F_{n+1}(x) = -\frac{1}{a_{n}} \, \left((b_{n}-x) F_{n}(x) + a_{n-1} F_{n-1}(x) \right), \quad n \in \left\{0, \ldots, \left\lfloor \beta/4 \right\rfloor \right\},
\end{equation*}
where
\begin{equation*}
a_{n} = \frac{2 (\beta - 2n)}{\alpha (4n - \beta) (2 + 4n - \beta)} \, \sqrt{-\frac{(n + 1) (2n + \alpha) \beta^{2} (\alpha + \beta - 2n - 2)}{(2 + 4n - \beta) (4 + 4n - \beta)}},
\end{equation*}
\begin{equation*}
b_{n} = \frac{\beta n (2n + \alpha - 2)}{\alpha (4n - \beta - 2)} - \frac{\beta (n + 1) (2n + \alpha)}{\alpha (4n - \beta + 2)}, \quad \quad \mu_{n} =
\frac{\beta (2n + \alpha - 2)}{\alpha (4n - \beta - 2)}.
\end{equation*}

\medskip

\textbf{A.5. Sturm-Liouville equation}

Both non-normalized and normalized, with respect to the density \eqref{fs}, Fisher-Snedecor polynomials satisfy the second order self-adjoint equation
\begin{equation*}
\frac{2(\alpha x + \beta)^{\frac{\alpha + \beta}{2}}}{x^{\frac{\alpha}{2} - 1}} \, \frac{d}{dx}\left\{ \frac{x^{\frac{\alpha}{2}}}{(\alpha x + \beta)^{\frac{\alpha + \beta}{2} - 1}} y^{\prime}(x) \right\} + \lambda_{n} y(x) = 0,
\end{equation*}
which could be written as
\begin{equation}
2x (\alpha x + \beta) y^{\prime \prime}(x) + \alpha \left(\beta - (\beta - 2)x \right) y^{\prime}(x) + \lambda_{n} y(x) = 0, \tag{A.7} \label{saeq}
\end{equation}
where the spectral parameter $\lambda_{n}$ is given by the expression
\begin{equation*}
\lambda_{n} = -n \left(s^{\prime}(x) + q(x) \right)^{\prime} - \frac{n(n-1)}{2} s^{\prime\prime}(x) = \alpha n \left(\beta - 2n \right), \quad n \in \left\{0, \ldots, \left\lfloor \beta/4 \right\rfloor \right\}.
\end{equation*}
Fisher-Snedecor polynomials $F_{0}(x), \ldots, F_{\left\lfloor \frac{\beta}{4} \right\rfloor}(x)$ are the polynomial eigenfunctions of the Sturm-Liouville equation \eqref{saeq}, while real numbers $\lambda_{n}$ are the corresponding eigenvalues.

\medskip

\begin{center}
\textbf{Appendix B}

\textbf{Gauss Hypergeometric Functions}
\end{center}

Gauss hypergeometric series is the power series
\begin{equation}
_{2}F_{1}(a, b; c; z) = \sum_{k=0}^{\infty} \frac{(a)_{k} \, (b)_{k}}{(c)_{k}} \frac{z^{k}}{k!} = 1 + \frac{ab}{c} \, z + \frac{a(a+1) \, b(b+1)}{c(c+1)} \, \frac{z^{2}}{2} + \ldots, \label{GaussSeries}
\end{equation}
where $z$ is a complex variable, $a$, $b$ and $c$ are real or complex parameters and $(a)_{k}$ is the Pochhammer symbol which denotes the quantity
\begin{equation*}
(a)_{0} = 1, \quad (a)_{k} = \frac{\Gamma(a + k)}{\Gamma(a)} = a(a + 1) \cdot \ldots \cdot (a + k - 1), \quad k \in \mathbb{N}.
\end{equation*}
The series is not defined for $c = -m$, $m \in \mathbb{N}_{0}$, provided that $a$ or $b$ is not the negative integer $n$ such that $n < m$. Furthermore, if the series \eqref{GaussSeries} is defined but $a$ or $b$ is equal to $(-n)$, $n \in \mathbb{N}_{0}$, then it terminates in a finite number of terms and its sum is then the polynomial of degree $n$ in variable $z$. Except for this case, in which the series is absolutely convergent for $|z| <\infty$, the radius of absolute convergence of the series \eqref{GaussSeries} is the unit circle, i.e. $|z| < 1$. In this case it is said that the series \eqref{GaussSeries} defines the Gauss or hypergeometric function
\begin{equation}
g(z) = \,_{2}F_{1}(a, b; c; z). \label{GaussFunction}
\end{equation}
It can be verified (see Slater \cite{Slater}, Section 1.2., page 6) that the function $g(z)$ is the solution of the second order differential equation
\begin{equation}
z(1 - z) g^{\prime \prime}(z) + \left(c - (a + b + 1)z \right) g^{\prime}(z) - ab g(z) = 0, \label{GaussEquation}
\end{equation}
in the region $|z| < 1$. However, the function \eqref{GaussFunction} can be analytically continued to the other parts of the complex plane, i.e. solutions of the equation \eqref{GaussEquation} are also defined outside the unit circle. These solutions are provided by following substitutions in the equation \eqref{GaussEquation}:
\begin{itemize}
  \item substitution $z = 1 - y$ yields solutions valid in the region $|1 - z| < 1$,
  \item substitution $z = 1/y$ yields solutions valid in the region $|z| > 1$.
\end{itemize}
Furthermore, solutions in these three regions of the complex plane yield solutions valid in the regions $|1 - z| > 1$, $\textrm{Re}(z) > 1/2$ and $\textrm{Re}(z) < 1/2$. The number of solutions in these six regions is $24$ and they are known as Kummer's system of solutions. This system of solutions is very important since it provides the analytic continuation of the hypergeometric function to the whole complex plane with an exception of one branch cut which is always the subinterval of the real axis. For example, the analytic continuation of the function \eqref{GaussFunction} from the region $|z| < 1$ to the whole complex plane with an exception of the branch cut situated along the interval $[1, \infty \rangle$ is provided by the Mellin-Barnes integral (see Slater \cite{Slater}, Section 1.6.1., Expression 1.6.1.6.):
\begin{equation}
\frac{\Gamma(c)}{2 \pi i \Gamma(a) \Gamma(b)} \int\limits_{-i \infty}^{i \infty} \frac{\Gamma(a + s) \Gamma(b + s) \Gamma(-s)}{\Gamma(c + s)} \, (-z)^{s} \, ds. \label{AnCont1}
\end{equation}
Namely, Mellin-Barnes integral defines the analytic function on the complex plane with the branch cut along the positive real axis, i.e. along the interval $\langle 0, \infty \rangle$. Analytic function \eqref{GaussFunction} defined by the Gauss hypergeometric series and the analytic function defined by the Mellin-Barnes integral coincide on the region $\{z \in \mathbb{C}: |z| < 1, z \notin \langle 0, 1 \rangle \}$. Therefore, the Mellin-Barnes integral provides the analytic continuation of the function \eqref{GaussFunction} to the region $|z| > 1$ cutted along the interval $[1, \infty \rangle$. Hence, the hypergeometric function $\,_{2}F_{1}(a, b; c; z)$ is an analytic function defined on the complex plane with the branch cut situated along the interval $[1, \infty \rangle$. Similarly, the hypergeometric function $\,_{2}F_{1}(a, b; c; -z)$ is an analytic function defined on the complex plane with the branch cut situated along the interval $\langle -\infty, -1]$.

From the Mellin-Barnes integral follows one of the most important relations for analytic continuation of the hypergeometric functions (cf. Whittaker and Watson \cite{Whittaker} or Luke \cite{Luke}):
\begin{equation}
g(z) = \, _{2}F_{1}(a, b; c; z) = \frac{\Gamma(c) \Gamma(b - a)}{\Gamma(b) \Gamma(c - a)} \, (-z)^{-a} \, _{2}F_{1} \left(a, 1 - c + a; 1 - b + a; \frac{1}{z} \right) + \label{AnCont2}
\end{equation}
\begin{equation*}
\hspace*{3.5cm} + \frac{\Gamma(c) \Gamma(a - b)}{\Gamma(a) \Gamma(c - b)} \, (-z)^{-b} \, _{2}F_{1} \left(b, 1 - c + b; 1 - a + b; \frac{1}{z} \right),
\end{equation*}
where functions $(-z)^{-a} \, _{2}F_{1} \left(a, 1 - c + a; 1 - b + a; \frac{1}{z} \right)$ and $(-z)^{-b} \, _{2}F_{1} \left(b, 1 - c + b; 1 - a + b; \frac{1}{z} \right)$ are solutions of the equation \eqref{GaussEquation} in the region $|z| > 1$. A list of useful relations between solutions of equation \eqref{GaussEquation} which providing analytic continuation of the hypergeometric function, can be found in Luke \cite{Luke} (see Section 3.9., page 69).

\medskip

\begin{center}
\textbf{Appendix C}

\textbf{Feller's and Weyl's classification of boundaries of the state space of Fisher-Snedecor diffusion}
\end{center}

In this section we present classical Feller's, Weyl's limit-point/limit-circle (LP/LC) and oscillatory/non-oscillatory (O/NO) classification of boundaries of the diffusion state space. For general information on Feller's classification scheme we refer to Karlin and Taylor \cite{Karlin2} and Linetsky \cite{Linetsky1}, while for Weyl's LP/LC and O/NO classification we refer to Fulton et al. \cite{Fulton}) and Linetsky \cite{Linetsky1}. The nature of each boundary is determined according to the behavior of the Sturm-Liouville equation \eqref{sleq} near it.

\begin{thm} \label{boundary_class}
Boundaries of the state space of Fisher-Snedecor diffusion with parameters $\alpha > 0$, $\alpha \notin \{2m, \, m \in \mathbb{N}\}$, and $\beta > 2$ are classified as follows:
\begin{itemize}
  \item[(i)] Boundary $e_{1} = 0$ is regular for $\alpha < 2$ and entrance otherwise, while $e_{2} = \infty$ is natural for $\alpha > 0$, $\alpha \notin \{2m, \, m \in \mathbb{N}\}$.
  \item[(ii)] For $\alpha > 0$, $\alpha \notin \{2m, \, m \in \mathbb{N}\}$, boundary $e_{1} = 0$ is non-oscillatory, while $e_{2} = \infty$ is oscilattory/non-oscillatory with unique positive cutoff $$\Lambda = \frac{\theta \beta^{2}}{8(\beta - 2)}.$$ Boundary $e_{2} = \infty$ is non-oscillatory for $\lambda \leq \Lambda$ and oscillatory for $\lambda > \Lambda$.
  \item[(iii)] Boundary $e_{1} = 0$ is of limit-circle type for $\alpha < 4$ and of limit-point type otherwise, while boundary $e_{2} = \infty$ is of limit-point type for $\alpha > 0$, $\alpha \notin \{2m, \, m \in \mathbb{N}\}$.
\end{itemize}
\end{thm}

\begin{proof}
The proof consists of three parts, as we observe three boundary classification schemes.

\textit{(i)}
Feller's boundary classification is based on the behavior of the scale function $$S[x, y] = \int\limits_{x}^{y} \mathfrak{s}(z) \, dz$$ near both boundaries, where $\mathfrak{s}(x)$ is the scale density defined by \eqref{scale}. In case of the Fisher-Snedecor diffusion, this non-standard integral is evaluated using Mathematica:
\begin{equation*}
S[x, y] = \frac{2 \beta^{\frac{\alpha}{2} + \frac{\beta}{2} - 1}}{2 - \alpha} \left[y^{1 - \frac{\alpha}{2}} \, _{2}F_{1} \left( 1 - \frac{\alpha}{2}, 1 - \frac{\alpha}{2} - \frac{\beta}{2}; 2 - \frac{\alpha}{2}; -\frac{\alpha}{\beta} \, y \right) \, - \right.
\end{equation*}
\begin{equation*}
\left. \hspace{2.8cm} - \, x^{1 - \frac{\alpha}{2}} \, _{2}F_{1} \left( 1 - \frac{\alpha}{2}, 1 - \frac{\alpha}{2} - \frac{\beta}{2}; 2 - \frac{\alpha}{2}; -\frac{\alpha}{\beta} \, x \right) \right].
\end{equation*}
Since $\lim\limits_{x \to 0} \,_{2}F_{1} \left( 1 - \frac{\alpha}{2}, 1 - \frac{\alpha}{2} - \frac{\beta}{2}; 2 - \frac{\alpha}{2}; -\frac{\alpha}{\beta} \, x \right) = 0$ and $\,_{2}F_{1} \left( 1 - \frac{\alpha}{2}, 1 - \frac{\alpha}{2} - \frac{\beta}{2}; 2 - \frac{\alpha}{2}; -\frac{\alpha}{\beta} \, y \right)$ does not converge as $y \to \infty$, for $\alpha \notin \{2m, \, m \in \mathbb{N}\}$ the function $S[x, y]$ has the following properties:
$$S \langle 0, y] = \lim\limits_{x \to 0} S[x, y] < \infty, \quad \alpha < 2,$$
$$S \langle 0, y] = \lim\limits_{x \to 0} S[x, y] = \infty, \quad \alpha > 2,$$
$$S [x, \infty \rangle = \lim\limits_{y \to \infty} S[x, y] = \infty, \quad \alpha > 0.$$
These results imply that for $\alpha \notin \{2m, \, m \in \mathbb{N}\}$ and arbitrary $\varepsilon > 0$ we have:
$$\int\limits_{0}^{\varepsilon} S\langle 0, y] \mathfrak{m}(y) \, dy  < \infty, \, \alpha < 2, \quad \int\limits_{0}^{\varepsilon} S\langle 0, y] \mathfrak{m}(y) \, dy  = \infty, \, \alpha > 2, \quad \int\limits_{\varepsilon}^{\infty} S[x, \infty \rangle \mathfrak{m}(x) \, dx = \infty, \, \alpha > 0,$$
$$\int\limits_{\varepsilon}^{\infty} S[\varepsilon, y] \mathfrak{m}(y) \, dy = \infty, \, \alpha > 0, \quad \int\limits_{0}^{\varepsilon} S[x, \varepsilon] \mathfrak{m}(x) \, dx < \infty, \, \alpha > 0.$$ Therefore, according to the standard Feller's boundary classification scheme, $e_{1} = 0$ is regular boundary for $\alpha < 2$ and entrance boundary for $\alpha > 2$, while $e_{2} = \infty$ is natural boundary for all positive values of $\alpha$, $\alpha \notin \{2m, \, m \in \mathbb{N}\}$. \vspace*{2mm}

\textit{(ii)} Since $e_{1} = 0$ is either regular or entrance boundary, it is necessary NO (cf. Linetsky \cite{Linetsky1}). For O/NO classification of the natural boundary $e_{2} = \infty$ the standard procedure is transformation of the Sturm-Liouville equation \eqref{sleq} to the Liouville normal form (cf. Fulton et al. \cite{Fulton}). If we observe the Sturm-Liouville equation \eqref{sleq} written in the form
$$-\left(x^{\frac{\alpha}{2}}(\alpha x + \beta)^{1 - \frac{\alpha}{2} - \frac{\beta}{2}} f^{\prime}(x) \right)^{\prime} = \lambda \, \frac{\alpha(\beta - 2)}{2 \theta} \, x^{\frac{\alpha}{2} - 1}(\alpha x + \beta)^{-\frac{\alpha}{2} - \frac{\beta}{2}} \, f(x),$$ then
$$u(x) = \sqrt{\frac{2(\beta-2)}{\theta}} \, \arcsinh{\left( \sqrt{\frac{\alpha}{\beta} \, x} \right)}, \quad x(u) = \frac{\beta}{\alpha} \, {\sinh}^{2}{\left( u \sqrt{\frac{\theta}{2(\beta - 2)}} \right)},$$
$$h(u) = \sqrt[4]{\frac{\alpha (\beta - 2)}{2 \theta} \, \left[ \beta {\cosh}^{2}\left( u \sqrt{\frac{\theta}{2(\beta-2)}} \right) \right]^{1-\alpha-\beta} \, \left[ \frac{\beta}{\alpha} {\sinh}^{2}\left( u \sqrt{\frac{\theta}{2(\beta-2)}} \right) \right]^{\alpha-1}}$$
and
\begin{equation*}
Q(u) = \frac{\theta}{8(\beta-2)} \left( \beta^{2} + (\alpha-1)(\alpha-3)\csch{^{2}}{\left(u\sqrt{\frac{\theta}{2(\beta-2)}}\right)} - \left[ (\alpha + \beta)^{2} - 1 \right]\sech{^{2}}{\left(u\sqrt{\frac{\theta}{2(\beta-2)}}\right)} \right).
\end{equation*}
Therefore, the Liouville normal form of the Sturm-Liouville equation \eqref{sleq} is
\begin{equation}
-g^{\prime \prime}(u) + \left( Q(u) - \lambda \right) g(u) = 0. \label{lnf}
\end{equation}
In this particular case the natural boundary $e_{2} = \infty$ remains unchanged under the Liouville transformation, i.e. the corresponding boundary of the equation \eqref{lnf} is $e_{2}^{*} = u(e_{2}) = \infty$. O/NO classification of the boundary $e_{2}^{*} = \infty$ depends on the behavior of the potential function $Q(u)$ near that endpoint. Since $$\lim\limits_{u \to \infty} \csch{^{2}}{\left(u\sqrt{\frac{\theta}{2(\beta-2)}}\right)} = 0 \quad \textrm{and} \quad \lim\limits_{u \to \infty} \sech{^{2}}{\left(u\sqrt{\frac{\theta}{2(\beta-2)}}\right)} = 0,$$ it follows that $$\lim\limits_{u \to \infty} Q(u) = \frac{\theta \beta^{2}}{8(\beta - 2)}, \quad \beta > 2.$$ Since $\lim\limits_{u \to \infty} Q(u) < \infty$, according to Fulton et al. (cf. \cite{Fulton}, Theorem 6) $e_{2}^{*} = \infty$ is O/NO boundary of the equation \eqref{lnf} with the unique positive cutoff $$\Lambda = \frac{\theta \beta^{2}}{8(\beta - 2)}, \quad \beta > 2,$$ i.e. it is NO for $\lambda < \Lambda$ and O for $\lambda > \Lambda$. According to Linetsky (cf. \cite{Linetsky1}, Theorem 3), the boundary $e_{2}^{*} = \infty$ is non-oscillatory for $\lambda = \Lambda$. This classification of boundaries remains invariant under the Liouville transformation (cf. \cite{Fulton}, Lemma 2), i.e. $e_{2} = \infty$ is O/NO boundary of the Sturm-Liouville equation \eqref{sleq} with the unique positive cutoff $\Lambda$. Furthermore, it is NO for $\lambda \leq \Lambda$ and O for $\lambda > \Lambda$. \vspace*{2mm}

\textit{(iii)} Since $e_{1} = 0$ is the regular boundary for $\alpha < 2$, it is necessary of the LC type. For LP/LC classification of the boundary $e_{1} = 0$ for $\alpha > 2$, $\alpha \notin \{2m, \, m \in \mathbb{N}\}$, we observe its Liouville trasformation $e_{1}^{*} = u(e_{1}) = 0$ and $\liminf\limits_{u \to e_{1}^{*}} u^{2} Q(u)$ and $\limsup\limits_{u \to e_{1}^{*}} u^{2} Q(u)$ (cf. Fulton \cite{Fulton}, Theorem 5). Since $$\lim\limits_{u \to 0} u^{2} \csch{^{2}}{\left(u\sqrt{\frac{\theta}{2(\beta-2)}}\right)} = \frac{2(\beta - 2)}{\theta} \quad \textrm{and} \quad \lim\limits_{u \to 0} u^{2} \sech{^{2}}{\left(u\sqrt{\frac{\theta}{2(\beta-2)}}\right)} = 0,$$ it follows that $$\lim\limits_{u \to 0} u^{2} Q(u) = \frac{(\alpha- 1)(\alpha - 3)}{4}$$ and $$\liminf\limits_{u \to 0} u^{2} Q(u) = \limsup\limits_{u \to 0} u^{2} Q(u) = \frac{(\alpha- 1)(\alpha - 3)}{4}.$$ Now it follows that $$\liminf\limits_{u \to 0} u^{2} Q(u) < \frac{3}{4}$$ for $\alpha \in \langle 2, 4 \rangle$, while $$\limsup\limits_{u \to 0} u^{2} Q(u) \geq \frac{3}{4}$$ for $\alpha > 4$, $\alpha \notin \{2m, \, m \in \mathbb{N}\}$. According to Fulton et al., $e_{1}* = 0$ is the boundary of LC type for $\alpha \in \langle 2, 4 \rangle$ and of LP type for $\alpha > 4$, $\alpha \notin \{2m, \, m \in \mathbb{N}\}$.

For LP/LC classification of the boundary $e_{2} = \infty$ we observe LP/LC classification of its Liouville transformation $e_{2}^{*} = \infty$ which depends on the behavior of the function $Q(u)/u^{2}$ near that endpoint. Since $$\lim\limits_{u \to \infty} u^{-2} \csch{^{2}}{\left(u\sqrt{\frac{\theta}{2(\beta-2)}}\right)} = 0 \quad \textrm{and} \quad \lim\limits_{u \to \infty} u^{-2} \sech{^{2}}{\left(u\sqrt{\frac{\theta}{2(\beta-2)}}\right)} = 0,$$ it follows that $$\lim\limits_{u \to \infty} \frac{Q(u)}{u^{2}} = 0.$$ According to Fulton et al. (cf. \cite{Fulton}, Theorem 7) $e_{2}^{*}=\infty$ is the boundary of LP type for all $\alpha > 0$, $\alpha \notin \{2m, \, m \in \mathbb{N}\}$. This classification of boundaries remains invariant under the Liouville transformation (cf. \cite{Fulton}, Lemma 1), i.e. the boundary $e_{1} = 0$ is of LC type for $\alpha < 4$ and of LP type for $\alpha > 4$, while the boundary $e_{2} = \infty$ is of LP type for all $\alpha > 0$, $\alpha \notin \{2m, \, m \in \mathbb{N}\}$.
\end{proof}

\begin{rem}
\emph{Since in the Liouville normal form \eqref{lnf} of the Sturm-Liouville equation \eqref{sleq} the first derivative term vanishes, according to Buchholz \cite{Buchholz} Wronskian of its linearly independent solutions $$g_{4}(u) = h(u) f_{4}(x(u), -\lambda) = h(u) f_{4}(u) \quad \textrm{and} \quad g_{1}(u) = h(u) f_{1}(x(u), -\lambda) = h(u) f_{1}(u)$$ is constant. According to the relation \eqref{f1=Bf3+Af4}, it follows that $$g_{1}(u) = h(u) \left( B_\l f_{3}(x(u), -\lambda) + A_\l f_{4}(x(u), -\lambda) \right) = B_\l g_{3}(u) + A_\l g_{4}(u).$$ Therefore, the Wronskian of the solutions $g_{1}(u)$ and $g_{4}(u)$ is given by
\begin{equation*}
W(g_{4}(u), g_{1}(u)) = B_\l W(g_{4}(u), g_{3}(u)) = B_\l h^{2}(u) W(f_{4}(u), f_{3}(u)),
\end{equation*}
where $B_\l$ is given by the expression \eqref{B}. Expression for $h(u)$ given in the proof of the Theorem \ref{boundary_class}. implies that
\begin{equation}
W(g_{4}(u(x)), g_{1}(u(x))) \sim 2 B_\l \alpha^{1 - \frac{\alpha}{2}} \beta^{-\frac{\beta}{2}} \, \sqrt{\frac{\beta^{2}}{16} - \frac{\lambda(\beta - 2)}{2 \theta}} = 2 \alpha^{1 - \frac{\alpha}{2}} \beta^{-\frac{\beta}{2}} B_{\l} \D_{\l}, \quad x \to \infty, \tag{C.2} \label{Wg1g4}
\end{equation}
where $\D_{\l}$ is given by \eqref{Delta}. Since $W(g_{4}(u), g_{1}(u))$ must be constant, this is its value. On the other hand, we know that $W(g_{4}(u), g_{1}(u)) = h^{2}(u) x^{\prime}(u) W(f_{4}(x(u), -\lambda), f_{1}(x(u), -\lambda))$.
It implies that
\begin{equation*}
W(f_{4}(x(u), -\lambda), f_{1}(x(u), -\lambda)) = \frac{W(g_{4}(u), g_{1}(u))}{h^{2}(u) x^{\prime}(u)}
\end{equation*}
i.e.
\begin{equation}
W(f_{4}(x, -\lambda), f_{1}(x, -\lambda)) = \frac{W(g_{4}(u(x)), g_{1}(u(x)))}{h^{2}(u(x)) x^{\prime}(u(x))}, \tag{C.3} \label{Wf1f4gen}
\end{equation}
where $$h^{2}(u(x)) x^{\prime}(u(x)) = \frac{1}{\mathfrak{s}(x)}$$ and $\mathfrak{s}(x)$ is the speed density given by \eqref{speed}.}
\end{rem}

\begin{rem} \label{Solf1f4}
\emph{According to Linetsky (see \cite{Linetsky1}, Lemma 1 and Lemma 2) it follows:
\begin{itemize}
\item[(i)] for NO boundary $e_{1}=0$ there exist the unique (up to the constant factor) non-trivial solution $\varphi(x, \lambda)$ of the Sturm-Liouville equation \eqref{sleq} which is for all $\lambda \geq 0$ square integrable with respect to the speed density $\mathfrak{m}(x)$ in the neighborhood of $e_{1}$, satisfies the appropriate boundary condition at $e_{1}$ and such that $\varphi(x, \lambda)$ and $\varphi^{\prime}(x, \lambda)$ are continuous in $x$ and $\lambda$ in $\langle 0, \infty \rangle^{2}$ and entire in $\lambda$ for each fixed $x \in \langle 0, \infty \rangle$. In our case such a solution is the solution $f_{1}(x, -\lambda)$ given by \eqref{f1l}, under the condition $\alpha > 2$ which ensures its square integrability in the neighborhood of the boundary $0$. This solution satisfies the boundary condition $$\lim\limits_{x \to 0} \frac{f_{1}^{\prime}(x, -\lambda)}{\mathfrak{s}(x)} = 0.$$
\item[(ii)] for O/NO boundary $e_{2}=\infty$ there exist the unique (up to the constant factor) non-trivial solution $\psi(x, \lambda)$ of the Sturm-Liouville equation \eqref{sleq} which is for all $\lambda < \Lambda$ square integrable with respect to the speed density $\mathfrak{m}(x)$ in the neighborhood of $e_{2}$, satisfies the appropriate boundary condition at $e_{2}$ and such that $\psi(x, \lambda)$ and $\psi^{\prime}(x, \lambda)$ are continuous in $x$ and $\lambda$ in $\langle 0, \infty \rangle \times \langle 0, \Lambda \rangle$ and analytic in $\lambda < \Lambda$ for each fixed $x \in \langle 0, \infty \rangle$. In our case such a solution is the solution $f_{4}(x, -\lambda)$ given by \eqref{f34l}, which satisfies the boundary condition $$\lim\limits_{x \to \infty} \frac{f_{4}^{\prime}(x, -\lambda)}{\mathfrak{s}(x)} = 0.$$
\end{itemize}}
\end{rem}

\end{document}